\documentclass[12pt]{amsart}
\usepackage{amsmath, amsthm, amssymb,cite,enumerate,slashed}
\usepackage{fullpage}
\usepackage{color}
\usepackage{hyperref}
\usepackage[english]{babel}
\usepackage{framed}

\allowdisplaybreaks

\newtheorem{thm}{Theorem}
\newtheorem{lem}{Lemma}[section]

\newtheorem{prop}[lem]{Proposition}
\newtheorem{cor}[lem]{Corollary}

\newtheorem{rmkk}[lem]{Remark}
\newtheorem{rmk}{Remark}[lem]

\newtheorem{defn}[lem]{Definition}
\numberwithin{equation}{section}

\newcommand{\R}{\mathbb{R}}

\newcommand{\Cp}{\mathbb{C}}

\newcommand{\Z}{\mathbb{Z}}

\newcommand{\I}{\mathcal{I}}

\newcommand{\F}{\mathcal{F}}
\newcommand{\eps}{\varepsilon}
\newcommand{\norm}[1]{\left\lVert#1\right\rVert}

\newcommand{\mcl}[1]{\mathcal{ #1 }}
\newcommand{\fm}[1]{\begin{align*} #1 \end{align*}}
\newcommand{\eq}[1]{\begin{equation}\begin{aligned} #1 \end{aligned}\end{equation}}
\newcommand{\lra}[1]{\langle #1 \rangle}
\newcommand{\abs}[1]{\left| #1 \right|}
\newcommand{\kh}[1]{\left( #1 \right)}

\newcommand{\wt}[1]{\widetilde{ #1 }}
\newcommand{\wh}[1]{\widehat{ #1 }}

\DeclareMathOperator{\Real}{Re}
\DeclareMathOperator{\Imag}{Im}
\DeclareMathOperator{\supp}{supp}

\DeclareMathOperator{\sgn}{sgn}

\DeclareMathOperator{\sech}{sech}
\allowdisplaybreaks[2]

\begin{document}

\title{Asymptotic stability of the sine-Gordon kinks under perturbations in  weighted Sobolev norms}

\author{Herbert Koch}
\address{Mathematical Institute, Hausdorff Center for Mathematics, University of Bonn}
\email{koch@math.uni-bonn.de}
\author{Dongxiao Yu}
\address{Department of Mathematics, Vanderbilt University}
\email{dongxiao.yu@vanderbilt.edu}

\begin{abstract}

We study the asymptotic stability of the sine-Gordon kinks under small perturbations in weighted Sobolev norms. Our main tool is the B\"acklund transform which reduces the study of the asymptotic stability of the kinks to the study of the asymptotic decay of solutions near zero. Our results consist of two parts. First, we prove an asymptotic stability result similar to the local results in \cite{MR4616683,chen2020long}. Our assumptions are the same as those in the local result in \cite{chen2020long}.  In its proof, we apply a result obtained by the inverse scattering method on the local decay of the solutions with sufficiently small and localized initial data. Moreover, we derive an asymptotic formula for the perturbations, i.e.\ the difference between solutions and kinks. This result is similar to that in  \cite{luhrmann2021asymptotic} and the full asymptotic stability result in \cite{chen2020long}. In its proof, we apply a result obtained by the method of testing by wave packets on the pointwise decay of the solutions with small and localized data. 
\end{abstract}

\maketitle
\tableofcontents\addtocontents{toc}{\protect\setcounter{tocdepth}{1}}

\section{Introduction}

This paper is devoted to the study of the one-dimensional sine-Gordon equation
\eq{\label{sg} f_{tt}-f_{xx}+\sin f=0}
where $f=f(t,x)$ is an $\R$-valued function defined on $\R\times\R$.  The sine-Gordon equation was originally introduced by  Bour \cite{ebour1862} in 1862 in the study of surfaces with constant negative curvature. It has been widely studied since it appears in several different fields such as differential geometry, relativistic field theory, quantum field theory, etc. For more background, we refer our readers to the monographs \cite{MR3642602,MR591458,MR2068924,MR3290305}. We also notice that the sine-Gordon equation is a completely integrable system \cite{MR406175,MR0406174,MR591458}. As a result, one can solve the initial value problem of \eqref{sg} by following the inverse scattering algorithm introduced by Zakharov and Shabat \cite{MR0406174}; we also refer to \cite{MR406175,takhtadzhyan1974exact,MR450817,MR0389050}.

The sine-Gordon equation admits an infinite number of conserved quantities, including the energy
\eq{\label{energy}E=\int_\R\frac{1}{2}(f_t^2+f_x^2)+1-\cos f\ dx,}
and the momentum
\eq{P=\int_\R \frac{1}{2}f_tf_x\ dx.}
Using the energy conservation and a fixed point argument, one can show that \eqref{sg} is globally well-posed in the natural energy space $H^1_{\sin}\times L^2(\R)$ where 
\eq{H^1_{\sin}(\R)=\{\phi\in \dot{H}^1(\R):\ \sin(\phi/2)\in L^2(\R)\}.}
We refer to \cite{MR3955109,MR1674843,MR3906859} for more details.

In this paper, we study the stability of a special type of global solutions to \eqref{sg} called the \emph{kinks}. For each fixed velocity $\beta\in(-1,1)$ and center $x_0\in\R$, a kink is of the form
\eq{\label{kink} Q(t,x;\beta,x_0)=4\arctan(e^{\gamma(x-\beta t-x_0)}),\qquad \gamma=(1-\beta^2)^{-1/2}.}
One notices that the kinks are not localized in space and that they connect two constant solutions $0$ and $2\pi$. 
Since the equation \eqref{sg} admits symmetries in translations $f(t,x)\mapsto f(t-t_0,x -x_0)$ and in Lorentz boosts \fm{f(t,x)\mapsto f_\beta(t,x)=f(\gamma(t-\beta x),\gamma(x-\beta t)),\qquad \forall \beta\in(-1,1),\ \gamma=(1-\beta^2)^{-1/2},}
we can obtain a general kink by applying translations and Lorentz boosts to the static kink centered at $0$ (i.e.\ $\beta=x_0=0$).

It has been proved in Henry-Perez-Wreszinski \cite{MR678151} that the kinks are \emph{orbitally stable} under small perturbations in $H^1\times L^2(\R)$. Precisely, if $f$ is a global solution to \eqref{sg} with
\fm{\norm{(f,f_t)|_{t=0}-(Q,Q_t)(0,\cdot;\beta,x_0)}_{H^1\times L^2}\leq \eps}
for a sufficiently small $0<\eps\ll1$, then for some path of center $x(t)\in\R$ we have
\fm{\norm{(f,f_t)(t)-(Q,Q_t)(t,\cdot;\beta,x(t))}_{H^1\times L^2}\lesssim \eps,\qquad \forall t\in\R.}
Here by $Q_t(t,x;\beta,x(t))$ we mean first taking the time derivative of \eqref{kink}  and then setting $x_0=x(t)$. Later, Hoffman and Wayne  \cite{MR3059166}  made use of the \emph{B\"acklund transform} and presented a different proof of the orbital stability result for the kinks. The B\"acklund transform, which will be defined in Section \ref{sec1ci},  is closely related to the complete integrability of the sine-Gordon equation. It is very useful in the study of the stability of kinks since it maps a zero solution to a kink. 

Now we consider the \emph{asymptotic stability} of the kinks. Roughly speaking, assuming that $(f,f_t)|_{t=0}$ is sufficiently close to a kink $(Q,Q_t)(0,\cdot;\beta_0,x_0)$ in a certain function space, we seek to find a final velocity $\beta$ and a path of center $x(t)$ so that \fm{f(t)\to Q(t,\cdot;\beta,x(t))\qquad\text{as $|t|\to\infty$ in some sense.}}
The convergence here could be in $L^\infty(\R)$ or  in a local  form: similar to the results in \cite{MR4616683}, we may obtain  
\eq{\label{s1weakasystaconv}\norm{(f,f_t)(t)- (Q,Q_t)(t,\cdot;\beta,x(t))}_{H^1(I)\times 
L^2(I)}\to 0\qquad \text{as }|t|\to\infty}
for all fixed bounded intervals $I$. Note that the kinks are not asymptotically stable under perturbations in $H^1\times L^2$ because of the existence of the \emph{wobbling kinks}; see \cite[Remark 1.3]{MR3630087} and \cite[Proposition 1.11]{chen2020long}. Some assumptions on the symmetry or decay of the perturbations of the initial data (i.e.\ the difference $(f,f_t)|_{t=0}-(Q,Q_t)(0,\cdot;\beta_0,x_0)$) are thus necessary. For example, Alejo, Mu$\tilde{\rm n}$oz, and Palacios \cite{MR4616683} constructed a smooth manifold of initial data close to the static kink centered at the origin, for which there is asymptotic stability in the energy space. More precisely, they showed that if the path $x(t)$ remains bounded for all $t\in\R$, then \eqref{s1weakasystaconv} holds; if there exists a sequence $t_n\to\pm\infty$ such that the path of center satisfies $|x(t_n)|\to\infty$, then \eqref{s1weakasystaconv} holds with $t$ replaced by $t_n$. Their tool is also the B\"acklund transform. In the paper, they assumed that the initial data has zero momentum and spatial symmetry: $f|_{t=0}-4\arctan e^x $ is odd and  $f_t|_{t=0}$ is even. Moreover, Chen, Liu, and Lu \cite{chen2020long} proved local asymptotic stability of the sine-Gordon kinks under small perturbations in the weighted Sobolev space $(H^{1,s}\times L^{2,s})(\R)$ with $s>1/2$ where
\fm{H^{m,s}(\R)=\{\phi\in L^1_{\rm loc}(\R):\ \lra{\cdot}^s\phi\in H^m(\R)\},\ m\geq 0;\quad L^{2,s}(\R)=H^{0,s}(\R).}
Here $\lra{x}:=\sqrt{1+x^2}$ for all $x\in\R$.
The authors applied the inverse scattering method, and no assumption on spatial symmetry is needed. Their result is also in the sense of \eqref{s1weakasystaconv} but a larger collection of intervals are considered: for every fixed $\mathfrak{v}\in\R$ and $\mathfrak{L}>0$, they proved that the limit \eqref{s1weakasystaconv} holds with $I$ replaced by 
\fm{I_{\mathfrak{v},\mathfrak{L}}(t)=\{y\in\R:\ |y-\mathfrak{v}t|\leq \mathfrak{L}\},\qquad  t\in\R}
as $t$ tends to infinity. They also showed that the path of center $x(t)$ remains bounded for all time.
In the same paper, Chen, Liu, and Lu also proved full asymptotic stability of the sine-Gordon kinks under perturbations in $(H^{2,s}\times H^{1,s})(\R)$ with $s>1/2$. They were able to derive the asymptotics of the difference $f(t,x)-Q(t,x;\beta,x(t))$ and obtain an $L^\infty$-type asymptotic stability result.
Besides, L\"uhrmann and Schlag \cite{luhrmann2021asymptotic} established an $L^\infty$-type asymptotic stability of the static kink centered at the origin under sufficiently small odd perturbations in $(H^{3,1}\times H^{2,1})(\R)$. Their approach does not rely on the complete integrability of \eqref{sg} and can be applied to some nonintegrable models. They also presented the asymptotics of the difference $f(t,x)-4\arctan e^x$. 

In this paper, we seek to prove the asymptotic stability of the sine-Gordon kinks under perturbations in weighted Sobolev norms via the B\"acklund transform. Our results consist of two parts.
\begin{enumerate}[(I)]
    \item When the perturbations are sufficiently small in $H^{1,1/2+}\times L^{2,1/2+}$, the difference $f(t)-Q(t,\cdot;\beta,x(t))$ converges to $0$ in $L^\infty$ and that the first derivatives of this difference converges to $0$ in $L^2+L^\infty$ (and thus in a local sense).
    \item When the perturbations are sufficiently small in $H^{5/2+,1}\times H^{3/2+,1}$, the $L^\infty$ norms of the difference $f(t)-Q(t,\cdot;\beta,x(t))$ and its first derivatives are bounded by $\eps \lra{t}^{-1/2}$. Moreover, we derive the asymptotics of this difference and its first derivatives.
\end{enumerate}
The result (I) is similar to those in \cite{MR4616683,chen2020long} and the result (II) is similar to those in \cite{chen2020long,luhrmann2021asymptotic}. Similar to \cite{chen2020long}, no assumption on the spatial symmetry of perturbations is required here. Detailed discussions on the similarities and differences among these results will be given in Remarks \ref{rmk:compare1} and \ref{rmk:compare2} after Theorem \ref{thmasysta}.

\subsection{Special explicit solutions to \eqref{sg}}
The sine-Gordon equation admits a large family of nontrivial explicit solutions. In addition to kinks, we notice that $-Q(t,x;\beta,x_0)$, where $Q(t,x;\beta,x_0)$ is given by \eqref{kink}, is also a global solution to \eqref{sg}. Such a solution is called an \emph{antikink}. Kinks and antikinks are not localized in space, and they are two simplest examples of the \emph{topological solitons}, i.e.\ solutions that are homotopically distinct from the zero solution. We refer to \cite{MR2068924} for an introduction to the topological solitons.

In addition to kinks and antikinks, the sine-Gordon equation admits another family of explicit solutions called \emph{breathers}:
\eq{\label{s1breather}B_v(t,x;\beta,x_1,x_2)=4\arctan\kh{\frac{\beta\cos(\alpha(t-vx-x_1))}{\alpha\cosh(\beta(x-vt-x_2))}}.}
Here we have fixed parameters $x_1,x_2\in\R$, $v\in(-1,1)$, $\gamma=(1-v^2)^{-1/2}$, $\beta\in(0,\gamma)$ and $\alpha=\sqrt{\gamma^2-\beta^2}$.
Unlike kinks and antikinks, the breathers are spatially localized. However, for each fixed $y_0\in\R$, we notice that $B_v(t,vt+y_0;\beta,x_1,x_2)$ is a periodic function of $t$. As a result, a breather does not decay as time goes to infinity. 

The kinks, antikinks, and breathers can be combined into more complicated structures. For example, by superposing a  kink and a breather with the same speed, we obtain a \emph{wobbling kink}, $W_\beta$, for each $\beta\in(-1,1)$:
\fm{W_\beta(t,x)&=4\arg(U_\beta+iV_\beta),\\U_\beta&=\cosh(\beta x)+\beta\sinh(\beta x)-\beta e^x\cos(t\sqrt{1-\beta^2}),\\
V_\beta&=e^x(\cosh(\beta x)-\beta\sinh(\beta x)-\beta e^{-x}\cos(t\sqrt{1-\beta^2})).}
The wobbling kinks were discovered by Segur \cite{MR708660}, and are another type of explicit solutions to \eqref{sg}. We also recall from \cite{MR3630087,chen2020long} that the existence of the wobbling kinks is the reason why a kink is not asymptotically stable in $H^1\times L^2$.
For other more complicated structures, we refer to \cite{MR4616683,MR2695612,MR3955109,MR767915,chen2020long}. 

\subsection{Complete integrability}\label{sec1ci}
The sine-Gordon equation is completely integrable in the sense that it admits a Lax pair (see, e.g.,\ \cite{MR406175,MR2348643,MR2695612,MR1697487,chen2020long}). In fact, the sine-Gordon equation is the compatibility condition for the system
\eq{\label{laxpaireqn}\Psi_x=A\Psi,\qquad \Psi_t=B\Psi}
where
\eq{A&=\frac{-iz}{4}\begin{pmatrix}
1&0\\0&-1
\end{pmatrix}+\frac{i}{4z}\begin{pmatrix}
\cos f&\sin f\\\sin f&-\cos f
\end{pmatrix}+\frac{f_t+f_x}{4}\begin{pmatrix}
0&-1\\1&0
\end{pmatrix},}
\eq{B&=\frac{-iz}{4}\begin{pmatrix}
1&0\\0&-1
\end{pmatrix}-\frac{i}{4z}\begin{pmatrix}
\cos f&\sin f\\\sin f&-\cos f
\end{pmatrix}+\frac{f_t+f_x}{4}\begin{pmatrix}
0&-1\\1&0
\end{pmatrix}.} 
Here $z\in\R\setminus\{0\}$ and $\Psi=\Psi(t,x;z)\in\Cp^2$. By compatibility condition, we mean that there exist two unique solutions to \eqref{laxpaireqn} with initial data $(1,0)$ and $(0,1)$ respectively if and only if $f$ solves \eqref{sg}. We also refer our readers to \cite{1974TMP} for a different way to formulate the complete integrability of \eqref{sg}. There the Lax pair for \eqref{sg} was given in a more classical sense:  two $4\times 4$ matrices $L$ and $M$  were given such that $\partial_tL=[L,M]$.

The kinks, antikinks, and breathers introduced in the previous subsection are closely related to the spectral problem 
$\Psi_x=A\Psi$. The eigenvalues here are not in the usual sense, so we quickly recall the definition from \cite{MR406175}.
Under reasonable assumptions on $f$ (e.g. it was assumed in \cite[Section 2]{chen2020long}  that the initial data  belongs to $H^{2,s}_{\sin}\times H^{1,s}(\R)$ for some $s>1/2$), for all $z\in\Cp\setminus \{0\}$ with $\Imag z\geq 0$ (not just for $z\in\R\setminus\{0\}$), we can  find two unique solutions $\Psi_l,\Psi_r$ to $\Psi_x=A\Psi$ such that
\fm{\Psi_l\exp(\frac{i}{4}(z-z^{-1})x)\to \begin{pmatrix}
    1\\0
\end{pmatrix},\qquad\text{as }x\to-\infty;\\
\Psi_r\exp(-\frac{i}{4}(z-z^{-1})x)\to \begin{pmatrix}
    0\\1
\end{pmatrix},\qquad\text{as }x\to+\infty.}These two solutions $\Psi_l$ and $\Psi_r$ are called the \emph{Jost solutions}.
For $z\in\R\setminus \{0\}$, if $\Psi=(\psi_1,\psi_2)$ is a solution to $\Psi_x=A\Psi$, then so is $\Psi^*:=(\overline{\psi_2},-\overline{\psi_1})$. Since $\Psi_r^*$ and $\Psi_r$ are linearly independent (otherwise $\Psi_r=0$), we can write
\fm{\Psi_l=a(z)\Psi_r^*+b(z)\Psi_r,\qquad z\in\R\setminus \{0\}.}
The coefficient $a(z)$ can be analytically continued to the upper half plane $\Imag{z}>0$. The zeros of this analytic continuation in the upper plane are defined as the eigenvalues of the spectral problem $\Psi_x=A\Psi$.
Because of \eqref{laxpaireqn}, the eigenvalues are independent of time $t$.  As remarked in \cite{MR406175}, the eigenvalues must be either purely imaginary or arise as complex conjugate pairs: $z$ and $-\bar{z}$. In the latter case, neither $z$ nor $-\bar{z}$ is purely imaginary.

For example, if $f$ is a kink or an antikink with a velocity $\beta\in(-1,1)$, then  the spectral problem above has exactly one eigenvalue: $2i\sqrt{\frac{1-\beta}{1+\beta}}$. If $f$ is a breather, then the spectral problem has exactly two eigenvalues: $z$ and $-\bar{z}$. Neither $z$ nor $-\bar{z}$ is purely imaginary. If $f$ is a wobbling kink, then the spectral problem has three distinct eigenvalues: an imaginary number and a pair $z$ and $-\bar{z}$. Again, neither $z$ nor $-\bar{z}$ is purely imaginary. The values of these eigenvalues can be computed explicitly from the parameters in the definitions of breathers and wobbling kinks. We skip such formulas here for simplicity and refer to \cite{MR406175}.

Conversely, the eigenvalues and the long time dynamics of $\Psi_l,\Psi_r$ determine the long time dynamics of $f$. For example,  Cheng \cite{MR2695612} and Cheng, Venakides, and Zhou \cite{MR1697487} studied the long time asymptotics of the sine-Gordon
equation with soliton-free initial data, i.e.\ when there is no eigenvalue. Recently, Chen, Liu, and Lu \cite{chen2020long} studied the general case and proved a soliton resolution conjecture for the sine-Gordon equation. The methods used in these papers are the inverse scattering method and the nonlinear steepest descent method.

The B\"acklund transform, which is the main tool used in this paper, also results from the complete integrability of \eqref{sg}. For a fixed constant $a>0$ and two globally defined functions $\phi=\phi(t,x)$ and $f=f(t,x)$, we say that $f$ is the B\"acklund transform of $\phi$ by the parameter $a$ if 
\eq{\label{s1backeqn}\left\{
\begin{array}{l}
  \displaystyle f_x-\phi_t-\frac{1}{a}\sin(\frac{f+\phi}{2})-a\sin(\frac{f-\phi}{2})=0,\\[1em]
 \displaystyle
 f_t-\phi_x-\frac{1}{a}\sin(\frac{f+\phi}{2})+a\sin(\frac{f-\phi}{2})=0.
\end{array}
\right.}
Note that the kink with a velocity $\beta\in(-1,1)$ and an arbitrary center is the B\"acklund transform of the zero solution by the parameter $a=\sqrt{\frac{1+\beta}{1-\beta}}$. 
This transform has been widely used in the stability problems for \eqref{sg}; we refer to \cite{MR3059166,MR4616683,MR3955109} for some examples. We also refer to \cite{koch2024multisolitons,MR3059166,MR2920823,MR3180019,MR3353827} for the applications of the B\"acklund transform to other models.

We now explain how the B\"acklund transform is related to the Lax pair. We follow the discussions in \cite[Section 5.4]{MR1146435}. To apply the formulas in this reference directly, we make a change of variables $(q,p)=\frac{1}{4}(x-t,x+t)$. Now we  rewrite \eqref{sg} as
\fm{f_{pq}=4\sin(f)}
and have the Lax pair
\eq{\label{newlaxpair}\begin{pmatrix}
\psi\\\varphi
\end{pmatrix}_q=\begin{pmatrix}
if_q/2&\lambda \\
\lambda&-if_q/2
\end{pmatrix}\begin{pmatrix}
\psi\\\varphi
\end{pmatrix},\qquad \begin{pmatrix}
\psi\\\varphi
\end{pmatrix}_p=\begin{pmatrix}
0&\lambda^{-1}e^{if} \\
\lambda^{-1}e^{-if} &0
\end{pmatrix}\begin{pmatrix}
\psi\\\varphi
\end{pmatrix}.}
Here $\lambda$ is a fixed constant, and $\psi,\varphi,f$ are  functions of $(p,q)$ (or $(t,x)$).
If $(\psi_1,\phi_1,\lambda_1,f)$ and $(\psi,\phi,\lambda,f)$ are two special solutions to \eqref{newlaxpair}, then the transform
\eq{\label{darboux}\begin{pmatrix}
\wt{\psi}\\\wt{\varphi}\\
\wt{\lambda}\\\wt{f}
\end{pmatrix}=\begin{pmatrix}
\lambda \varphi-\lambda_1\varphi_1\psi_1^{-1}\psi\\\lambda \psi-\lambda_1\varphi_1^{-1}\psi_1\varphi\\\lambda\\f-2i\ln(\varphi_1\psi_1^{-1})
\end{pmatrix}}
generates a new solution  to \eqref{newlaxpair}. We can check that  $\wt{f}$ is a new solution to the sine-Gordon equation. The formula \eqref{darboux} is the \emph{Darboux transform}, and we refer our readers to \cite{MR1146435} for a general introduction. As a consequence of the Darboux transform \eqref{darboux},  $\wt{f}$ is the B\"acklund transform of $f$ by the parameter $-\lambda_1$.  To see this, we set $h=\varphi_1/\psi_1$. From \eqref{darboux} we have $e^{i(\wt{f}-f)/2}=h$ and thus
\fm{h-h^{-1}=2i\sin(\frac{\wt{f}-f}{2}),\qquad  e^{if}h-e^{-if}h^{-1}=2i\sin(\frac{\wt{f}+f}{2}).}
By \eqref{newlaxpair} and \eqref{darboux}, we have 
\fm{(\wt{f}+f)_q&=2f_q-2i h^{-1}h_q=2i\lambda_1(h-h^{-1})=-4\lambda_1\sin(\frac{\wt{f}-f}{2}),\\
(\wt{f}-f)_p&=-2i h^{-1}h_p=2i\lambda_1^{-1}(e^{if}h-e^{-if}h^{-1})=-4\lambda_1^{-1} \sin (\frac{\wt{f}+f}{2}).}
Also recall that $(\partial_t,\partial_x)=\frac{1}{4}(\partial_p-\partial_q,\partial_p+\partial_q)$, so 
\fm{\wt{f}_x-f_t&=\frac{1}{4}(\wt{f}_p-f_p+\wt{f}_q+f_q)=-\lambda_1^{-1}\sin(\frac{\wt{f}+f}{2})-\lambda_1\sin(\frac{\wt{f}-f}{2}),\\
\wt{f}_t-f_x&=\frac{1}{4}(\wt{f}_p-f_p-\wt{f}_q-f_q)=
-\lambda_1^{-1}\sin(\frac{\wt{f}+f}{2})+\lambda_1\sin(\frac{\wt{f}-f}{2}).}
These two formulas coincide with \eqref{s1backeqn} with $(\phi,f,a)$ replaced by $(f,\wt{f},-\lambda_1)$.

\subsection{Nonlinear Klein-Gordon equations in 1D}
There have also been several studies on the long time asymptotics for the sine-Gordon equation not relying on the Lax pair. As a result, many results listed here apply not only to the sine-Gordon equation but also to several nonintegrable wave-type models. In this subsection, we collect and summarize several results in this direction.

We start with a simple case when the initial data is small and localized. In this case, one expects that there is no soliton and that the global solution decays to zero as time goes to infinity. To see this, we make use of the Taylor expansion of $\sin(x)$ at $x=0$ and rewrite \eqref{sg} as
\fm{f_{tt}-f_{xx}+f=\frac{1}{6}f^3+O(|f|^5).}
When the solution is small, the nonlinear effect of the remainder term $O(|f|^5)$ is expected to be much weaker than the cubic term, so the leading order behavior of the solution to \eqref{sg} is given by the 1D cubic Klein-Gordon equation\fm{f_{tt}-f_{xx}+f=\beta_0f^3,\qquad \beta_0\in\R.}
The global existence and long time dynamics for this equation with small and localized initial data have been extensively studied; we refer to \cite{MR2457221,MR2572684,MR2134954,MR2188297,MR3494176}.
One could also study the global existence and long time dynamics for a general 1D nonlinear Klein-Gordon equation of the form
\fm{f_{tt}-f_{xx}+f=F(f,f_t,f_x,f_{tx},f_{xx}).}
Delort \cite{MR1833089,MR2245535} proved a global existence result for this equation when the nonlinearity satisfies a null condition. We also refer to \cite{MR3864873} for another result related to the null condition, and to  \cite{MR4693097,MR2561936,MR3801822,MR3050628} for a collection of results for specifically given nonlinearities $F$.

We move on to the stability of kinks. Suppose that $f$ solves \eqref{sg}. From a perturbative perspective, we consider the difference $u(t,x)=f(t,x)-4\arctan(e^x)$ of $f$ and the static kink.  Now $u$ can be viewed as a small solution to 
\fm{u_{tt}-u_{xx}+(1-2\sech(x)^2)u=-\sech(x)\tanh(x)u^2+(\frac{1}{6}-\frac{1}{3}\sech^2(x))u^3+\text{higher order}.}
Similarly, if we study the stability of the kink $Q=\tanh(x/\sqrt{2})$ of the $\phi^4$ model \fm{f_{tt}-f_{xx}=f-f^3}
by considering the difference $u=f-\tanh(x/\sqrt{2})$,  we obtain\fm{u_{tt}-u_{xx}+(2-3\sech(x/\sqrt{2}))u=-3\tanh(x/\sqrt{2})u^2-u^3.} These two examples motivate the study of the long time dynamics of the 1D Klein-Gordon equation with variable coefficient quadratic and cubic nonlinearities and with potential
\eq{\label{eq:1dnlkg_gen}u_{tt}-u_{xx}+(m^2+V(x))u=\alpha(x)u^2+\beta(x)u^3.}
Here $m>0$ is a fixed constant.
Again we are interested in the case when the initial data is small and localized. We hope to show that $u$ decays to zero and exhibits a modified scattering behavior as time goes to infinity.  In the case when $V\equiv 0$,  we  refer to \cite{MR3403074,MR3449234,MR4284529,MR4189725}. For general potentials, we refer to \cite{MR4451296,MR4565686,luhrmann2021asymptotic,MR4529848,MR3630087,MR4242134,MR4506085,MR4553942,MR4710855,kowalczyk2022kink} and references therein, among which several asymptotic stability results of kinks or other solitons for several 1D wave-type models (including the sine-Gordon equation and the $\phi^4$ model) have been established. 

We end this subsection with a remark on the difficulties in studying the long time dynamics of \eqref{eq:1dnlkg_gen}. There are obvious difficulties arising from the nonlinear terms, due to the slow dispersive decay for the 1D linear Klein-Gordon equation and the variable coefficients. Meanwhile, the linear operator $L=\partial_t^2-\partial_x^2+(m^2+V)$ could exhibit a threshold resonance (i.e.\ a resonance at the bottom of the continuous spectrum of $L$) and an internal mode (i.e.\ a positive gap eigenvalue below the continuous spectrum). For example, a threshold resonance occurs if we study the linearized equation around the kink $4\arctan(e^x)$ for the sine-Gordon equation ($m=1$, $V=-2\sech(x)^2$), or around the kink $\tanh(x/\sqrt{2})$ for the $\phi^4$ model ($m=\sqrt{2}$, $V=-3\sech(x/\sqrt{2})$). In fact, it even occurs for the 1D linear Klein-Gordon equation ($V=0$). Here we refer to \cite{MR4701407} (in particular \cite[Remark 1.7]{MR4701407}) where Chen and Pusateri introduced a general method to study the decay properties and asymptotics of small solutions when an odd or even zero energy resonance is present.   For the $\phi^4$ model, an internal mode also occurs. See \cite{MR4529848,kowalczyk2022kink,MR3630087}, etc. Such spectral properties of $L$  complicate the analysis of the long time dynamics of solutions to \eqref{eq:1dnlkg_gen}. We refer to, for example,  \cite[Section 1.1]{MR4565686} and \cite[Section 1.2]{luhrmann2021asymptotic} for detailed discussions on this topic. 

\subsection{Main results}\label{secmthm}
We now state the main results of this paper. In Theorem \ref{thmasysta}, we state the asymptotic stability results which will be proved in this paper. We ask our readers to compare this theorem with the results in \cite{chen2020long,luhrmann2021asymptotic,MR4616683}.  Then, to explain the role of the B\"acklund transform in our proofs, we present Theorem \ref{mthm}. Despite the similarity between the two main theorems, our goal here is to illustrate that the B\"acklund transform reduces the study of the asymptotic stability of kinks to the study of the asymptotic decay of solutions near zero.

Let us start with the statements of asymptotic stability. Because of the time reversal symmetry of the sine-Gordon equation, we only state our results for $t\geq 0$ for simplicity.  Define the energy space \fm{{\bf E}^0_{\sin}&=\{f\in C_tH^{1}_{\sin}(\R\times\R):\ \partial_tf\text{ exists and belongs to }C_tL^2(\R\times\R)\}.}

\thm{\label{thmasysta}
{\rm(a) (Asymptotic stability).} 
Fix $\beta_0\in(-1,1)$ and $x_0\in\R$. Consider the kink $K(t,x)=Q(t,x;\beta_0,x_0)$ defined by \eqref{kink}. Fix $s>1/2$ and let $f\in {\bf E}^0_{\sin}$ be a global solution to \eqref{sg} with 
\fm{\norm{(f,f_t)(0)-(K,K_t)(0)}_{(H^{1,s}\times L^{2,s})(\R)}\leq\eps\ll1.}Then, for sufficiently small  $\eps$, there exist $\beta\in(-1,1)$ with $|\beta-\beta_0|\lesssim\eps$ and a $C^1$ path $x=x(t)\in\R$ defined for all $t\geq 0$ satisfying the following properties:
\begin{enumerate}[\rm i)]
    \item $|x(0)-x_0|\lesssim\eps$ and  $\lim_{t\to\infty} |x'(t)|=0$; 
    \item we have \eq{\label{asystacenter}\int (f(t,x)-\wt{K}(t,x))\sech(\gamma(x-\beta t-x(t)))\ dx=0, \qquad \forall t\geq 0}
    where $\wt{K}(t,x)=Q(t,x;\beta,x(t))$ and $\gamma=(1-\beta^2)^{-1/2}$;
    \item   we have
    \eq{\label{thmasystakeyc0}\norm{(f-\wt{K},f_t-\wt{K}_0)(t)}_{(H^1\times L^2)(\R)}\lesssim\eps,} \eq{\label{thmasystakeyc}\lim_{t\to\infty}\kh{\norm{(f-\wt{K})(t)}_{L^\infty(\R)}+\norm{(f_t-\wt{K}_0,f_x-\wt{K}_x)(t)}_{(L^2+L^\infty)(\R)}}=0,}     \eq{\label{thmasystakeyc2}\lim_{t\to\infty}\norm{(f-\wt{K},f_t-\wt{K}_0,f_x-\wt{K}_x)(t)}_{L^2(\{x\in\R:\ |x|\geq t+R\})}=0,\qquad \text{for each fixed }R\in\R}
    where $\wt{K}(t,x)$ is defined in part ii) and $\wt{K}_0(t,x)=\kh{(\partial_tQ)(t,x;\beta,y)}|_{y=x(t)}$.
\end{enumerate}

{\rm(b)  (Asymptotics of the difference).}
If we assume that the initial data satisfies a stronger assumption: for a fixed $m>3/2$, we have \fm{\norm{(f,f_t)(0)-(K,K_t)(0)}_{(H^{m+1,1}\times H^{m,1})(\R)}\leq\eps\ll1,}
then in addition to the results listed in part (a), we also have
\eq{\norm{(f-\wt{K},f_t-\wt{K}_0,f_x-\wt{K}_x)(t)}_{L^\infty(\R)}\lesssim \eps \lra{t}^{-1/2},\qquad \forall t\geq 0.}
Besides, there exists a complex-valued function $W=W(\xi)$ defined on $\R$ and two sufficiently small constants $\delta,\kappa\in(0,1)$, such that  $|W(\xi)|\lesssim\eps\lra{\xi}^{-1-2\delta}$ and that
\eq{\label{thmasysta:asyc}&(f-\wt{K})(t,x)\\
&=\int_{\beta t+x(t)}^x \frac{\cosh(\gamma(y-\beta t-x(t)))}{\cosh(\gamma(x-\beta t-x(t)))}\mcl{A}_W(t,y;\beta,x(t))\ dy\\
&\quad-\frac{1}{2\cosh(\gamma(x-\beta t-x(t)))}\int_\R \sgn(y-\beta t-x(t))e^{-\gamma|y-\beta t-x(t)|}\mcl{A}_W(t,y;\beta,x(t))\ dy\\
&\quad+O_\kappa(\eps t^{-1/2-\kappa}),}
\eq{\label{thmasysta:asyc2}&(f_x-\wt{K}_x,f_t-\wt{K}_0)(t,x)\\
&=\kh{\mcl{A}_W(t,x;\beta,x(t))+\gamma\cos(\frac{\wt{K}}{2})(f-\wt{K}),\mcl{B}_W(t,x;\beta,x(t))-\beta\gamma\cos(\frac{\wt{K}}{2})(f-\wt{K})}(t,x)\\
&\quad+O_\kappa(\eps t^{-1/2-\kappa}).}
Here
\eq{\label{thmasysta:asyc3}&\mcl{A}_W(t,x;\beta,x(t))\\
&=t^{-1/2}1_{|x|<t}\cdot \Real\kh{(-i\lra{x/\rho}-\beta\gamma\cos(\frac{\wt{K}}{2})) W(x/\rho)\exp(-i\rho+\frac{i}{32\lra{x/\rho}}|W(x/\rho)|^2\ln t)},\\
&\mcl{B}_W(t,x;\beta,x(t))\\
&=t^{-1/2}1_{|x|<t}\cdot \Real\kh{(ix/\rho+\gamma\cos(\frac{\wt{K}}{2})) W(x/\rho)\exp(-i\rho+\frac{i}{32\lra{x/\rho}}|W(x/\rho)|^2\ln t)},}
and $\rho=\sqrt{t^2-x^2}$. Note that both $\mcl{A}_W$ and $\mcl{B}_W$ are supported in the light cone $\{|x|<t\}$, and that $|\mcl{A}_W|+|\mcl{B}_W|\lesssim \eps t^{-1/2}\lra{x/\rho}^{-2\delta}$.

{\rm(c) (Boundedness of the center).} If we further assume that for a fixed $m>3/2$,  \fm{\norm{(f,f_t)(0)-(K,K_t)(0)}_{(H^{m+1,2}\times H^{m,2})(\R)}\leq\eps\ll1,}
then in addition to the results listed in parts (a) and (b), we also have $|x(t)-x_0|\lesssim \eps$ for all $t\geq 0$.
}

\rmkk{\label{rmk:compare1}\rm Let us compare part (a) of Theorem \ref{thmasysta} with the local asymptotic stability results in \cite{MR4616683,chen2020long}. First, our assumptions on the initial data are different from those in \cite{MR4616683}. There the authors only required the data to be close to a kink in the usual energy space. However, they added some assumptions on the spatial symmetry of the perturbations, and their results hold for a subset of the data with the properties stated above. On the other hand, our assumptions on the initial data are the same as those in \cite{chen2020long}.

Next, we notice that  \eqref{thmasystakeyc} implies  \eqref{s1weakasystaconv} with $I$ replaced by $I_{\mathfrak{v},\mathfrak{L}}(t)$ for fixed $\mathfrak{v}\in\R$ and $\mathfrak{L}>0$ (which then implies \eqref{s1weakasystaconv} in the original sense). To see this, we set
\fm{\mcl{Q}(t,x):=|(f-\wt{K})(t,x)|+|(f_x-\wt{K}_x)(t,x)|+|(f_t-\wt{K}_0)(t,x)|.}
For an arbitrary weight function $w\in (L^1\cap L^\infty)(\R)$ with $w\geq 0$, we apply H\"older's inequality to obtain
\fm{\int_\R w\mcl{Q}^2 \ dx&\lesssim \norm{w}_{L^1\cap L^\infty}\norm{\mcl{Q}^2}_{L^1+L^\infty}\lesssim\norm{w}_{L^1\cap L^\infty}\norm{\mcl{Q}}^2_{L^2+L^\infty} .}
The implicit constants are independent of $w$. Thus, for all $M>0$, we obtain a uniform convergence
\eq{\label{thmasystakeyc:cor}\lim_{t\to\infty}\sup_{\norm{w}_{L^1\cap L^\infty}\leq M}\int_\R w(x)\mcl{Q}(t,x)^2 \ dx=0.}
To recover the local asymptotic stability result in \cite{chen2020long}, we simply apply~\eqref{thmasystakeyc:cor} with $w$ replaced by $1_{I_{\mathfrak{v,L}}(t)}$.

On the other hand, we remark that part (a) is weaker than the result in \cite{chen2020long} in view of the choice of center $x(t)$. We will discuss this difference later in Remark \ref{rmkcenterdiff}.}

\rmkk{\label{rmk:compare2}\rm We also compare part (b) of Theorem \ref{thmasysta} with the full asymptotic stability result in \cite{chen2020long} and the $L^\infty$-type asymptotic stability result in \cite{luhrmann2021asymptotic}. First, our assumptions on the initial data lie between those in \cite{chen2020long} and those in \cite{luhrmann2021asymptotic}. In \cite{chen2020long}, it was assumed that the perturbations belong to $H^{2,1/2+}\times H^{1,1/2+}$. In part (b) of Theorem \ref{thmasysta}, we assume that the perturbations belong to $H^{5/2+,1}\times H^{3/2+,1}$. In \cite{luhrmann2021asymptotic}, it was assumed that the perturbations are odd and belong to $H^{3,1}\times H^{2,1}$.

Second, the asymptotics derived in this paper are different from those in the full asymptotic stability result in \cite{chen2020long}. For example, our asymptotic formula  \eqref{thmasysta:asyc} involves an integral of $\mcl{A}_W$. In \cite[Corollary 10.3]{chen2020long}, the asymptotics of $f-\wt{K}$ are implicitly given by the radiation terms $f_{r,\cos}$ and $f_{r,\sin}$. These radiation terms can be computed using the scattering data (see, e.g., \cite[(8.10) and (8.12)]{chen2020long}), and those formulas look very different from ours. On the other hand,  our asymptotic formula  \eqref{thmasysta:asyc} is similar to that in \cite{luhrmann2021asymptotic}. If one considers odd perturbations, then we have $\wt{K}=4\arctan e^x$, $\beta=0$, and $x(t)\equiv 0$ in Theorem \ref{thmasysta}. Moreover,  the function $W$ obtained in part (b) is an even function. In summary, we reduce \eqref{thmasysta:asyc} to
\fm{f(t,x)-4\arctan e^x&=\int_0^x\frac{\cosh(y)}{\cosh(x)} \mcl{A}_W(t,y;0,0)\ dy+O(\eps t^{-1/2-}).}
This formula is essentially the same as \cite[(1.12)]{luhrmann2021asymptotic}. It is not a coincidence. In fact, the authors of \cite{luhrmann2021asymptotic} showed that
\eq{\label{s1wfkink}w=(-\partial_x-\tanh(x))(f-4\arctan e^x)}
solves a 1D nonlinear Klein-Gordon equation. Since $f(t)-4\arctan e^x$ is expected to be an odd function of $x$, by \eqref{s1wfkink} one has
\eq{\label{s1wfkink2}f(t,x)-4\arctan e^x=\int_0^x\frac{\cosh(y)}{\cosh(x)}w(t,y)\ dy.}
After deriving an asymptotic for $w$, they applied \eqref{s1wfkink2} and obtained an asymptotic for $f-4\arctan e^x$. Meanwhile, by the B\"acklund transform \eqref{s1backeqn} (with $a=1$), we have
\fm{\phi_t&=(\partial_x+\tanh(x))(f-\arctan e^x)+O(|f-4\arctan e^x|^2+|\phi|^2).}
This estimate will be derived in Section \ref{s3pfmthmiiv}; see the first equation in \eqref{s3varphieqnpre} (with $\beta=0$, $\gamma=1$, and $Q=4\arctan e^x$). As a result, we expect that $w\approx -\phi_t$. To derive the asymptotic for $f-4\arctan e^x$, we derive an asymptotic for $\phi_t$ and also apply \eqref{s1wfkink2}.

Besides, in Theorem \ref{thmasysta}, we present the asymptotics not only for $f-\wt{K}$ but also for $(f_x-\wt{K}_x,f_t-\wt{K}_0)$. The asymptotics for the difference of the first derivatives seem to be new and do not appear in \cite{chen2020long,luhrmann2021asymptotic}.

We finally remark that part (a) of Theorem \ref{thmasysta} contains an $L^\infty$-type asymptotic stability result. By \eqref{thmasystakeyc}, we have $(f-\wt{K})(t)\to0$ in $L^\infty$ as $t\to\infty$. In part (a), the assumptions are much weaker than those mentioned above, but we cannot obtain a quantitative decay rate for the convergence.}

\rmkk{\label{rmkcenterdiff}\rm 
In Theorem \ref{thmasysta}, the most important conclusion is the limit \eqref{thmasystakeyc}. One difficulty in its proof is that the choice of the center $x(t)$ that makes \eqref{thmasystakeyc} hold is not unique. For example,  if we have the limit \eqref{thmasystakeyc} for a path $x(t)$, then  \eqref{thmasystakeyc} also holds for any path of the form $x(t)+o_{t\to\infty}(1)$. Besides, the B\"acklund transform itself does not determine this center; we recall that the kink with a velocity $\beta\in(-1,1)$ and an arbitrary center is the B\"acklund transform of the zero solution by the parameter $a=\sqrt{\frac{1+\beta}{1-\beta}}$. 

In \cite{chen2020long}, the authors made use of information from the inverse scattering method to choose the center. Such information is not available in our proof. In \cite{luhrmann2021asymptotic}, since only odd perturbations are considered, the center of the moving kink would stay at the origin; i.e.\  $x(t)\equiv 0$. In \cite{MR4616683},  an extra orthogonality condition \cite[Lemma 8.4]{MR4616683} which determines the center was introduced.

In this paper, we make use of an orthogonality condition  \eqref{asystacenter} to choose the center.  This orthogonality condition seems different from that in \cite[Lemma 8.4]{MR4616683} since there the condition also involves the time derivative of $f$. However,  in \cite[Theorem 6.1]{MR4616683}, the authors only considered the data with zero momentum, in which case their orthogonal condition coincides with ours. From the computations in this paper, we also find it more convenient to not include the time derivative $f_t$ in \eqref{asystacenter}. We will explain how this condition is used in the next subsection. Unfortunately, our current method does not seem to allow us to fully recover all the bounds proved in \cite{chen2020long}, where the authors proved that $|x(t)-x_0|\lesssim\eps$ for all time under the assumptions in part (a) of Theorem~\ref{thmasysta}. It seems impossible to get this bound by only using the B\"acklund transform unless we put a stronger assumption on the initial data. This is reflected in part (c) in Theorem \ref{thmasysta}.}\\

\rm

The asymptotic stability of the sine-Gordon kinks follows from Theorem \ref{mthm} below. Note that the meanings of the formulas \eqref{thmasysta:asyc}-\eqref{thmasysta:asyc3} in part (b) of Theorem \ref{thmasysta} are also clear from Theorem \ref{mthm}.

\thm{\label{mthm}Fix  $\beta\in(-1,1)$ and $x_0\in \R$. Consider the kink $K(t,x)=Q(t,x;\beta,x_0)$ defined by \eqref{kink}. Let $f\in{\bf E}_{\sin}^0$ be a global solution to the sine-Gordon equation \eqref{sg} with
\fm{\norm{(f,f_t)(0)-(K,K_t)(0)}_{(H^1\times L^2)(\R)}\leq\eps\ll1.}
Let $\phi\in{\bf E}_{\sin}^0$ be a global solution to \eqref{sg} such that 
\fm{\norm{(\phi,\phi_t)(0)}_{(H^1\times L^2)(\R)}\leq\eps\ll1}
and that $f$ is the B\"acklund transform of $\phi$ by the parameter $a=\sqrt{\frac{1+\beta}{1-\beta}}$. Our conclusion is that, whenever $\eps$ is sufficiently small, there  exists a $C^1$ path $x=x(t)\in\R$ defined for all $t\geq 0$ satisfying the following properties:
\begin{enumerate}[\rm i)]
    \item $|x(0)-x_0|\lesssim\eps$ and  $|x'(t)|\lesssim \norm{(\phi,\phi_t,\phi_x)(t)}_{(L^2+L^\infty)(\R)}\lesssim\eps$; 
    \item we have \eq{\label{mthmc1}\int (f(t,x)-\wt{K}(t,x))\sech(\gamma(x-\beta t-x(t)))\ dx=0, \qquad \forall t\geq 0}
    where $\wt{K}(t,x)=Q(t,x;\beta,x(t))$ and $\gamma=(1-\beta^2)^{-1/2}$;
    \item for all $t\geq 0$, we have 
    \eq{\label{mthmc11}&\norm{(f-\wt{K},f_t-\wt{K}_0)(t)}_{(H^1\times L^2)(\R)}+\norm{(\phi(t),\phi_t(t))}_{(H^1\times L^2)(\R)}\lesssim\eps,}
    \eq{\label{mthmc2}&\norm{(f-\wt{K})(t)}_{L^\infty(\R)}+\norm{(f_t-\wt{K}_0,f_x-\wt{K}_x)(t)}_{(L^2+L^\infty)(\R)}\lesssim\norm{(\phi,\phi_t,\phi_x)(t)}_{(L^2+L^\infty)(\R)},}
    \eq{\label{mthmc3}&\norm{(f-\wt{K},f_t-\wt{K}_0,f_x-\wt{K}_x)(t)}_{L^\infty(\R)}\lesssim \norm{(\phi,\phi_t,\phi_x)(t)}_{L^\infty(\R)},} 
    where $\wt{K}(t,x)$ is defined in part ii) and $\wt{K}_0(t,x)=\kh{(\partial_tQ)(t,x;\beta,y)}|_{y=x(t)}$; in \eqref{mthmc3} we also need to assume that the right hand side is finite and sufficiently small;
    \item for all $t\geq 0$ and $x\in\R$, we have 
    \eq{\label{mthmc4}&(f-\wt{K})(t,x)\\
    &=\Phi(t,x;\beta,x(t))-\frac{\gamma}{2\cosh(\gamma(x-\beta t-x(t)))}\int_\R \Phi(t,y;\beta,x(t))\sech(\gamma(y-\beta t-x(t)))\ dy\\
    &\quad+O(\norm{(\phi,\phi_t)(t)}^2_{(L^2+L^\infty)(\R)}),\\
    &(f_x-\wt{K}_x,f_t-\wt{K}_0)(t,x)\\
    &=\kh{\phi_t-\beta\gamma\cos(\frac{\wt{K}}{2})\phi+\gamma\cos(\frac{\wt{K}}{2})(f-\wt{K}),\phi_x+\gamma\cos(\frac{\wt{K}}{2})\phi-\beta\gamma\cos(\frac{\wt{K}}{2})(f-\wt{K})}(t,x)\\
    &\quad+O(\norm{\phi(t)}^2_{L^\infty(\R)}+\norm{\phi_t(t)}^2_{(L^2+L^\infty)(\R)})}
    where \fm{\Phi(t,x;\beta,x(t))&=\int_{\beta t+x(t)}^x\frac{\cosh(\gamma(y-\beta t-x(t)))}{\cosh(\gamma(x-\beta t-x(t)))}\cdot(\phi_t-\beta\gamma\cos(\frac{\wt{K}}{2})\phi)(t,y)\ dy;}
    \item assuming that for some $s>1/2$,
    \eq{\label{mthmc21asu}\norm{(f,f_t)(0)-(K,K_t)(0)}_{(H^{1,s}\times L^{2,s})(\R)}<\infty,}
    we have for each $R\in\R$,  whenever $t\gtrsim_{R,s}1$ and $|x|\geq t+R$,
    \eq{\label{mthmc21} &|f(t,x)-\wt{K}(t,x)|+|f_x(t,x)-\wt{K}_x(t,x)|+|f_t(t,x)-\wt{K}_0(t,x)|\\
    &\lesssim_{R,s} \min\{t^{-1/4}\lra{|x|-t}^{-1/4},\lra{|x|-t}^{-s}\}+|\phi(t,x)|+|\phi_t(t,x)|+|\phi_x(t,x)|.}
\end{enumerate}
}

\rmkk{\rm The second equation in \eqref{mthmc4} is closely related to the B\"acklund transform \eqref{s1backeqn}. In fact, if we subtract 
\fm{\wt{K}_x-(\frac{1}{a}+a)\sin(\frac{\wt{K}}{2})=\wt{K}_0-(\frac{1}{a}-a)\sin(\frac{\wt{K}}{2})&=0
}
from the B\"acklund transform \eqref{s1backeqn} and notice that \fm{\sin(\frac{f\pm\phi}{2})-\sin(\frac{\wt{K}}{2})&=\frac{1}{2}\cos(\frac{\wt{K}}{2})(f-\wt{K}\pm\phi)+O(|f-\wt{K}|^2+|\phi|^2),}
we obtain the second equation in \eqref{mthmc4} with a different remainder.
See \eqref{s3varphieqnpre} in Section~\ref{s3pfmthmiiv}.

We also notice that the estimate for $f_x-\wt{K}_x$ in \eqref{mthmc4} can be viewed as an ordinary differential equation of $f-\wt{K}$. By solving this ODE and applying \eqref{mthmc1}, we can derive the estimate for $f-\wt{K}$ in \eqref{mthmc4}. }

\rmkk{\label{mthmrmkfst}\rm As mentioned previously, Theorem \ref{mthm} shows the advantages of the usage of the B\"acklund transform in the stability problem. To apply Theorem \ref{mthm} in the proof of Theorem \ref{thmasysta}, we start with a global solution $f$ to \eqref{sg} whose initial data is close to a kink and find a corresponding solution $\phi$ whose initial data is sufficiently close to zero. Once this $\phi$ is found, we can put all our focus on the asymptotic decays of $\phi$. Studying the asymptotic behaviors of $\phi$ is expected to be easier than studying $f$ directly because $\phi$ is close to $0$ while $f$ is close to a family of kinks.

For example, as long as we have\fm{\lim_{t\to\infty}\norm{(\phi,\phi_t,\phi_x)(t)}_{(L^2+L^\infty)(\R)}=0,}
we can apply \eqref{mthmc2} to conclude that $(f-\wt{K})(t)\to 0$ in $L^\infty$. Besides, as long as we derive the asymptotic for \eq{\label{s1::asyphi}(\phi_t-\beta\gamma\cos(\wt{K}/2)\phi,\phi_x+\gamma\cos(\wt{K}/2)\phi)} (which are given by the functions $\mcl{A}_W$ and $\mcl{B}_W$ in \eqref{thmasysta:asyc3}), we can apply \eqref{mthmc4} to derive the asymptotics \eqref{thmasysta:asyc} and \eqref{thmasysta:asyc2}. }

\rmkk{\rm
In part (b) of Theorem \ref{thmasysta}, we assume that the perturbations are in $H^{5/2+,1}\times H^{3/2+,1}$. We choose such a weighted Sobolev space because one can apply the method of testing by wave packets to derive the asymptotics for $\phi$ and $\phi_t$ provided that the data of $\phi$ is sufficiently small in $H^{5/2+,1}\times H^{3/2+,1}$. See Proposition \ref{thmwpt} below. Then, by \eqref{mthmc4}, we obtain \eqref{thmasysta:asyc}.

Meanwhile, because of the full asymptotic stability result in \cite{chen2020long} under small perturbations in $H^{2,1/2+}\times H^{1,1/2+}$, the assumptions of part (b) in Theorem~\ref{thmasysta} are likely to be not optimal. For simplicity, we do not seek to find the weakest assumptions in this paper, but we remark that the key steps in the proof are to obtain an asymptotic for \eqref{s1::asyphi} above and to apply \eqref{mthmc4}. For example, one may apply the inverse scattering method to prove Proposition~\ref{thmwpt} under weaker assumptions (e.g., for sufficiently small data in $H^{2,1/2+}\times H^{1,1/2+}$) and then apply \eqref{mthmc4} to finish the proof.}
\rm

\subsection{Ideas of the proofs}\label{s1ideapf}
In this subsection, we outline the proofs of both the main theorems. Before doing so, we first summarize several properties of the sine-Gordon kinks \eqref{kink} in Lemma \ref{lemQ}. Next, we study the B\"acklund transform in Section \ref{secbackt}. It is not clear from the definition itself whether the B\"acklund transform of a given function exists or not and whether a given function is the B\"acklund transform of another function or not. We answer these questions in Section \ref{secbackt} and obtain several estimates.

We now discuss the proof of Theorem \ref{mthm} in Section \ref{secpfmthm}. The structure of the proof here is similar to that in \cite[Section 6-11]{MR4616683} since the main tools in our paper and \cite{MR4616683} are both the B\"acklund transform. Meanwhile, the details of our proof are different from theirs. This is because we do not have assumptions on spatial symmetry and our final convergence results \eqref{thmasystakeyc} and \eqref{thmasystakeyc2} are different from theirs. In Section \ref{secpfmthm:ift}, we focus on the choice of the center $x(t)$. With the help of several results in Section \ref{secbackt} and orbital stability, we prove Proposition \ref{proporbsta} which shows the existence and uniqueness of the $C^1$ path $x(t)$ satisfying \eqref{mthmc1}.  Then, in Section \ref{s3pfmthmiiv}, we prove part i)-iv) of Theorem \ref{mthm}. The main idea of this proof is as follows. Given the $C^1$ path $x(t)$ obtained in Section \ref{secpfmthm:ift}, we define
\eq{\label{varphithetadefn}(\varphi,\theta)(t,x)&=(f(t,x)-Q(t,x;\beta,x(t)),f_t(t,x)-(\partial_tQ)(t,x;\beta,x(t))).}
We now seek to estimate $(\varphi,\varphi_x,\theta)$.
To achieve this goal, we derive the following equations for $(\varphi,\theta)$ from the B\"acklund transform \eqref{s1backeqn}:
\fm{\varphi_x&=\phi_t+\frac{1}{a}(\sin(\frac{\varphi+Q+\phi}{2})-\sin(\frac{Q}{2}))+a(\sin(\frac{\varphi+Q-\phi}{2})-\sin(\frac{Q}{2})),\\
\theta&=\phi_x+\frac{1}{a}(\sin(\frac{\varphi+Q+\phi}{2})-\sin(\frac{Q}{2}))-a(\sin(\frac{\varphi+Q-\phi}{2})-\sin(\frac{Q}{2})).}
With an abuse of notation, we write $Q(t,x;\beta,x(t))$ as $Q$. The first equation can be viewed as an ODE for $\varphi$. By solving it, we obtain 
\eq{\label{s1varphiid}\varphi(t,x)&=\frac{\varphi(t,\beta t+x(t))}{\cosh(\gamma(x-\beta t-x(t)))}+\int_{\beta t+x(t)}^x\frac{\cosh(\gamma(y-\beta t-x(t)))}{\cosh(\gamma(x-\beta t-x(t)))}F(t,y)\ dy}
where
\fm{F(t,x)&=\phi_t+\frac{1}{a}(\sin(\frac{\varphi+Q+\phi}{2})-\sin(\frac{Q}{2}))+a(\sin(\frac{\varphi+Q-\phi}{2})-\sin(\frac{Q}{2}))-\gamma\cos(\frac{Q}{2})\varphi\\
&=O(|\phi_t|+|\phi|+|\varphi|^2).}
Note that the orthogonality condition \eqref{mthmc1} has not been used yet. That is, we can start with an arbitrary $C^1$ path $x(t)$, define $(\varphi,\theta)$ by \eqref{varphithetadefn}, and 
derive \eqref{s1varphiid}. The identity \eqref{s1varphiid} itself does not give a good estimate for $\varphi(t,\beta t+x(t))$, because we only get a trivial identity when we set $x=\beta t+x(t)$ in \eqref{s1varphiid}. This explains why we introduce the orthogonality condition \eqref{mthmc1}. 
By plugging  \eqref{s1varphiid} into  \eqref{mthmc1}, we   express $\varphi(t,\beta t+x(t))$ explicitly in terms of  $F,\beta,x(t)$; see \eqref{s3varphiid}. This explicit formula along with the equations for $(\varphi_x,\theta)$ above allow us to prove estimates for $(\varphi,\varphi_x,\theta)$. For example, by noticing that \fm{F&=\phi_t-\beta \gamma \cos(\frac{Q}{2})\phi+O(|\varphi|^2+|\phi|^2)=\phi_t-\beta \gamma \cos(\frac{Q}{2})\phi+\text{lower order terms},}
we can derive the first equation in \eqref{mthmc4} from \eqref{s3varphiid}. In Section \ref{secpfmthm:final}, we finish the proof by showing \eqref{mthmc21}. The main tool used there is the finite speed of propagation and the energy-momentum tensor for the sine-Gordon equation. Part of our proof (e.g.\ the proof of Lemma \ref{lem4.7weiyang}) is essentially the same as that in \cite[Section 2]{MR4242133}.

We move on to the proof of Theorem \ref{thmasysta} in Section \ref{secpfasysta}. To apply Theorem \ref{mthm}, we first find the solution $\phi$ such that the data of $\phi$ is close to zero and that $f$ is the B\"acklund transform of $\phi$. These are proved in Corollary~\ref{corsolvebackthigh}. Then, as discussed in Remark \ref{mthmrmkfst}, we study the asymptotic decays of $\phi$; see Propositions \ref{thmlowdata} and \ref{thmwpt}.

To prove part (a), we apply Proposition \ref{thmlowdata}. The proof of this proposition can be found in Appendix \ref{seclow}. In its proof, we make use of the following two results obtained by applying the inverse scattering method. Note that both of them are stated in the soliton-free case. First, we make use of an approximation result (Proposition \ref{chen2020thm}) from \cite{chen2020long}. This theorem allows us to reduce the proof of Proposition \ref{thmlowdata} to the case when we have $C_c^\infty$ data.  Besides, we need the $t^{-1/2}$-pointwise decay from \cite{MR2695612,MR1697487}. Note that our assumption on the initial data seems too weak for us to apply nonintegrable methods such as energy methods. As a result, the inverse scattering method seems unavoidable here.  Fortunately, applying the inverse scattering method in the soliton-free case should be simpler than applying the same method to the case when a kink is present. 

To derive the asymptotics of the difference $f-\wt{K}$ in part (b), we apply Proposition \ref{thmwpt}. One can prove this proposition by the method of testing by wave packets which was first introduced by Ifrim-Tataru \cite{MR3382579}. We refer our readers to Section \ref{spfthmhrc} for more references on this method, and to Appendix \ref{secthmwpt} for a sketch of the proof.

Finally, to prove the boundedness of $x(t)$ in part (c), we introduce a new way to choose the center. Let us temporarily forget the previous choice of the center $x(t)$ and return to \eqref{s1varphiid}. To handle the term involving $\varphi(t,\beta t+x(t))$, previously we used the orthogonality condition \eqref{asystacenter}. There is however an easier way; we can simply choose $x(t)$ so that $\varphi(t,\beta t+x(t))=0$. This is equivalent to $f(t,\beta t+x(t))=\pi$ because by \eqref{varphithetadefn} one has
\fm{\varphi(t,\beta t+x(t))=f(t,\beta t+x(t))-Q(t,\beta t+x(t);\beta,x(t))=f(t,\beta t+x(t))-\pi.}
In the proof of part (c), we will show that there exists a $C^1$ path $\wt{x}(t)$ such that $f(t,\beta t+\wt{x}(t))=\pi$ and that $|\wt{x}(t)-x_0|\lesssim\eps$ for all $t\geq 0$. Moreover, if $x(t)$ is the center determined by \eqref{asystacenter}, we have $|x(t)-\wt{x}(t)|\lesssim\eps$, so both $x(t)$ and $\wt{x}(t)$ remain bounded for all the time. We even have $\lim_{t\to\infty}|x(t)-\wt{x}(t)|=0$. A proof will be given in Section \ref{secthmvariant}.

\rmkk{\rm 
Here we have some remarks on part (c). There we stated that $|x(t)-x_0|\lesssim\eps$. If we hope to choose the center by solving $f(t,\beta t+x(t))=\pi$, then we need $f\in C_{t,x}^1$ to apply the implicit function theorem. This is impossible if we only assume that  $(f-K,f_t-K_t)|_{t=0}$ is small in $H^{1,1/2+}\times L^{2,1/2+}$. As a result, we prefer to use \eqref{asystacenter} instead of $f(t,\beta t+x(t))=\pi$ to choose the center in Theorem \ref{thmasysta}. On the other hand, there is a variant of parts (b) and (c) of Theorem \ref{thmasysta}: the same conclusions hold if we replace \eqref{asystacenter} with $f(t,\beta t+x(t))=\pi$. The proof of this variant is essentially the same as the proof of Theorem \ref{thmasysta}.

Let us also explain why we need $H^{5/2+,2}\times H^{3/2+,2}$ in part (c). Later, we will see that the boundedness of $\wt{x}(t)$ and $x(t)$ results from the bound $|Z\phi(t,x)|\lesssim \eps t^{-c}$ for a constant $0<c<1/2$. Here $Z=x\partial_t+t\partial_x$ is the Lorentz boost. To obtain such a pointwise estimate, we need to estimate the Sobolev norm of $Z^2\phi$. At $t=0$, one notices that $Z^2\approx x^2\partial_t^2+l.o.t$ which is then related to $H^{*,2}$ norm.}
\rm

\subsection{Acknowledgement} The authors would like to thank Gong Chen and Daniel Tataru for many helpful discussions on this paper. The authors would also like to thank the anonymous reviewers for their valuable comments and suggestions on this paper. Both of the authors have been funded by the Deutsche
Forschungsgemeinschaft (DFG, German Research Foundation) through the Hausdorff Center
for Mathematics under Germany's Excellence Strategy - GZ 2047/1, Projekt-ID 390685813. This material is also based upon work supported by the National Science Foundation under Grant No.\ DMS-1928930 while both the authors were in residence at the Simons Laufer Mathematical Sciences Institute (formerly MSRI) in Berkeley, California, during the summer of 2023. Besides, the first author has also been funded by  Project ID 211504053 - SFB 1060. The second author has also been funded by a Simons Investigator grant of Daniel Tataru from the Simons Foundation and a VandyGRAF Fellowship from Vanderbilt University.

\subsection{Data availability statement}
Data sharing not applicable to this article as no datasets were generated or analysed during the current study.

\subsection{Declaration of competing interests} The authors have no competing interests to declare that are relevant to the content of this article. 

\section{Setup}
\subsection{Notations}\label{secnotation}
We use $C$ to denote universal positive constants, and the values of $C$ could vary from place to place. We write $A\lesssim B$ or $A=O(B)$ if $|A|\leq CB$ for some  $C>0$. We write $A\sim B$ if $A\lesssim B$ and $B\lesssim A$. We use $C_{v}$ or $\lesssim_v$ if we want to emphasize that the constant depends on a parameter $v$.

We say that a conclusion holds for all $c\ll1$ if there exists a small constant $c_0$, such that this conclusion holds for all $0\leq c<c_0$. Similarly for $C\gg1$. We also use $\ll_v,\gg_v$ to emphasize the dependence on a certain parameter $v$.

In addition to the $L^2$-based Sobolev space $H^m$ for each $m\in\R$, we make the following definitions. For each integer $k\geq 1$ and $p\in[1,\infty]$, define
\fm{W^{k,p}(\R)=\{f\in L^1_{\rm loc}(\R):\ \partial_x^jf\in L^p(\R),\ \forall j=0,1,\dots,k\}}
and
\fm{H^{k}_{\sin}(\R)=\{f\in L^1_{\rm loc}(\R):\ \partial_xf\in H^{k-1}(\R),\sin(f/2)\in L^2(\R) \}.}
For each  real numbers  $m,s\geq 0$, we define the weighted Sobolev space
\fm{H^{m,s}(\R)=\{f\in L^2(\R):\ \norm{\lra{x}^sf}_{H^m(\R)}<\infty\}.}
When $m$ is a nonnegative integer, $f\in H^{m,s}(\R)$ if and only if
\fm{\lra{x}^s\partial_x^jf\in L^2(\R),\qquad \text{for each }0\leq j\leq m.}

We also define
\fm{{\bf E}^{m,s}&= C_tH^{m+1,s}_x(\R\times\R)\cap C_t^1H^{m,s}_x(\R\times\R),\qquad m,s\in[0,\infty);\\
{\bf E}^k_{\sin}&=\{f\in C_tH^{k+1}_{\sin}(\R\times\R):\ \partial_tf\text{ exists and belongs to }C_tH^{k}_x(\R\times\R)\},\qquad k\in\Z_{\geq 0}.}
We remark that ${\bf E}^0_{\sin}$ is the usual energy space for the sine-Gordon equation with the energy~\eqref{energy}.

\subsection{Properties of the kinks}
In the next lemma, we state some important properties of the kinks which will be used in the rest of this paper.
\lem{\label{lemQ}Fix $\beta\in(-1,1)$ and $x_0\in\R$. Set $\gamma=(1-\beta^2)^{-1/2}$. Then, the kink $Q=Q(t,x;\beta,x_0)$ defined by \eqref{kink} satisfies the following properties.
\begin{enumerate}[\rm a)]
    \item We have \fm{\sin(\frac{Q}{2})=\sech(\gamma(x-\beta t-x_0)),\quad \cos(\frac{Q}{2})=-\tanh(\gamma(x-\beta t-x_0));}
    \fm{Q_x= 2 \gamma \sech(\gamma(x-\beta t-x_0)),\quad Q_t=-2\beta \gamma\sech(\gamma(x-\beta t-x_0)).}
    \item Fix  $t_1,t_2,x_1,x_2\in\R$. Fix an integer $k\geq 0$ and a real number $s\geq 0$. Then, for each $c\in(0,1)$, whenever $|\beta_1|,|\beta_2|\leq c$, we have
    \fm{&\norm{(Q,Q_t)(t_1,\cdot;\beta_1,x_1)-(Q,Q_t)(t_2,\cdot;\beta_2,x_2)}_{H^{k+1,s}\times H^{k,s}}\\
    &\lesssim_{k,s,c}(\lra{x_1+\beta_1t_1}^s+\lra{x_2+\beta_2t_2}^s)(\lra{t_1}|\beta_1-\beta_2|+|x_1-x_2|+|t_1-t_2|).}
    The implicit constant can be chosen to be uniform in $\beta_1,\beta_2,x_1,x_2,t_1,t_2$.
    \item Fix an integer $k\geq 0$ and a real number $s\geq 0$.  For each $w\in H^{k,s}$ and each $t\in\R$, we have
    \fm{\norm{\sin(\frac{w+Q(t)}{2})-\sin(\frac{Q(t)}{2})}_{H^{k,s}}&\lesssim_{k,s}(\norm{w}_{H^{k}}+1)^k\norm{w}_{H^{k,s}}}
    and therefore $\sin((w+Q(t))/2)\in H^{k,s}$. The implicit constant here is independent of $x_0$ and $t$.
\end{enumerate}
}
\begin{proof}
Part a) follows from direct computations. For part b), we first notice that
\fm{&Q(t_1,x;\beta_1,x_1)-Q(t_2,x;\beta_2,x_2)=Q(0,x;\beta_1,x_1+\beta_1t_1)-Q(0,x;\beta_2,x_2+\beta_2t_2)\\
&=\int_0^1((\beta_1-\beta_2)\partial_\beta Q+(x_1-x_2+\beta_1t_1-\beta_2t_2)\partial_{x_0} Q)(0,x;\wt{\beta}_\tau,\wt{x}_\tau)\ d\tau}
where \fm{(\wt{\beta}_\tau,\wt{x}_\tau)=(\tau\beta_1+(1-\tau)\beta_2,\tau(x_1+\beta_1t_1)+(1-\tau)(x_2+\beta_2t_2)),\qquad \tau\in [0,1].}
The same identity holds if we replace $Q$ with $\partial_t^j\partial_x^kQ$ with $j,k\geq 0$. Note that $|x_1-x_2+\beta_1t_1-\beta_2t_2|\leq |x_1-x_2|+|t_1||\beta_1-\beta_2|+|t_1-t_2|$. By the Minkowski inequality, we only need to estimate the $L^{2,s}_x$ norm (with a fixed $s\geq 0$) of $\partial_\beta \partial_t^j\partial_x^k Q(0,x;\wt{\beta}_\tau,\wt{x}_\tau)$  and $\partial_{x_0} \partial_t^j\partial_x^k Q(0,x;\wt{\beta}_\tau,\wt{x}_\tau)=-\partial_t^j\partial_x^{k+1}Q(0,x;\wt{\beta}_\tau,\wt{x}_\tau)$. 

Let us start with $\partial_t^j\partial_x^{k+1}Q(0,x;\wt{\beta}_\tau,\wt{x}_\tau)$. We have
\fm{\partial_t^j\partial_x^{k+1}Q(0,x;\wt{\beta}_\tau,\wt{x}_\tau)&=\partial_t^j\partial_x^{k}(2\wt{\gamma}_\tau \sech(\wt{\gamma}_\tau(x-\wt{\beta}_\tau t-\wt{x}_\tau)))|_{t=0}\\
&=2\wt{\gamma}_\tau^{j+k+1}(-\wt{\beta}_\tau)^j(\sech)^{(j+k)}(\wt{\gamma}_\tau(x-\wt{x}_\tau)).}
Here $\wt{\gamma}_\tau=(1-\wt{\beta}_\tau^2)^{-1/2}$. Recall that we assume $|\beta_1|,|\beta_2|\leq c$. It follows that $1\leq\gamma_1,\gamma_2,\wt{\gamma}_\tau\leq (1-c^2)^{-1/2}$. Thus, we have
\fm{&\norm{\partial_t^j\partial_x^{k+1}Q(0,x;\wt{\beta}_\tau,\wt{x}_\tau)}_{L^{2,s}_x}\lesssim_{c,j+k} \norm{\lra{x}^{s}(\sech)^{(j+k)}(\wt{\gamma}_\tau(x-\wt{x}_\tau))}_{L^2_x}\\
&\lesssim_{c,j+k,s} \norm{\lra{x-\wt{x}_\tau}^{s}(\sech)^{(j+k)}(\wt{\gamma}_\tau(x-\wt{x}_\tau))}_{L^2_x}+\norm{\lra{\wt{x}_\tau}^{s}(\sech)^{(j+k)}(\wt{\gamma}_\tau(x-\wt{x}_\tau))}_{L^2_x}\\
&\lesssim_{c,j+k,s}\lra{x_1+\beta_1t_1}^s+\lra{x_2+\beta_2t_2}^s. }
In the last two estimates, we use $\lra{x+y}^s\lesssim_s\lra{x}^s+\lra{y}^s$ for all $s\geq 0$ and $x,y\in\R$.

Now, to estimate $\partial_\beta \partial_t^j\partial_x^k Q(0,x;\wt{\beta}_\tau,\wt{x}_\tau)$, we notice that
\fm{&\partial_\beta \partial_t^j\partial_x^k Q(0,x;\wt{\beta}_\tau,\wt{x}_\tau)=\partial_t^j\partial_x^k(2\sech(\wt{\gamma}_\tau(x-\wt{\beta}_\tau t-\wt{x}_\tau))\cdot \partial_\beta(\gamma x-\beta \gamma t-\gamma \wt{x}_\tau))|_{t=0,\beta=\wt{\beta}_\tau}.}
Here $\partial_\beta(\gamma,\gamma\beta)|_{\beta=\wt{\beta}_\tau}=(\wt{\beta}_\tau\wt{\gamma}_\tau^3,\wt{\gamma}_\tau^3)$. We now apply the Leibniz rule. For those terms with a derivative falling on $\partial_\beta(\gamma x-\beta \gamma t-\gamma \wt{x}_\tau)$, we can estimate them by following the proof in the previous paragraph. The only remaining term is 
\fm{&[(\beta\gamma^3 (x-\wt{x}_\tau)- \gamma^3 t)\cdot \partial_t^j\partial_x^k(2\sech(\gamma_\tau(x-\wt{\beta}_\tau t-\wt{x}_\tau)))]|_{t=0,\beta=\wt{\beta}_\tau}\\
&=-2\wt{\gamma}_\tau^{j+k+3}(-\wt{\beta}_\tau)^{j+1}(x-\wt{x}_\tau)\cdot (\sech)^{(j+k)}(\wt{\gamma}_\tau(x-\wt{x}_\tau)).}
By applying $\lra{x+y}^s\lesssim_s\lra{x}^s+\lra{y}^s$ again, we conclude that the $L^{2,s}$ norm of this term is bounded by $\lra{x_1+\beta_1t_1}^s+\lra{x_2+\beta_2t_2}^s$.

For part c), we notice that
\fm{|\sin(\frac{w+Q}{2})-\sin(\frac{Q}{2})|\leq|w|.}
This finishes the proof when $k=0$. For general $k>0$, by Leibniz's  rule and the chain rule, we can write  
\fm{\partial_x^k(\sin(\frac{w+Q}{2})-\sin(\frac{Q}{2}))}
as a linear combination of terms of the form
\fm{&(\sin)^{(r)}(\frac{w+Q}{2})\prod_{l=1}^r\partial_x^{j_l}(\frac{w+Q}{2})-(\sin)^{(r)}(\frac{Q}{2})\prod_{l=1}^r\partial_x^{j_l}(\frac{Q}{2})\\
&=(\sin)^{(r)}(\frac{w+Q}{2})\prod_{l=1}^r\partial_x^{j_l}(\frac{w}{2})+[(\sin)^{(r)}(\frac{w+Q}{2})-(\sin)^{(r)}(\frac{Q}{2})]\prod_{l=1}^r\partial_x^{j_l}(\frac{Q}{2}).}
Here $r>0$, $\sum j_*=k$ and each $j_*>0$. Thus the absolute value of the left side is controlled by 
\fm{C\cdot \prod_{l=1}^r|\partial_x^{j_l}w|+|(\sin)^{(r)}(\frac{w+Q}{2})-(\sin)^{(r)}(\frac{Q}{2})|\prod_{l=1}^r|\partial_x^{j_l}Q|\leq C\prod_{l=1}^r|\partial_x^{j_l}w|+C|w|.}
Here we use $|\partial_x^jQ|\lesssim_j1$ everywhere for each $j\geq 1$ where the implicit constants are independent of $t,x_0$. Since at most one of $j_*$ exceeds $k/2$,   we conclude that
\fm{|\partial_x^k(\sin(\frac{w+Q}{2})-\sin(\frac{Q}{2}))|\lesssim \sum_{j\leq k}|\partial_x^jw|\cdot (\sum_{j\leq k/2}|\partial_x^jw|+1)^k.}
The constant here is uniform for all  $x_0,t\in\R$.
It follows that
\fm{\norm{\sin(\frac{w+Q}{2})-\sin(\frac{Q}{2})}_{H^{k,s}}&\lesssim_{k,s}(\sum_{j\leq k/2}\norm{\partial_x^jw}_{L^{\infty}}+1)^k\norm{w}_{H^{k,s}}\lesssim_{k,s}(\norm{w}_{H^{k}}+1)^k\norm{w}_{H^{k,s}}.}
In the last estimate, we use the inequality $k\geq \lfloor k/2\rfloor+1$ whenever $k\geq 1$ and the Sobolev embedding.
\end{proof}\rm 

\section{The B\"acklund transform}\label{secbackt}

In this section, we study the B\"acklund transform. Here we are mainly interested in solving the B\"acklund transform. In other words, we seek to answer the following two questions:
\begin{enumerate}[(Q1)]
    \item Given a global solution $\phi\in {\bf E}_{\sin}^0$ to \eqref{sg}, can we find a global solution $f\in {\bf E}_{\sin}^0$ such that $f$ is the B\"acklund transform of $\phi$?
    \item Given a global solution $f\in {\bf E}_{\sin}^0$ to \eqref{sg}, can we find a global solution $\phi\in {\bf E}_{\sin}^0$ such that $f$ is the B\"acklund transform of $\phi$?
\end{enumerate}
In this paper, we are only interested in the cases when the initial data of $\phi$ is close to zero and when the initial data of $f$ is close to a kink.

In Section \ref{secbackt:basic}, we define the B\"acklund transform and an associated functional $\F$. This $\F$ is also related to our choice of the center \eqref{asystacenter} and \eqref{mthmc1}. Next, we state a key proposition, Proposition \ref{propbackt}, which states that $\F=0$ induces a local $C^1$ diffeomorphism. For each fixed time $t$, velocity $\beta_0$, and center $x_0$, the local $C^1$ diffeomorphism maps between a neighborhood of the kink $(Q,Q_t)(t,\cdot;\beta_0,x_0)$ and a neighborhood of  the zero solution, the velocity $\beta_0$ and the center $x_0$. With the help of this diffeomorphism, we are able to answer (Q1) and (Q2); see Corollary \ref{corsolvebackt}.

Here is a brief summary of the answers. Given a solution $f$ whose data is close to a certain kink (in the sense that their difference is sufficiently small in $H^1\times L^2$), we are able to obtain a solution $\phi$ along with a velocity $\beta$, such that $f$ is the B\"acklund transform of $\phi$ by the parameter $\sqrt{\frac{1+\beta}{1-\beta}}$. Here the data of $\phi$ is close to zero in $H^1\times L^2$. This answers (Q2). In addition to $\phi$ and $\beta$, from this $f$ we also obtain a center $y_0$ in the sense that the orthogonal condition \eqref{mthmc1} holds at $t=0$ with $x(0)=y_0$. It turns out that the map $f\mapsto (\phi,\beta,y_0)$ is locally invertible. This answers (Q1).

Finally, since Theorem \ref{thmasysta} involves several weighted Sobolev norms, we also study the restrictions of the diffeomorphism above to these weighted Sobolev spaces; see Corollary \ref{corsolvebackthigh}.

\subsection{Definitions}\label{secbackt:basic}
We start with the definition of the B\"acklund transform.

\defn{\rm Fix a  constant $a>0$. Let $\phi=\phi(t,x)$ and $f=f(t,x)$ be two functions defined in $\R\times\R$. We say that  $f$  is the \emph{B\"acklund transform} of $\phi$ by the parameter $a$ if
\eq{\label{backeqn}\left\{
\begin{array}{l}
  \displaystyle f_x-\phi_t-\frac{1}{a}\sin(\frac{f+\phi}{2})-a\sin(\frac{f-\phi}{2})=0,\\[1em]
 \displaystyle
 f_t-\phi_x-\frac{1}{a}\sin(\frac{f+\phi}{2})+a\sin(\frac{f-\phi}{2})=0.
\end{array}
\right.}
}
\rm

We now define a functional which is closely related to \eqref{backeqn}. A similar version of this definition was given in \cite[Definition 3.3]{MR4616683}. 
\defn{\label{defnF0}\rm Fix $\beta_0\in(-1,1)$ and $t,x_0\in\R$. Set $\gamma_0=\frac{1}{\sqrt{1-\beta_0^2}}$ and $a_0=\sqrt{\frac{1+\beta_0}{1-\beta_0}}$. Set $(Q_0,Q_1)(x)=(Q,Q_t)(t,x;\beta_0,x_0)$. We define  $\F=(\F_1,\F_2,\F_3)(\delta,y,v_0,v_1,u_0,u_1)$ for \fm{(\delta,y,v_0,v_1,u_0,u_1)\in (-a_0,\infty)\times \R\times  W^{1,p}\times L^p\times W^{1,p}\times L^p,\qquad p\in[1,\infty] } by the following formula:
\eq{\label{defnF}&(\F_1(\delta,v_0,v_1,u_0),\F_2(\delta,v_0,u_0,u_1),\F_3(\delta,y,u_0))\\
&:=\begin{pmatrix}
\displaystyle Q_{0,x}+u_{0,x}-v_1-\frac{1}{a(\delta)}\sin(\frac{u_0+Q_0+v_0}{2})-a(\delta)\sin(\frac{u_0+Q_0-v_0}{2})\\[1em]
\displaystyle Q_{1}+u_{1}-v_{0,x}-\frac{1}{a(\delta)}\sin(\frac{u_0+Q_0+v_0}{2})+a(\delta)\sin(\frac{u_0+Q_0-v_0}{2})\\[1em]
\displaystyle \int_\R (u_0(x)+Q_0(x))\sech(\gamma(\delta)\cdot (x-\beta_0t-x_0-y))\ dx-\frac{\pi^2}{\gamma(\delta)}
\end{pmatrix}
.}
Here $a(\delta)=a_0+\delta>0$, $\beta(\delta)=\frac{a(\delta)^2-1}{a(\delta)^2+1}\in(-1,1)$ and $\gamma(\delta)=\frac{1}{\sqrt{1-\beta(\delta)^2}}$. 

Later, we  write $\F=\F_{(\beta_0,t,x_0)}$ to emphasize the dependence of $\F$ on the three parameters.
}

\rmk{\rm The functional $\F$ is related to the B\"acklund transform and the orthogonal conditions \eqref{asystacenter} and \eqref{mthmc1} in the main theorems. Suppose that $f,\phi$ are two solutions to \eqref{sg} such that $f$ is the B\"acklund transform of $\phi$ by the parameter $a>0$ (which is associated to a velocity $\beta\in(-1,1)$ by $a=\sqrt{\frac{1+\beta}{1-\beta}}$). Also suppose that for some path $x(t)$,
\fm{\int_\R (f(t,x)-Q(t,x;\beta,x(t)))\sech(\gamma(x-\beta t-x(t)))\ dx=0,\qquad \forall t\in\R.} Then, for all $t,x_0\in\R$ we have \eq{\label{Frelation}\F_{(\beta,t,x_0)}(0,x(t)-x_0,\phi(t),\phi_t(t),f(t)-Q(t,\cdot;\beta,x_0),f_t(t)-Q_t(t,\cdot;\beta,x_0))=0;}
for all $x_0\in\R$ and $\beta_0\in(-1,1)$, by setting $a_0=\sqrt{\frac{1+\beta_0}{1-\beta_0}}$ we have
\eq{\label{Frelation2}\F_{(\beta_0,0,x_0)}(a-a_0,x(0)-x_0,\phi(0),\phi_t(0),f(0)-Q(0,\cdot;\beta_0,x_0),f_t(0)-Q_t(0,\cdot;\beta_0,x_0))=0.}
Here we use the identity \eq{\label{intidentity}\int_\R Q(t,x;\beta,x_0)\sech(\gamma(x-\beta t-x_0)) \ dx=\gamma^{-1}\pi^2,\qquad \forall t,x_0\in\R.}

Since there are several different ways to choose the center $x(t)$, there might be also several different ways to define $\F_3$. We refer our readers to Remark \ref{rmkcenterdiff} on this topic.}

\rmk{\rm By part a) of Lemma \ref{lemQ}, we have 
\eq{\label{defnFcor}\F_1&=u_{0,x}-v_1-\frac{1}{a(\delta)}(\sin(\frac{u_0+Q_0+v_0}{2})-\sin(\frac{Q_0}{2}))+\frac{\delta}{a_0a(\delta)}\sin(\frac{Q_0}{2})\\
&\quad -a(\delta)(\sin(\frac{u_0+Q_0-v_0}{2})-\sin(\frac{Q_0}{2}))-\delta\sin(\frac{Q_0}{2}),\\
\F_2&=u_{1}-v_{0,x}-\frac{1}{a(\delta)}(\sin(\frac{u_0+Q_0+v_0}{2})-\sin(\frac{Q_0}{2}))+\frac{\delta}{a_0a(\delta)}\sin(\frac{Q_0}{2})\\
&\quad +a(\delta)(\sin(\frac{u_0+Q_0-v_0}{2})-\sin(\frac{Q_0}{2}))+\delta\sin(\frac{Q_0}{2}).}
Thus,  $\F$ is a continuously differentiable (i.e.\ $C^1$) map from $(-a_0,\infty)\times\R\times W^{1,p}\times L^p\times W^{1,p}\times L^p$ to $L^p\times L^p\times\R$ for each $p\in[1,\infty]$. It is also clear that $\F(0,0,0,0,0,0)=0$.

Note that \eqref{defnFcor} also implies the following fact. Fix $p\in(1,\infty]$. Suppose that for some $(\delta,v_0,v_1,u_0,u_1)\in (-a_0/2,\infty)\times W^{1,p}\times L^p\times W^{1,p}\times L^p$ such that $(\F_1,\F_2)(\delta,v_0,v_1,u_0,u_1)=0$ (note that $\F_1,\F_2$ do not depend on the variable $y$). By the Gagliardo-Nirenberg  inequality, we have $W^{1,p}(\R)\subset L^q(\R)$ for each $1<p\leq q\leq\infty$.  By \eqref{defnFcor}, we thus have $u_{0,x}-v_1,u_1-v_{0,x}\in L^q$. We also obtain
\eq{\label{defnFcor:cor}\norm{u_{0}}_{L^q}+\norm{v_0}_{L^q}+\norm{u_{0,x}-v_1}_{L^q}+\norm{v_{0,x}-u_1}_{L^q}&\lesssim_{p,q}\norm{u_0}_{W^{1,p}}+\norm{v_0}_{W^{1,p}}+|\delta|.}Recall that here $1<p\leq q\leq\infty$.
This estimate will be useful later.
}

\rmk{\label{rmkdefnF0tr}\rm In this section we seek to solve the equation $\F=0$. One notices that $\F$ admits translation invariance:
\fm{\F_{(\beta_0,t,x_0)}(\delta,y,v_0,v_1,u_0,u_1)&=\mcl{T}_{x_0+\beta_0t}\circ\F_{(\beta_0,0,0)}\circ\mcl{T}_{-x_0-\beta_0t}(\delta,y,v_0,v_1,u_0,u_1).}
Here for each function $h=h(x)$, we have $(\mcl{T}_{z}h)(x)=h(x-z)$. For each real number $C$, we have  $\mcl{T}_*C=C$ since a real number can be viewed as a constant function. 
Later we assume that $v_0,v_1,u_0,u_1$  belong to some  Banach spaces whose norms are translation-invariant (e.g.\ $H^m$ for $m\geq 0$). Thus, the results induced from $\F_{(\beta_0,t,x_0)}=0$ usually follow from those induced from  $\F_{(\beta_0,0,0)}=0$.
}

\rmk{\rm In the definition of $\F_3$, we use an unmatched pair of coefficients $(\gamma(\delta),\beta_0)$ in the weight $\sech(\gamma(\delta)(x-\beta_0t-x_0-y))$. Instead, the weight $\sech(\gamma(\delta)(x-\beta(\delta)t-x_0-y))$ seems to be a more natural option. However,  the translation invariance stated in the previous remark fails if we use the new weight. This would lead to a more complicated statement of Proposition \ref{propbackt}. Besides, in this paper, we are only interested in studying $\F$ in two cases: when $t=0$ or when $\delta=0$. In both cases, the two weights above coincide.  }
\rm

\subsection{A key proposition on $\F$}
We now state the key proposition of Section \ref{secbackt} on the functional $\F$. It states that the equation $\F=0$ induces a local $C^1$ diffeomorphism between a neighborhood of the kink $(Q,Q_t)(t,\cdot;\beta_0,x_0)$ and a neighborhood of the zero solution, the velocity $\beta_0$ and the center $x_0$.

\prop{\label{propbackt} Fix $\beta_0\in(-1,1)$, $t,x_0\in\R$ and $p\in[1,\infty]$. Define $\gamma_0,a_0,\beta(\delta),Q_0,Q_1$ as in Definition \ref{defnF0}. For each $r>0$, we define \fm{Y_{0}(r)&:=\{(u_0,u_1)\in W^{1,p}\times L^p:\ \norm{u_0}_{W^{1,p}}+\norm{u_1}_{L^p}<r\},\\
Y_{1}(r)&:=\{(\delta, y,v_0,v_1)\in \R\times\R\times W^{1,p}\times L^p:\ \norm{(\delta,y,v_0,v_1)}_{\R\times\R\times W^{1,p}\times L^p}<r\}.}
Then, there exist two small constants $r_0,r_1\in(0,1)$ and two $C^1$ maps
$\Phi=\Phi_{(\beta_0,t,x_0)}:Y_{0}(r_0)\to \R\times \R\times W^{1,p}\times L^p$ and $\Psi=\Psi_{(\beta_0,t,x_0)}:Y_1(r_1)\to  W^{1,p}\times L^p$, such that
\eq{\label{defnPhiPsi}\F_{(\beta_0,t,x_0)}(\Phi(u_0,u_1),u_0,u_1)&=0,\qquad \forall (u_0,u_1)\in Y_0(r_0),\\
\F_{(\beta_0,t,x_0)}(\delta,y,v_0,v_1,\Psi(\delta,y,v_0,v_1))&=0,\qquad \forall (\delta,y,v_0,v_1)\in Y_1(r_1)}
and that
\eq{\label{estPhiPsi}&\norm{\Phi(u_0,u_1)}_{\R\times\R\times W^{1,p}\times L^p}\lesssim \norm{(u_0,u_1)}_{W^{1,p}\times L^p},\quad \forall (u_0,u_1)\in Y_0(r_0),\\
&\norm{\Psi(\delta,y,v_0,v_1)}_{W^{1,p}\times L^p}\lesssim \norm{(\delta,y,v_0,v_1)}_{\R\times\R\times W^{1,p}\times L^p},\quad \forall (\delta,y,v_0,v_1)\in Y_1(r_1).}
Here $r_0,r_1$ and the implicit constants in \eqref{estPhiPsi} are independent of the choice of $(t,x_0)$.

Moreover, there exists a small constant $r_2\in(0,1)$ which is  independent of the choice of $(t,x_0)$,  such that $(\Psi\circ \Phi)|_{Y_0(r_2)}=\text{Id}$ and $(\Phi\circ \Psi)|_{Y_1(r_2)}=\text{Id}$. Thus, $\Phi,\Psi$ are $C^1$ diffeomorphisms between an open neighborhood of $(0,0)$ in $W^{1,p}\times L^p$ and an open neighborhood of $(0,0,0,0)$ in $\R\times\R\times W^{1,p}\times L^p$.

Finally, the maps $\Phi,\Psi$ are independent of the choice of $p$ in the following sense. Fix $\beta_0,t,x_0$. For $1<p\leq q\leq\infty$, we let $(\Phi_p,\Psi_p)$ and $(\Phi_q,\Psi_q)$ be the maps associated to $p$ and $q$, respectively. Then, there exists a small constant $r_{p,q}\in(0,1)$, depending on $p,q$ but not on $(t,x_0)$, such that
\eq{\label{PhiPsiindep}\Phi_p(u_0,u_1)=\Phi_q(u_0,u_1),\quad&\text{if }\norm{(u_0,u_1)}_{W^{1,p}\times L^p}+\norm{(u_0,u_1)}_{W^{1,q}\times L^q}<r_{p,q};\\
\Psi_p(\delta,y,v_0,v_1)=\Psi_q(\delta,y,v_0,v_1),\quad&\text{if }|\delta|+|y|+\norm{(v_0,v_1)}_{W^{1,p}\times L^p}+\norm{(v_0,v_1)}_{ W^{1,q}\times L^q}<r_{p,q}.}
}

\rmk{\label{rmkbacktimp}\rm To construct $\Phi$ and $\Psi$, we apply the implicit function theorem that implies the following fact.  For each fixed $(\beta_0,t,x_0)$, there exists a small constant $\wt{r}\in(0,1)$ (independent of $(t,x_0)$), such that if
\fm{\F(\delta,y,v_0,v_1,u_0,u_1)=0\text{ and }\norm{(\delta,y,v_0,v_1,u_0,u_1)}_{\R\times\R\times W^{1,p}\times L^p\times W^{1,p}\times L^p}<\wt{r},}
then we must have
\fm{(\delta,y,v_0,v_1)=\Phi(u_0,u_1),\qquad (u_0,u_1)=\Psi(\delta,y,v_0,v_1).}}
\rmk{\rm
\label{rmktx0general}If we have proved Proposition \ref{propbackt} in the case  $t=x_0=0$, then for general $(t,x_0)$, we can set\fm{\Phi_{(\beta_0,t,x_0)}(u_0,u_1)&=\mcl{T}_{\beta_0t+x_0}\Phi_{(\beta_0,0,0)}(\mcl{T}_{-x_0-\beta_0t}u_0,\mcl{T}_{-x_0-\beta_0t}u_1),\\
\Psi_{(\beta_0,t,x_0)}(\delta,y,v_0,v_1)&=\mcl{T}_{\beta_0t+x_0}\Psi_{(\beta_0,0,0)}(\delta,y,\mcl{T}_{-x_0-\beta_0t}v_0,\mcl{T}_{-x_0-\beta_0t}v_1).}
We refer our readers to Remark \ref{rmkdefnF0tr}. Here recall that the $W^{1,p}$ and $L^p$ norms are preserved by spatial translations, so  $r_0,r_1,r_2$ and all the implicit constants in \eqref{estPhiPsi} are independent of the choice of $(t,x_0)$.
}
\\
\rm

From now on we assume $t=x_0=0$. Let us first construct the map $\Phi$.
\lem{\label{lemexistPhi}Fix $\beta_0\in(-1,1)$ and $p\in[1,\infty]$. Set $t=x_0=0$. Define $\gamma_0,a_0,Q_0,Q_1,Y_0(r)$ as in Definition \ref{defnF0}. Then, there exist a small constant $r_0\in(0,1)$ and a  $C^1$ map $\Phi:Y_0(r_0)\to \R\times\R\times W^{1,p}\times L^p$ such that $\F(\Phi(u_0,u_1),u_0,u_1)=0$ and $\norm{\Phi(u_0,u_1)}_{\R\times\R\times W^{1,p}\times L^p}\lesssim \norm{(u_0,u_1)}_{W^{1,p}\times L^p}$ for all $(u_0,u_1)\in Y_0(r_0)$.}
\begin{proof}
Given $(u_0,u_1)\in Y_0(r_0)$ for some small $r_0>0$, we seek to find $(\delta,y,v_0,v_1)\in\R\times\R\times W^{1,p}\times L^p$ such that $\F(\delta,y,v_0,v_1,u_0,u_1)=0$.  First, we solve the equation $\F_2(\delta,v_0,u_0,u_1)=0$ by the implicit function theorem. Once this step is finished, we obtain a pair $(\delta,v_0)$. Next, we solve the equation $\F_1(\delta,v_0,v_1,u_0)=0$. Note that $\delta,v_0$ are known from the previous step, so $v_1$ is the only unknown in this step. Finally, we solve the equation $\F_3(\delta,y,u_0)=0$ for $y$ by the implicit function theorem. 

Let us start by solving $\F_2=0$ for $(\delta,v_0)$. Note that $\F_2$ is a $C^1$ map from $(-a_0,\infty)\times W^{1,p}\times W^{1,p}\times L^p$ to $L^p$ and  that $\F_2(0,0,0,0)=0$. Moreover, the differential $D_{(\delta,v_0)}\F_2(0,0,0,0)$ is a linear map
\fm{(\lambda,w)\mapsto-w_x-(\frac{1}{2a_0}+\frac{a_0}{2})\cos(\frac{Q_0}{2})w+(\frac{\lambda}{a_0^2}+\lambda)\sin(\frac{Q_0}{2}).}
By part a) of  Lemma \ref{lemQ}, we easily check that this map is bounded from $\R\times W^{1,p}$ to $L^p$. We claim that it is invertible and has a bounded inverse. By explicitly solving the ODE\eq{\label{leminitsoleqn}-w_x-(\frac{1}{2a_0}+\frac{a_0}{2})\cos(\frac{Q_0}{2})w+(\frac{\lambda}{a_0^2}+\lambda)\sin(\frac{Q_0}{2})=g\in L^p,}
we have 
\fm{w(x)&=\cosh(\gamma_0 x)(C+\int_0^x\kh{\lambda(1+a_0^{-2})\sech(\gamma_0 y)^2-g(y)\sech(\gamma_0 y)}\ dy).}
The integral on the right hand side converges absolutely because $g\in L^p$ and $\sech(\gamma_0\cdot)\in \mcl{S}$ (the space of Schwartz functions).
In order to guarantee $w\in W^{1,p}$, we choose $\lambda\in\R$ so that
\fm{\int_\R\kh{\lambda(1+a_0^{-2})\sech(\gamma_0 y)^2-g(y)\sech(\gamma_0 y)}\ dy=0.}
It is clear that such a $\lambda$ exists and that $|\lambda|\lesssim \norm{g}_{L^2}$. Now if $x\geq 0$, we write
\fm{w(x)&=-\cosh(\gamma_0 x) \int_x^{\infty}\kh{\lambda(1+a_0^{-2})\sech(\gamma_0 y)^2-g(y)\sech(\gamma_0y)}\ dy.}
Whenever $y\geq x\geq 0$ we have $0<\cosh(\gamma_0x)/\cosh(\gamma_0y)\lesssim e^{\gamma_0 (x-y)}$. 
It follows that
\fm{|w|&\lesssim \int_x^\infty |\lambda| \sech(\gamma_0 x)e^{2\gamma(x-y)}\ dy+\int_{x}^\infty|g(y)|e^{-\gamma_0|x-y|}\ dy\lesssim |\lambda|\sech(\gamma_0x)+(e^{-\gamma_0|\cdot|}*|g|)(x).}
Similarly, we obtain the same estimate for $x\leq 0$.
Young's convolution inequality yields 
\fm{\norm{w}_{L^p}\lesssim |\lambda|+\norm{e^{-\gamma_0|\cdot|}}_{L^1}\norm{g}_{L^p}\lesssim \norm{g}_{L^p}.}
Using \eqref{leminitsoleqn}, we have
\fm{\norm{w_x}_{L^p}&\lesssim \norm{w}_{L^p}+|\lambda|+\norm{g}_{L^p}\lesssim \norm{g}_{L^p}.}
The claim above is thus proved. By the implicit function theorem, we obtain a small $r_0\in(0,1)$ and a  $C^1$ map $\Phi_1$. Set $(\delta,v_0)=\Phi_1(u_0,u_1)$. Since a  $C^1$ map is locally Lipschitz, we obtain the bound $|\delta|+\norm{v_0}_{H^{1}}\lesssim \norm{u_0}_{H^1}+\norm{u_1}_{L^2}$.

Next, we define $v_1$ by the equation $\F_1=0$. The first formula in \eqref{defnFcor} yields 
\fm{\norm{v_1}_{L^p}&\lesssim \norm{u_{0,x}}_{L^p}+\norm{u_0}_{L^p}+\norm{v_0}_{L^p}+|\delta|\lesssim \norm{u_0}_{W^{1,p}}+\norm{u_1}_{L^p}.}

Finally, we solve $\F_3=0$ for $y$. Here $\F_3$ is a $C^1$ map from $(-a_0,\infty)\times\R\times W^{1,p}$ to $\R$, and we have $\F_3(0,0,0)=0$. Moreover, we have 
\fm{\partial_y\F_3(0,0,0)&=\int_\R Q_0(x)\partial_y(\sech(\gamma_0(x-y)))|_{y=0}\ dx=-\int_\R Q_0(x)\partial_x(\sech(\gamma_0 x))\ dx\\
&=\int_\R 2\gamma_0\sech(\gamma_0x)^2\ dx=4.}
By the implicit function theorem, we obtain a small constant $\wt{r}_0\in(0,1)$ and a  $C^1$ map $\Phi_2:\{(\kappa,w)\in\R\times W^{1,p}:\ |\kappa|+\norm{w}_{W^{1,p}}<\wt{r}_0\}\to\R$ such that $\F_3(\kappa,\Phi_2(\kappa,w),w)=0$. For $(\delta,v_0)=\Phi_1(u_0,u_1)$, we have $|\delta|+\norm{u_0}_{W^{1,p}}\lesssim r_0$. By choosing $r_0\ll1$, we have $|\delta|+\norm{u_0}_{W^{1,p}}< \wt{r}_0$ which allows us to set $y=\Phi_2(\delta,u_0)$. We also have $|y|\lesssim |\delta|+ \norm{u_0}_{W^{1,p}}\lesssim \norm{u_0}_{W^{1,p}}+\norm{u_1}_{L^p}$. 
\end{proof}
\rmk{\rm In order to solve $\F(\delta,y,v_0,v_1,u_0,u_1)=0$ for $(\delta,y,v_0,v_1)$, one can also compute the differential $D_{(\delta,y,v_0,v_1)}\F(0,0,0,0,0,0)$ and then show that the implicit function theorem is applicable to $\F$. However, we believe that our proof above is more efficient. By first solving $\F_2(\delta,v_0,u_0,u_1)=0$, we avoid computing the differentials $D_{(\delta,v_0)}\F_1(0,0,0,0)$ and $\partial_\delta\F_3(0,0,0)$. This makes our proof simpler.

A similar remark can be made for Lemma \ref{lemexistPsi} below.} 

\rmk{\label{rmkexistPhi}\rm The last paragraph in the proof and Remark \ref{rmktx0general} indicate that, for each fixed $\beta_0\in(-1,1),t,x_0\in\R$ and $p\in[1,\infty]$, there exist a small constant $r_0\in(0,1)$ and a unique $C^1$ map \fm{\zeta=\zeta_{(\beta_0,t,x_0)}:\{u_0\in W^{1,p}:\ \norm{u_0}_{W^{1,p}}<r_0\}\to\R,}
such that $|\zeta(u_0)-x_0|\lesssim \norm{u_0}_{W^{1,p}}$ and 
\fm{\int_\R (u_0(x)+Q(t,x;\beta_0,x_0)-Q(t,x;\beta_0,\zeta(u_0)))\sech(\gamma_0(x-\beta _0t-\zeta(u_0)))\ dx=0.}
Moreover, the constant $r_0$ and the implicit constant in the estimate are independent of $(t,x_0)$.
In fact, using the notations above, we can take  $\zeta_{(\beta_0,t,x_0)}(\cdot)=x_0+(\Phi_2)_{(\beta_0,t,x_0)}(0,\cdot)$.
}\\
\rm

We then construct the map $\Psi$. 

\lem{\label{lemexistPsi} Fix $\beta_0\in(-1,1)$ and $p\in[1,\infty]$. Set $t=x_0=0$. Define $\gamma_0,a_0,Q_0,Q_1,Y_1(r)$ as in Definition \ref{defnF0}. Then, there exist a small constant $r_1\in(0,1)$ and a  $C^1$ map $\Psi:Y_1(r_1)\to  W^{1,p}\times L^p$ such that $\F(\delta,y,v_0,v_1,\Psi(\delta,y,v_0,v_1))=0$ and $\norm{\Psi(\delta,y,v_0,v_1)}_{ W^{1,p}\times L^p}\lesssim \norm{(\delta,y,v_0,v_1)}_{\R\times\R\times W^{1,p}\times L^p}$ for all $(\delta,y,v_0,v_1)\in Y_1(r_1)$.}
\begin{proof}
Given $(\delta,y,v_0,v_1)\in Y_1(r_1)$ for some small $r_1>0$, we seek to find $(u_0,u_1)\in W^{1,p}\times L^p$ such that $\F(\delta,y,v_0,v_1,u_0,u_1)=0$. We start by solving the equation $(\F_1,\F_3)(\delta,y,v_0,v_1,u_0)=0$ for $u_0$. We will apply the implicit function theorem in this step. Then, we  define $u_1$ by the equation $\F_2(\delta,v_0,u_0,u_1)=0$.

Let us first solve $(\F_1,\F_3)=0$ for $u_0$.  Note that $(\F_1,\F_3)$ is a $C^1$ map from $(-a_0,\infty)\times \R\times W^{1,p}\times L^p\times W^{1,p}$ to $L^p\times\R$ and that the differential $D_{(u_0)}(\F_1,\F_3)(0,0,0,0,0)$ is a linear map
\fm{w\mapsto (w_x-\frac{1}{2}(a_0+\frac{1}{a_0})\cos(\frac{Q_0}{2})w,\int_\R w(x)\sech(\gamma_0x)\ dx).}This map is bounded from $W^{1,p}$ to $L^p\times\R$. We claim that it is invertible and has a bounded inverse. To see this, we need to solve
\eq{\label{lemexistPsif1}(w_x-\frac{1}{2}(a_0+\frac{1}{a_0})\cos(\frac{Q_0}{2})w,\int_\R w(x)\sech(\gamma_0x)\ dx)=(g,\lambda)\in L^p\times\R.}
The identity on the first component indicates that
\fm{\partial_x(\cosh(\gamma_0x)w)=\cosh(\gamma_0x)g\Longrightarrow w(x)=\sech(\gamma_0x)(w(0)+\int_0^x \cosh(\gamma_0z)g(z)\ dz).}
Since $0<\sech(\gamma_0x)\cosh(\gamma_0z)\leq e^{-\gamma_0|x-z|}$ whenever $z$ lies between $0$ and $x$, we have
\fm{|w(x)|&\lesssim |w(0)|\sech(\gamma_0x)+\abs{\int_0^xe^{-\gamma_0|x-z|}|g(z)|\ dz}\lesssim |w(0)|\sech(\gamma_0x)+(e^{-\gamma_0|\cdot|}*|g|)(x).}
We obtain $\norm{w}_{L^p}\lesssim |w(0)|+\norm{g}_{L^p}$ by  Young's convolution inequality. By \eqref{lemexistPsif1} we have
\fm{\norm{w_x}_{L^p}&\lesssim \norm{w}_{L^p}+\norm{g}_{L^p}\lesssim |w(0)|+\norm{g}_{L^p}.}
The identity on the second component indicates that
\fm{\lambda&=\int_\R w(x)\sech(\gamma_0x)\ dx=\int_\R w(0)\sech(\gamma_0x)^2+\int_\R\sech(\gamma_0x)^2 \int_0^x \cosh(\gamma_0z)g(z)\ dz dx.}
Note that $\int_\R\sech(\gamma_0x)^2 \ dx=2/\gamma_0$ and that the second integral on the right side is bounded by 
\fm{\int_\R\sech(\gamma_0x)(e^{-\gamma_0|\cdot|}*|g|)(x) \ dx\lesssim \norm{g}_{L^p}.} Thus, one can solve for $w(0)$  and obtain an estimate $|w(0)|\lesssim |\lambda|+\norm{g}_{L^p}$.
In summary, we  have  $\norm{w}_{W^{1,p}}\lesssim|\lambda|+\norm{g}_{L^2}$. This claim above is thus proved. By the implicit function theorem, we obtain a small $r_1\in(0,1)$ and a  $C^1$ map $\Psi_1:Y_1(r_1)\to W^{1,p}$, such that $(\F_1,\F_3)(\delta,y,v_0,v_1,\Psi_1(\delta,y,v_0,v_1))=0$ and $\norm{\Psi_1(\delta,y,v_0,v_1)}_{W^{1,p}}\lesssim \norm{(\delta,y,v_0,v_1)}_{\R\times\R\times W^{1,p}\times L^p}$.

We then use $\F_2=0$ to define $u_1$. By the second formula in \eqref{defnFcor}, we have
\fm{\norm{u_1}_{L^p}&\lesssim \norm{v_{0,x}}_{L^p}+\norm{u_0}_{L^p}+\norm{v_0}_{L^p}+|\delta|\lesssim \norm{(\delta,y,v_0,v_1)}_{\R\times\R\times W^{1,p}\times L^p}.}
This finishes the proof.
\end{proof}

\rm

Next, we check that $\Phi$ and $\Psi$ are locally invertible to each other.
\lem{There exists a small constant $r_2\in(0,1)$ such that $(\Psi\circ \Phi)|_{Y_0(r_2)}=\text{Id}$ and $(\Phi\circ \Psi)|_{Y_1(r_2)}=\text{Id}$. Thus, $\Phi$ is  a  $C^1$ diffeomorphism between an open neighborhood of $(0,0)$ in $W^{1,p}\times L^p$ and an open neighborhood of $(0,0,0,0)$ in $\R\times\R\times W^{1,p}\times L^p$.}
\begin{proof}
Let $0<r_2<\min\{r_0,r_1\}$ be a small constant to be chosen later.
Fix $(u_0,u_1)\in Y_0(r_2)$ and we set $(\delta,y,v_0,v_1):=\Phi(u_0,u_1)$. Since $\norm{\Phi(u_0,u_1)}_{\R\times \R\times W^{1,p}\times L^p}\lesssim\norm{(u_0,u_1)}_{W^{1,p}\times L^p}\lesssim r_2$, by choosing $r_2\ll1$ we have $\Phi(u_0,u_1)\in Y_1(r_1)$. Set $(w_0,w_1):=\Psi(\delta,y,v_0,v_1)$. By the definitions of $\Phi$ and $\Psi$, we have
\fm{\F(\delta,y,v_0,v_1,u_0,u_1)=\F(\delta,y,v_0,v_1,w_0,w_1)=0.} 
We have $\norm{(w_0,w_1)}_{W^{1,p}\times L^p}\lesssim \norm{\Phi(u_0,u_1)}_{\R\times \R\times W^{1,p}\times L^p}\lesssim \norm{(u_0,u_1)}_{W^{1,p}\times L^p}\lesssim r_2$. By choosing $r_2\ll1$, we force  $(w_0,w_1)=(u_0,u_1)$ by the uniqueness result in Remark \ref{rmkbacktimp}. That is, $(\Psi\circ \Phi)|_{Y_0(r_2)}=\text{Id}$. Similarly, by choosing $r_2\ll1$, we have $(\Phi\circ \Psi)|_{Y_1(r_2)}=\text{Id}$.
\end{proof}\rm

Finally, we check that the choices of $\Phi$ and $\Psi$ are independent of the parameter $p$.
\lem{Fix $(\beta_0,t,x_0)$ as above. Fix $\beta_0,t,x_0$. For $1<p\leq q\leq\infty$, we let $(\Phi_p,\Psi_p)$ and $(\Phi_q,\Psi_q)$ be the maps associated to $p$ and $q$, respectively. Then, there exists a small constant $r_{p,q}\in(0,1)$, depending on $p,q$ but not on $(t,x_0)$, such that
\fm{\Phi_p(u_0,u_1)=\Phi_q(u_0,u_1),\quad&\text{if }\norm{(u_0,u_1)}_{W^{1,p}\times L^p}+\norm{(u_0,u_1)}_{W^{1,q}\times L^q}<r_{p,q};\\
\Psi_p(\delta,y,v_0,v_1)=\Psi_q(\delta,y,v_0,v_1),\quad&\text{if }|\delta|+|y|+\norm{(v_0,v_1)}_{W^{1,p}\times L^p}+\norm{(v_0,v_1)}_{ W^{1,q}\times L^q}<r_{p,q}.}}
\begin{proof}
Let us assume $p<q$ since there is nothing to prove when $p=q$.
Fix $(u_0,u_1)$ with $\norm{(u_0,u_1)}_{W^{1,p}\times L^p}+\norm{(u_0,u_1)}_{W^{1,q}\times L^q}<r_{p,q}$. The value of $r_{p,q}\in(0,1)$ will be chosen later. Set $(\delta,y,v_0,v_1)=\Phi_p(u_0,u_1)$. Recall that $|\delta|+|y|+\norm{(v_0,v_1)}_{W^{1,p}\times L^p}\lesssim \norm{(u_0,u_1)}_{W^{1,p}}\lesssim r_{p,q}$. By \eqref{defnFcor:cor} in a previous remark, we have
\fm{\norm{(v_0,v_1)}_{W^{1,q}\times L^q}\lesssim \norm{(u_0,u_1)}_{W^{1,p}\times L^p}+\norm{(u_0,u_1)}_{W^{1,q}\times L^q}\lesssim r_{p,q}.}
The implicit constants here depend on $p,q$ but not on $t,x_0$. Thus, by choosing $r_{p,q}\ll1$, we can apply Remark \ref{rmkexistPhi} (with $p$ replaced by $q$) to conclude that $\Phi_p(u_0,u_1)=\Phi_q(u_0,u_1)$. The proof of the other half is similar. 
\end{proof}\rm

\subsection{Solving the B\"acklund transform}
With the definition of the B\"acklund transform, we now seek to solve the B\"acklund transform. In the case when the initial data of $\phi$ is small in $H^1\times L^2$, or when the initial data of $f$ is close to a kink in $H^1\times L^2$, we can answer the two questions (Q1) and (Q2) stated at the beginning of this section by applying Proposition \ref{propbackt} (with $p=2$ and $t=0$).

\cor{\label{corsolvebackt}
We have the following conclusions.
\begin{enumerate}[\rm a)]
    \item Let $\phi\in {\bf E}_{\sin}^0$ be a solution to \eqref{sg}. Suppose that \fm{\norm{(\phi,\phi_t)(0)}_{H^1\times L^2}<\eps_0\ll1.} Fix $\beta\in(-1,1)$ and $y_0\in\R$. Then, as long as $\eps_0\ll1$, there exists a unique solution $f\in {\bf E}_{\sin}^0$  to \eqref{sg}, such that $f$ is the B\"acklund transform of $\phi$ by the parameter $a=\sqrt{\frac{1+\beta}{1-\beta}}$, that
    \fm{\norm{(f,f_t)(0,\cdot)-(Q,Q_t)(0,\cdot;\beta,y_0)}_{H^1\times L^2}\lesssim\norm{(\phi,\phi_1)}_{H^1\times L^2}\lesssim\eps_0}
    and that \fm{\int_\R(f(0,x)-Q(0,x;\beta,y_0))\sech(\gamma(x-y_0))\ dx=0,\qquad \gamma=(1-\beta^2)^{-1/2}.}
    \item Let $f\in {\bf E}_{\sin}^0$ be a solution to \eqref{sg}. Suppose that for some $\beta_0\in(-1,1)$ and $x_0\in\R$,
    \fm{\norm{(f,f_t)(0,\cdot)-(Q,Q_t)(0,\cdot;\beta_0,x_0)}_{H^1\times L^2}<\eps_0\ll1.}
    Then, as long as $\eps_0\ll1$, there exists a unique triple $(\beta,y_0,\phi)\in (-1,1)\times\R\times {\bf E}_{\sin }^0$, such that $\phi$ is a solution of \eqref{sg}, that $f$ is the B\"acklund transform of $\phi$ by the parameter $a=\sqrt{\frac{1+\beta}{1-\beta}}$, that
    \fm{|\beta-\beta_0|+|y_0-x_0|+\norm{(\phi,\phi_t)(0)}_{H^1\times L^2}\lesssim \norm{(f,f_t)(0,\cdot)-(Q,Q_t)(0,\cdot;\beta_0,x_0)}_{H^1\times L^2}\lesssim\eps_0} and that \fm{\int_\R(f(0,x)-Q(0,x;\beta,y_0))\sech(\gamma(x-y_0))\ dx=0,\qquad \gamma=(1-\beta^2)^{-1/2}.}
\end{enumerate}}\rm

It suffices to prove this corollary at $t=0$. In fact, since \eqref{sg} is globally wellposed in the energy space, a solution in ${\bf E}_{\sin}^0$ to \eqref{sg} is uniquely determined by its data at $t=0$. The next lemma implies that, if $f,\phi\in{\bf E}_{\sin}^0$ are solutions to \eqref{sg}, and if \eqref{backeqn} is satisfied at $t=0$, then \eqref{backeqn} holds everywhere.

\lem{\label{leminitalltime}Suppose that $f,\phi\in {\bf E}_{\sin}^0$ are both solutions to \eqref{sg} and that  \eqref{backeqn} holds at time $t=0$. Then,  \eqref{backeqn} holds for all time $t\in\R$; in other words, $f$ is the B\"acklund transform of $\phi$ by the parameter $a$. }
\begin{proof}
We first claim that $f_t-\phi_x,f_x-\phi_t\in C_tH^1_x(\R\times\R)$. To see this, we recall from  
$f,\phi\in{\bf E}_{\sin}^0$that  $\sin(f/2),\sin(\phi/2),f_x,f_t,\phi_x,\phi_t\in C_tL_x^2(\R\times\R)$. By the chain rule, we have $\sin(f),\sin(\phi),\sin((f\pm\phi)/2)\in C_tH^1_x(\R\times\R)$. Moreover, the equation \eqref{sg} yields
\fm{(\partial_t-\partial_x)(f_t+f_x-\phi_t-\phi_x)&=-\sin f+\sin\phi\in C_tH_x^1(\R\times\R),\\
(\partial_t+\partial_x)(f_t-f_x+\phi_t-\phi_x)&=-\sin f+\sin\phi\in C_tH_x^1(\R\times\R).}
Since \eqref{backeqn} holds at $t=0$, we have
\fm{(f_x+f_t-\phi_t-\phi_x,f_t-f_x+\phi_t-\phi_x)=(\frac{2}{a}\sin(\frac{f+\phi}{2}),-2a\sin(\frac{f-\phi}{2}))\in H^1(\R)}
at $t=0$. As a result, at each $(t_0,x_0)\in\R\times\R$ one has 
\fm{(f_t+f_x-\phi_t-\phi_x)(t_0,x_0)&=\frac{2}{a}\sin(\frac{(f+\phi)(0,x_0+t_0)}{2})+\int_0^{t_0}(\sin\phi-\sin f)(t,x_0+t_0-t)\ dt.}
The right side is well defined pointwisely because we have $C_tH^1_x(\R\times\R)\subset C_{t,x}(\R\times\R)$ by the Gagliardo-Nirenberg inequality.  It is now easy to check that $f_t+f_x-\phi_t-\phi_x\in C_tH^1_x(\R\times\R)$. Similarly, we have $f_t-f_x+\phi_t-\phi_x\in C_tH^1_x(\R\times\R)$.

Next, we set
\fm{R_1&=f_t+f_x-\phi_t-\phi_x-\frac{2}{a}\sin(\frac{f+\phi}{2}),\\
R_2&= f_t-f_x+\phi_t-\phi_x+2a\sin(\frac{f-\phi}{2}).}
So far, we have proved
$R_1,R_2\in C_tH^1_x(\R\times\R)$. By \eqref{sg}, we have
\fm{\left\{
\begin{array}{l}
    \displaystyle (\partial_t-\partial_x)R_1=-\frac{1}{a}\cos(\frac{f+\phi}{2})R_2, \\[1em]
    \displaystyle (\partial_t+\partial_x)R_2=-a\cos(\frac{f-\phi}{2})R_1,\\[1em]
    R_1=R_2=0,\qquad \text{at }t=0.
\end{array}
\right.}
By writing $R_1$ as an integral of $\cos(\frac{f+\phi}{2})R_2$ along a curve $x+t=\textrm{constant}$ and writing $R_2$ as an integral of $\cos(\frac{f-\phi}{2})R_1$ along a curve $x-t=\textrm{constant}$, we obtain
\fm{\norm{R_1(t_0)}_{L^\infty}\lesssim \abs{\int_0^{t_0}\norm{R_2(t)}_{L^\infty}\ dt},\quad \norm{R_2(t_0)}_{L^\infty}\lesssim \abs{\int_0^{t_0}\norm{R_1(s)}_{L^\infty}\ ds},\qquad \forall t_0\in\R.}
Since $R_1,R_2\in C_tH^1_x(\R\times\R)$, the map $t\mapsto \norm{R_1(t)}_{L^\infty}+\norm{R_2(t)}_{L^\infty}$ is  continuous. By  Gronwall's inequality, we conclude  $R_1=R_2=0$ everywhere and thus \eqref{backeqn} everywhere.
\end{proof}
\rm

We now prove Corollary \ref{corsolvebackt} at $t=0$. To prove part a), we apply Proposition \ref{propbackt} (with $\beta_0$ replaced by $\beta$). We thus obtain a $C^1$ map $\Psi:Y_1(r_1)\to H^1\times L^2$. By choosing $\eps_0<r_1$, we have $(0,0,\phi(0),\phi_t(0))\in Y_1(r_1)$, so we obtain $(u_0,u_1)=\Psi(0,0,\phi(0),\phi_t(0))$. Set
\fm{(f,f_t)|_{t=0}=(u_0,u_1)+(Q,Q_t)(0,x;\beta,x_0).}
We  now finish the proof by  the equation $\F(0,0,\phi(0),\phi_t(0),u_0,u_1)=0$ and the estimates in Proposition \ref{propbackt}. 

To prove part b), we apply Proposition \ref{propbackt} and obtain a $C^1$ map $\Phi:Y_0(r_0)\to (-a_0,\infty)\times\R\times H^1\times L^2$. By choosing $\eps_0<r_0$, we have $(u_0,u_1)=(f,f_t)|_{t=0}-(Q,Q_t)(0,x;\beta_0,x_0)\in Y_0(r_0)$, so we obtain $(\delta,y,v_0,v_1)=\Phi(u_0,u_1)$. Set $y_0=x_0+y$, $a=a_0+\delta$ and $(\phi,\phi_t)|_{t=0}=(v_0,v_1)$. We also define $\beta,\gamma$ accordingly. We now finish the proof by the equation $\F(\delta,y,v_0,v_1,u_0,u_1)=0$ and the estimates in Proposition \ref{propbackt}.

\subsection{Results involving weighted Sobolev norms}
In Theorem \ref{thmasysta}, we assume that the initial data of $f$ is close to a certain kink in a weighted Sobolev space. This in fact implies that the initial data of the corresponding $\phi$ (constructed in  Corollary \ref{corsolvebackt}) is also close to zero in the same weighted Sobolev space. This result follows from the following lemma.

\lem{\label{lemhighreginit}Fix  $\beta_0\in(-1,1)$ and $x_0\in\R$. Also fix two real numbers $m,s\geq 0$. Let $\Phi:Y_0(r_0)\to \R\times\R\times H^1\times L^2$ be the $C^1$ map constructed in Lemma \ref{lemexistPhi}. Also define
\fm{Y_{m,s}(r)=\{(u_0,u_1)\in H^{m+1,s}\times H^{m,s}:\ \norm{u_0}_{H^{m+1,s}}+\norm{u_1}_{H^{m,s}}<r\},\qquad r>0.}
Then, there exists a small constant $r_{m,s}\in(0,1)$, depending on $m,s,x_0,\beta_0$, such that we have a restriction 
$\Phi:{Y_{m,s}(r_{m,s})}\to \R\times\R\times H^{m+1,s}\times H^{m,s}$
with 
\fm{\norm{\Phi(u_0,u_1)}_{\R\times\R\times H^{m+1,s}\times H^{m,s}}\lesssim\norm{(u_0,u_1)}_{H^{m+1,s}\times H^{m,s}}.}
}
\begin{proof}
Set $(\delta,y,v_0,v_1)=\Phi(u_0,u_1)$. Recall that $|\delta|+|y|\lesssim \norm{(u_0,u_1)}_{H^1\times L^2}\lesssim r_{m,s}$. It remains to control the $H^{m+1,s}\times H^{m,s}$ norm of $(v_0,v_1)$. By  \eqref{defnFcor}, we can write $\F_2=0$ as \fm{-\cosh(\gamma_0(x-x_0))\cdot\partial_x(\sech(\gamma_0(x-x_0))v_{0})=-v_{0,x}-\gamma_0\cos(\frac{Q_0}{2})v=G}
and thus
\fm{v_0(x)&=\cosh(\gamma_0(x-x_0))(v_0(x_0)-\int_{x_0}^x G(z)\sech(\gamma_0(z-x_0))\ dz).}
We recall that $Q_0(x)=Q(0,x;\beta_0,x_0)$. Since
\fm{\sin(\frac{z+Q_0}{2})-\sin(\frac{Q_0}{2})&=\frac{1}{2}\cos(\frac{Q_0}{2})z+O(|z|^2)}
and $\norm{v_0}_{L^\infty}+\norm{u_0}_{L^\infty}\lesssim\norm{v_0}_{H^1}+\norm{u_0}_{H^1}\lesssim r_{m,s}\ll1$,  we have
\fm{G&=(\gamma(\delta)-\gamma_0)\cos(\frac{Q_0}{2})v_0+(\frac{1}{2a(\delta)}-\frac{a(\delta)}{2})\cos(\frac{Q_0}{2})u_0-u_1\\
&\quad-\delta(1+\frac{1}{a_0 a(\delta)})\sin(\frac{Q_0}{2})+O(|v|^2+|u_0|^2)\\
&=O(r_{m,s}|v_0|+|u_0|+|u_1|+|\delta|\sech(\gamma_0(\cdot-x_0)))\in L^2.}
Since $v_0\in H^1$, we have
\fm{\int_\R G(z)\sech(\gamma_0(z-x_0))\ dz=0}and
\fm{v_0(x)&=\cosh(\gamma_0(x-x_0))\int_x^\infty G(z)\sech(\gamma_0(z-x_0))\ dz\\
&=-\cosh(\gamma_0(x-x_0))\int_{-\infty}^x G(z)\sech(\gamma_0(z-x_0))\ dz.}
Note that
\fm{0<\frac{\cosh(\gamma_0(x-x_0))}{\cosh(\gamma_0(z-x_0))}\leq e^{-\gamma_0|x-z|},\quad \lra{z}\gtrsim_{x_0}\lra{x},\qquad \forall z\geq x\geq x_0\quad\text{or}\quad \forall z\leq x\leq x_0.}
It follows that
\fm{&\lra{x}^s|v_0(x)|\lesssim \lra{x}^s\int_\R|G(z)|e^{-\gamma_0|x-z|}\ dz\\
&\lesssim_{x_0,s}\int_\R\lra{z}^s(r_{m,s}|v_0|+|u_0|+|u_1|)e^{-\gamma_0|x-z|}\ dz+\int_\R|\delta|\lra{x}^s\sech(\gamma_0(x-x_0))e^{-2\gamma_0|x-z|}\ dz\\
&\lesssim [\lra{\cdot}^s(r_{m,s}|v_0|+|u_0|+|u_1|)]*e^{-\gamma_0|\cdot|}(x)+|\delta|\lra{x}^s\sech(\gamma_0(x-x_0)).}
By  Young's convolution inequality, we conclude that 
\fm{\norm{v_0}_{L^{2,s}}&\lesssim_{x_0,s} r_{m,s}\norm{v_0}_{L^{2,s}}+\norm{u_0}_{L^{2,s}}+\norm{u_1}_{L^{2,s}}+|\delta|\lesssim r_{m,s}\norm{v_0}_{L^{2,s}}+\norm{u_0}_{H^{1,s}}+\norm{u_1}_{L^{2,s}}. }
By choosing $r_{m,s}\ll_{m,s,x_0}1$, the term $r_{m,s}\norm{v_0}_{L^{2,s}}$ can be absorbed by the left side. Thus, we have $\norm{v_0}_{L^{2,s}}\lesssim \norm{u_0}_{H^{1,s}}+\norm{u_1}_{L^{2,s}}$.

To continue, we temporarily assume that $m\geq 0$ is an integer. In this case,  by \eqref{defnFcor} and part c) of Lemma \ref{lemQ}, we have
\fm{&\norm{v_0}_{H^{m+1,s}}+\norm{v_1}_{H^{m,s}}\lesssim \norm{v_{0,x}}_{H^{m,s}}+\norm{v_0}_{L^{2,s}}+\norm{v_1}_{H^{m,s}}
\\&\lesssim \norm{u_0}_{H^{m+1,s}}+\norm{u_1}_{H^{m,s}}+|\delta|\norm{\sech(\gamma_0(\cdot+x_0))}_{H^{m,s}}\\
&\quad+\norm{\sin(\frac{u_0+Q_0+v_0}{2})-\sin(\frac{Q_0}{2})}_{H^{m,s}}+\norm{\sin(\frac{u_0+Q_0-v_0}{2})-\sin(\frac{Q_0}{2})}_{H^{m,s}}\\
&\lesssim \norm{u_0}_{H^{m+1,s}}+\norm{u_1}_{H^{m,s}}+(1+\norm{u_0}_{H^{m,s}}+\norm{v_0}_{H^{m,s}})^{m}(\norm{u_0}_{H^{m,s}}+\norm{v_0}_{H^{m,s}}).}
If $m=0$, we conclude that $\norm{v_0}_{H^{1,s}}+\norm{v_1}_{L^{2,s}}\lesssim  \norm{u_0}_{H^{1,s}}+\norm{u_1}_{L^{2,s}}$.
If $m>0$,  we induct on $m$. If we have $\norm{v_0}_{H^{m,s}}+\norm{v_1}_{H^{m-1,s}}\lesssim \norm{u_0}_{H^{m,s}}+\norm{u_1}_{H^{m-1,s}}$, then  the estimate above implies that
\fm{\norm{v_0}_{H^{m+1,s}}+\norm{v_1}_{H^{m,s}}&\lesssim_{m,s,x_0} \norm{u_0}_{H^{m+1,s}}+\norm{u_1}_{H^{m,s}}+(1+r_{m,s})^{m}(\norm{u_0}_{H^{m,s}}+\norm{u_1}_{H^{m-1,s}})\\
&\lesssim_{m,s,x_0} \norm{u_0}_{H^{m+1,s}}+\norm{u_1}_{H^{m,s}}.}
This finishes the proof when $m\geq 0$ is an integer.

For a general noninteger  $m>0$, we set $k=\lfloor m\rfloor+1$. Since $m<k<m+1$,  we have \fm{&\norm{v_0}_{H^{m+1,s}}+\norm{v_1}_{H^{m,s}}\lesssim \norm{v_{0,x}}_{H^{m,s}}+\norm{v_0}_{L^{2,s}}+\norm{v_1}_{H^{m,s}}\\
&\lesssim \norm{u_0}_{H^{m+1,s}}+\norm{u_1}_{H^{m,s}}+|\delta|\norm{\sech(\gamma_0(\cdot+x_0))}_{H^{k,s}}\\
&\quad+\norm{\sin(\frac{u_0+Q_0+v_0}{2})-\sin(\frac{Q_0}{2})}_{H^{k,s}}+\norm{\sin(\frac{u_0+Q_0-v_0}{2})-\sin(\frac{Q_0}{2})}_{H^{k,s}}\\
&\lesssim \norm{u_0}_{H^{m+1,s}}+\norm{u_1}_{H^{m,s}}+(1+\norm{u_0}_{H^{k,s}}+\norm{v_0}_{H^{k,s}})^{k}(\norm{u_0}_{H^{k,s}}+\norm{v_0}_{H^{k,s}})\\
&\lesssim \norm{u_0}_{H^{m+1,s}}+\norm{u_1}_{H^{m,s}}+(1+\norm{u_0}_{H^{k,s}}+\norm{u_1}_{H^{k-1,s}})^{k}(\norm{u_0}_{H^{k,s}}+\norm{u_1}_{H^{k-1,s}})\\
&\lesssim \norm{u_0}_{H^{m+1,s}}+\norm{u_1}_{H^{m,s}}.}
In the second last estimate, we use the estimate proved in the integer case.
\end{proof}
\cor{\label{corsolvebackthigh}Let $f\in {\bf E}_{\sin}^0$ be a solution to \eqref{sg}. Fix two real numbers $m,s\geq 0$. Suppose that for some $\beta_0\in(-1,1)$ and $x_0\in\R$,\fm{\norm{(f,f_t)(0,\cdot)-(Q,Q_t)(0,\cdot;\beta_0,x_0)}_{H^{m+1,s}\times H^{m,s}}<\eps_0\ll1.}
Then, as long as $\eps_0\ll_{m,s,x_0,\beta_0}1$, there exists a unique triple $(\beta,y_0,\phi)\in (-1,1)\times\R\times {\bf E}_{\sin }^0$, such that $\phi$ is a solution of \eqref{sg}, that $f$ is the B\"acklund transform of $\phi$ by the parameter $a=\sqrt{\frac{1+\beta}{1-\beta}}$, that \fm{|\beta-\beta_0|+|y_0-x_0|+\norm{(\phi,\phi_t)(0)}_{H^{m+1,s}\times H^{m,s}}\lesssim \eps_0,} 
that
\fm{\norm{(f,f_t)(0,\cdot)-(Q,Q_t)(0,\cdot;\beta,y_0)}_{H^{m+1,s}\times H^{m,s}}\lesssim\eps_0\ll1,}
and that \fm{\int_\R(f(0,x)-Q(0,x;\beta,y_0))\sech(\gamma(x-y_0))\ dx=0,\qquad \gamma=(1-\beta^2)^{-1/2}.}}
\begin{proof}
Everything follows directly from Lemma \ref{lemhighreginit} except the second last estimate on the difference of $f$ and the kink with a new velocity $\beta$ and a new center $y_0$. To show it, we only need to prove
\fm{\norm{(Q,Q_t)(0,\cdot;\beta_0,x_0)-(Q,Q_t)(0,\cdot;\beta,y_0)}_{H^{m+1,s}\times H^{m,s}}\lesssim\eps_0.}
This follows from  $|\beta-\beta_0|+|x_0-y_0|\lesssim\eps_0$ and part b) of Lemma \ref{lemQ} (with $k=\lceil m\rceil$) as long as $\eps_0\ll_{\beta_0}1$.
\end{proof}
\rm

\section{Proof of Theorem \ref{mthm}}\label{secpfmthm}
In this section, we prove Theorem \ref{mthm}. Fix $\beta\in(-1,1)$ and $x_0\in\R$.  Set $K(t,x)=Q(t,x;\beta,x_0)$. Let $f,\phi\in{\bf E}_{\sin}^0$ be two global solutions to the sine-Gordon equation \eqref{sg} with
\eq{\label{s3init}\norm{(f,f_t)(0)-(K,K_t)(0)}_{H^1\times L^2}+\norm{(\phi,\phi_t)(0)}_{H^1\times L^2}\lesssim\eps\ll1.}
Suppose that $f$ is the B\"acklund transform of $\phi$ by the parameter $a=\sqrt{\frac{1+\beta}{1-\beta}}$.

In Section \ref{secpfmthm:ift}, we focus on the choice of the center $x(t)$ by solving the orthogonal condition \eqref{mthmc1}. Our main result here is  Proposition \ref{proporbsta} which shows the existence and uniqueness of the $C^1$ path $x(t)$ satisfying \eqref{mthmc1}. Our proof relies on Proposition \ref{propbackt} and orbital stability of the sine-Gordon kinks. In Section \ref{s3pfmthmiiv}, we finish the proof of part i)-iv) of Theorem \ref{mthm}. The main idea is as follows. Define $\varphi,\theta$ by \fm{(\varphi,\theta)(t,x)&=(f(t,x)-Q(t,x;\beta,x(t)),f_t(t,x)-(\partial_tQ)(t,x;\beta,x(t))).}  By the B\"acklund transform \eqref{s1backeqn}, we have
\fm{\varphi_x&=\phi_t+\frac{1}{a}(\sin(\frac{\varphi+Q+\phi}{2})-\sin(\frac{Q}{2}))+a(\sin(\frac{\varphi+Q-\phi}{2})-\sin(\frac{Q}{2})),\\
\theta&=\phi_x+\frac{1}{a}(\sin(\frac{\varphi+Q+\phi}{2})-\sin(\frac{Q}{2}))-a(\sin(\frac{\varphi+Q-\phi}{2})-\sin(\frac{Q}{2})).}
From the first equation and the orthogonality condition \eqref{mthmc1}, we write down an explicit formula \eqref{s3varphiid} for $\varphi$. This formula leads to a key estimate \eqref{varphilinftykappat} for $\norm{\varphi(t)}_{L^\infty(\R)}$. One can then derive the asymptotic for $(\varphi,\varphi_x,\theta)$ using the B\"acklund transform; see \eqref{s3varphieqnpre1} and \eqref{s3varphieqnpre2}. Finally, in Section \ref{secpfmthm:final}, we prove part v) of Theorem \ref{mthm}. The main tools used there are the finite speed of propagation and the energy-momentum tensor for the sine-Gordon equation.

\subsection{The choice of the center}\label{secpfmthm:ift}

The first step in the proof is to determine the center $x(t)$. We seek to prove the following proposition.

\prop{\label{proporbsta} Fix $\beta\in(-1,1)$ and $x_0\in\R$.  Suppose the $f,\phi\in {\bf E}_{\sin}^0$ are two solutions to \eqref{sg} such that $f$ is the B\"acklund transform of $\phi$ by the parameter $a=\sqrt{\frac{1+\beta}{1-\beta}}$, that 
\fm{\norm{(f,f_t)(0,\cdot)-(Q,Q_t)(0,\cdot;\beta,x_0)}_{H^1\times L^2}+\norm{(\phi,\phi_t)(0,\cdot)}_{H^1\times L^2}\lesssim \eps\ll1.}
Then, there exists a unique $C^1$ path $x(t)$ defined for all $t\geq 0$, such that $|x(0)-x_0|\lesssim\eps$ and that  for all $t\geq 0$, we have $|x'(t)|\lesssim\eps$,
\eq{\label{orbstac1}\norm{(f,f_t)(t,\cdot)-(Q,Q_t)(t,\cdot;\beta,x(t))}_{H^1\times L^2}\lesssim \eps}and \eq{\label{orbstac2}\int_\R (f(t,x)-Q(t,x;\beta,x(t))\sech(\gamma(x-\beta t-x(t)))\ dx=0.}
Here by $Q_t(t,x;\beta,x(t))$, we mean $((\partial_tQ)(t,x;\beta,y))|_{y=x(t)}$.

Moreover, we have
\eq{\label{orbstac3}\norm{(f,f_t)(t,\cdot)-(Q,Q_t)(t,\cdot;\beta,x(t))}_{W^{1,\infty}\times L^\infty}\lesssim\norm{(\phi,\phi_t)(t,\cdot)}_{W^{1,\infty}\times L^\infty}}
provided that the right hand side of this estimate is finite and sufficiently small.}\\

\rm
We remark that the $L^2+L^\infty$ norms of $\phi,\phi_t,\phi_x$ are not used at all in this proposition.

To prove this proposition, we first construct a (not necessarily continuous) function $y(t)$ defined for all $t\geq 0$ such that $y(0)=x_0$ and that
\eq{\label{eq:osmstep}\norm{(f,f_t)(t,\cdot)-(Q,Q_t)(t,\cdot;\beta,y(t))}_{H^1\times L^2}\lesssim \eps,\qquad \forall t\geq 0.}
The existence of such a function is guaranteed by the orbital stability of the sine-Gordon kinks. We also remark that this function is not unique. Next, for each time $t\geq 0$, we find $x(t)$ such that $|x(t)-y(t)|\lesssim\eps$ and that both \eqref{orbstac1} and \eqref{orbstac2} hold. We finally prove that $x(t)$ is $C^1$, $|x'(t)|\lesssim \eps$, the path $x(t)$  is unique, and that  \eqref{orbstac3} holds.

\subsubsection{Orbital stability}
We first recall an orbital stability result from \cite[Theorem 1]{MR4242134}.
\prop[Orbital stability, \cite{MR4242134}]{\label{proporbstaold}Let $f\in {\bf E}_{\sin}^0$ be a global solution to \eqref{sg}. Suppose that for some $\beta\in(-1,1)$ and $x_0\in\R$ we have
\fm{\norm{(f,f_t)(0,\cdot)-(Q,Q_t)(0,\cdot;\beta,x_0)}_{H^1\times L^2}\leq \eps\ll1.}
Then, as long as $\eps\ll1$, for each time $t\in\R$ we have
\fm{\inf_{y\in\R}\norm{(f,f_t)(t,\cdot)-(Q,Q_t)(t,\cdot;\beta,y)}_{H^1\times L^2}\lesssim \norm{(f,f_t)(0,\cdot)-(Q,Q_t)(0,\cdot;\beta,x_0)}_{H^1\times L^2}.}}
\rmk{\rm Here we cite the orbital stability result from \cite{MR4242134} instead of those from the earlier papers \cite{MR678151,MR3059166}. This is merely for convenience. If one prefers to use the result from \cite{MR3059166}, for example, one can also prove Proposition \ref{proporbsta}. However, the conclusion in \cite{MR3059166} is slightly weaker than that of Proposition \ref{proporbstaold}: under the same assumption, the solution $f$ stays in an $\sqrt{\eps}$-neighborhood (instead of an $\eps$-neighborhood) of the family of kinks with a velocity $\beta$. Extra explanations are thus needed in our proof.

We remark that there is an alternative proof of Proposition \ref{proporbsta} without relying on the previous results on the orbital stability. Instead, we set up a continuity argument. Suppose that the $C^1$ path $x(t)$ satisfying \eqref{orbstac2} exists for $t\in[0,T]$ and that $(f,f_t)(t)$ lies in an $\sqrt{\eps}$-neighborhood of $(Q,Q_t)(t,\cdot;\beta,x(t))$ for $t\in[0,T]$. By checking that $\sqrt{\eps}$-neighborhood above can be replaced by a $C\eps$-neighborhood, we also conclude the orbital stability.}\\
\rm

As a result, for each $t\geq 0$, there exists a $y(t)\in\R$ such that \eqref{eq:osmstep} holds. This $y(t)$ might not be continuous at all, but we have the following lemma. 

\lem{\label{lem:proporbstaold}Under the assumptions of Proposition \ref{proporbstaold}, there exists a (not necessarily continuous) function $y(t)$ defined for all $t\geq 0$, such that $y(0)=x_0$ and that \eqref{eq:osmstep} holds. 

Besides, we have $|y(t_1)-y(t_2)|\lesssim\eps$ whenever $t_1,t_2\geq 0$ and $|t_1-t_2|\leq\eps\ll1$. Here the choice of $\eps$ and all the implicit constants are uniform in $t_1,t_2$.}
\begin{proof}
The existence of $y(t)$ for $t>0$ follows directly from orbital stability, and we set $y(0)=x_0$. It thus remains to prove the second half of the lemma.
For convenience, we set \fm{d(t,y):=\norm{(f,f_t)(t,\cdot)-(Q,Q_t)(t,\cdot;\beta,y)}_{H^1\times L^2},\qquad \forall t\geq0,\ y\in\R.}
Proposition \ref{proporbstaold} and part b) of Lemma \ref{lemQ} imply that $d(t,y)$ is globally Lipschitz (i.e.\ the Lipschitz constant is uniform in $(t,y)$) and finite everywhere.

Next, we claim that  $d(t,y_1)+d(t,y_2)\geq C_\gamma^{-1} \min\{|y_1-y_2|,1\}$ for some constant $C_\gamma>1$. In fact, if $|y_1-y_2|\leq 1$, 
\fm{|Q(t,x;\beta,y_2)-Q(t,x;\beta,y_1)|
&=|\int_{y_1}^{y_2}2\gamma\sech(\gamma(x-\beta t-y))\ dy|\geq C_\gamma^{-1} |y_1-y_2|\cdot 1_{|x-\beta t-y_1|\leq 1}.}
The last estimate holds because we have $|x-\beta t-y|\leq 1+|y_1-y_2|\leq 2$  whenever $y$ lies between $y_1,y_2$ and $|x-\beta t-y_1|\leq 1$. If $|y_1-y_2|\geq 1$, we notice that $y\mapsto Q(t,x;\beta,y)$ is strictly decreasing and thus \fm{|Q(t,x;\beta,y_2)-Q(t,x;\beta,y_1)|\geq |Q(t,x;\beta,y_1+\frac{y_2-y_1}{|y_2-y_1|})-Q(t,x;\beta,y_1)|\geq C_\gamma^{-1}\cdot 1_{|x-\beta t-y_1|\leq 1}.}
It follows that
\fm{d(t,y_1)+d(t,y_2)\geq\norm{Q(t,\cdot;\beta,y_2)-Q(t,\cdot;\beta,y_1)}_{L^2}\geq C_\gamma^{-1}\min\{|y_1-y_2|,1\}.} 

Now fix $t_1,t_2\geq 0$ with $|t_1-t_2|\lesssim\eps$. By the claim above, we have
\fm{\min\{|y(t_1)-y(t_2)|,1\}\lesssim_\gamma d(t_1,y(t_1))+d(t_1,y(t_2))\lesssim |t_1-t_2|+d(t_1,y(t_1))+d(t_2,y(t_2))\lesssim \eps.}
By choosing $\eps\ll1$, we forces $|y(t_1)-y(t_2)|\lesssim\eps$. Note that the choice of $\eps$ and all the implicit constants here are uniform in $t_1,t_2$.
\end{proof}\rm

\subsubsection{Construction of $x(t)$}
We recall the following result from Remark \ref{rmkexistPhi}.

\lem{\label{lemcenter}For each $(t,y_0)\in\R\times\R$, there exists a $C^1$ map $\zeta_{(\beta,t,y_0)}:\{g\in H^1:\ \norm{g}_{H^1}<r\}\to \R$ such that for $y=\zeta_{(\beta,t,y_0)}(g)$ we have $|y_0-y|\lesssim\norm{g}_{H^1}$ and 
\eq{\label{lemcenterc}\int_\R (g(x)+Q(t,x;\beta,y_0)-Q(t,x;\beta,y))\sech(\gamma (x-\beta t-y))\ dx=0.}
Here  $r\in(0,1)$ and the implicit constant in the estimate are independent of $(t,y_0)$.}
\rmk{\label{rmklemcenter}\rm Here we have a remark similar to Remark \ref{rmkbacktimp}. For each $(\beta,t,y_0)$, there exists a small constant $\wt{r}>0$ (independent of $(t,y_0)$), such that if
$\F_3(0,y-y_0,g)=0$, $|y-y_0|<\wt{r}$ and $\norm{g}_{H^1}<r$, then we must have $y=\zeta_{(\beta,t,y_0)}(g)$.}\\
\rm

Let $y(t)$ be the function constructed in Lemma \ref{lem:proporbstaold}. Now we define 
\fm{x(t)&=\zeta_{(\beta,t,y(t))}(f(t)-Q(t,\cdot;\beta,y(t))),\qquad t\geq 0.}
This is well defined because of \eqref{eq:osmstep}. It follows from Lemma \ref{lemcenter} that $|x(t)-y(t)|\lesssim \norm{f(t)-Q(t,\cdot;\beta,y(t))}_{H^1}\lesssim\eps$ and that \eqref{orbstac2} holds for all $t\geq 0$. In particular, we have $|x(0)-x_0|\lesssim\eps$. Besides, by the properties of the function $d$ defined in the proof of Lemma \ref{lem:proporbstaold}, we obtain \eqref{orbstac1}. To show  $x(t)$ is $C^1$, we fix an arbitrary $t_0\geq 0$. By the second half of Lemma \ref{lem:proporbstaold}, as long as $\eps\ll1$, we have $|y(t)-y(t_0)|\lesssim\eps$ whenever $|t-t_0|\leq\eps$ and $t\geq 0$. For these $t$, we have
\fm{&(\F_3)_{(\beta,t,y(t_0))}(0,x(t)-y(t_0),f(t)-Q(t,\cdot;\beta,y(t_0)))\\
&=\int_\R (f(t,x)-Q(t,x;\beta,x(t)))\sech(\gamma(x-\beta t-x(t)))\ dx=0.}
Using the notation from the proof of Lemma \ref{lem:proporbstaold}, we have $|x(t)-y(t_0)|+d(t,y(t_0))\lesssim |x(t)-y(t)|+d(t,y(t))+|y(t)-y(t_0)|\lesssim\eps$. By choosing $\eps\ll1$, we can apply Remark \ref{rmklemcenter} to conclude that $x(t)=\zeta_{(\beta,t,y(t_0))}(f(t)-Q(t,\cdot;\beta ,y(t_0)))$. The right side is $C^1$ at $t=t_0$.

We now prove that the $C^1$ path $x(t)$  with $|x(0)-x_0|\lesssim\eps$ and \eqref{orbstac1}, \eqref{orbstac2} is unique. Suppose that $\wt{x}(t)$ is another $C^1$ path with the same properties. For each $t\geq0$, we have
\fm{(\F_3)_{(\beta,t,y(t))}(0,\wt{x}(t)-y(t),f(t)-Q(t,\cdot;\beta,y(t)))=0.}
By Remark \ref{rmklemcenter} and \eqref{orbstac1} (with $x(t)$ replaced by $\wt{x}(t)$), we  have either $|\wt{x}(t)-x(t)|\geq \wt{r}$ for some constant $\wt{r}>0$ independent of $t$, or $x(t)=\wt{x}(t)$. However, the set  $\{ |\wt{x}(t)-x(t)|:\ t\geq 0\}$ contains $0$ (since $|\wt{x}(t)-y(0)|\lesssim\eps\ll1$) and is connected. This forces $\wt{x}(t)=x(t)$ everywhere.

Next we check $|x'(t)|\lesssim\eps$ for all $t\geq 0$. Since $x(t)$ is $C^1$, we can differentiate \eqref{orbstac2} with respect to $t$. This gives us
\eq{\label{diffthmc1}0&=\int_\R f_t(t,x)\sech(\gamma (x-\beta t-x(t)))+f(t,x)\cdot (-\beta-x'(t))\partial_x(\sech(\gamma (x-\beta t-x(t))))\ dx\\
&=\int_\R (f_t+(\beta+x'(t))f_x)(t,x)\cdot\sech(\gamma (x-\beta t-x(t)))\ dx\\&=x'(t)\int_\R f_x(t,x)\sech(\gamma(x-\beta t-x(t)))\ dx+\int_\R (f_t+\beta f_x)(t,x)\sech(\gamma(x-\beta t-x(t)))\ dx.}
By \eqref{orbstac1} and  H\"older's inequality, we have
\fm{&\int_\R f_x(t,x)\sech(\gamma(x-\beta t-x(t)))\ dx\\&=\int_\R 2\gamma\sech(\gamma(x-\beta t-x(t)))^2\ dx+O(\norm{f_x(t)-Q_x(t,\cdot;\beta,x(t)))}_{L^2})=4+O(\eps),\\
&\int_\R f_t(t,x)\sech(\gamma(x-\beta t-x(t)))\ dx\\
&=\int_\R -2\beta\gamma\sech(\gamma(x-\beta t-x(t)))^2\ dx+O(\norm{f_t(t)-Q_t(t,\cdot;\beta,x(t)))}_{L^2})=-4\beta+O(\eps).}
It follows from \eqref{diffthmc1} that $|x'(t)|\lesssim\eps$ for all $t\geq 0$.

\subsubsection{Proof of \eqref{orbstac3}} \label{secpforbstac3}Finally we check \eqref{orbstac3}. We first claim that $\norm{(\phi,\phi_t)(t)}_{H^1\times L^2}\lesssim \eps $  for all $t\geq 0$. By the energy conservation law, we have $\norm{(\sin(\phi/2),\phi_t,\phi_x)(t)}_{L^2}\lesssim\eps$. Since $|\sin(y)|\sim |y|$ whenever $|y|<1$, we can use the Sobolev embedding $H^1(\R)\subset L^\infty(\R)$ and a continuity argument to show that $|\sin(\phi/2)|\sim|\phi|$ everywhere. The proof is easy and we skip the details here. 

Now, by \eqref{orbstac1}, we have
\fm{\norm{(f(t)-Q(t,\cdot;\beta,x(t)),f_t(t)-Q_t(t,\cdot;\beta,x(t)))}_{H^1\times L^2}+\norm{(\phi,\phi_t)(t)}_{H^1\times L^2}\lesssim\eps\ll1,\quad t\geq 0.}
Since $f$ is the B\"acklund transform of $\phi$ and since \eqref{mthmc1} holds for $t\geq 0$, we have 
\fm{\F_{(\beta,t,x(t))}(0,0,\phi(t),\phi_t(t),f(t)-Q(t,\cdot;\beta,x(t)),f_t(t)-Q_t(t,\cdot;\beta,x(t)))=0,\qquad t\geq 0.}
By Remark \ref{rmkbacktimp}, as long as $\eps\ll1$ we have 
\fm{(f(t)-Q(t,\cdot;\beta,x(t)),f_t(t)-Q_t(t,\cdot;\beta,x(t)))=\Psi_{(\beta,t,x(t))}(0,0,\phi(t),\phi_t(t)),\qquad t\geq 0.}
In this identity above, the map $\Psi$ is associated to $p=2$ in Proposition \ref{propbackt}. If we apply the last part of Proposition \ref{propbackt} with $p=2$ and $q=\infty$, we notice that the left hand side is also the image of $(0,0,\phi(t),\phi_t(t))$ under the map $\Psi$ associated to $p=\infty$, as long as $\eps\ll1$. As a result, we also obtain \eqref{orbstac3} by applying \eqref{estPhiPsi} with $p=\infty$.

\subsection{Proof of part i)-iv) in Theorem \ref{mthm}}\label{s3pfmthmiiv}
By Proposition \ref{proporbsta}, we have proved the existence of the $C^1$ path $x(t)$ such that $|x(0)-x_0|+|x'(t)|\lesssim\eps$ and  \eqref{mthmc1} hold. We have also proved \eqref{mthmc11} and \eqref{mthmc3} in Theorem \ref{mthm}. Recall that \eqref{mthmc11} states that
\fm{\norm{(f(t)-\wt{K}(t),f_t(t)-\wt{K}_0(t))}_{H^1\times L^2}+\norm{(\phi(t),\phi_t(t))}_{H^1\times L^2}\lesssim\eps,\qquad \forall t\geq 0}
where $(\wt{K},\wt{K}_0)(t,x)=(Q,Q_t)(t,x;\beta,x(t))$.

Define $\varphi,\theta$ by \eq{\label{defnvarphitheta}(\varphi,\theta)(t,x)&=(f(t,x)-Q(t,x;\beta,x(t)),f_t(t,x)-(\partial_tQ)(t,x;\beta,x(t))).}  With an abuse of notation, we write $Q(t,x;\beta,x(t))$ as $Q$. By \eqref{backeqn} and Lemma \ref{lemQ}, we have
\eq{\label{varphixtheta:eq}\varphi_x&=\phi_t+\frac{1}{a}(\sin(\frac{\varphi+Q+\phi}{2})-\sin(\frac{Q}{2}))+a(\sin(\frac{\varphi+Q-\phi}{2})-\sin(\frac{Q}{2})),\\
\theta&=\phi_x+\frac{1}{a}(\sin(\frac{\varphi+Q+\phi}{2})-\sin(\frac{Q}{2}))-a(\sin(\frac{\varphi+Q-\phi}{2})-\sin(\frac{Q}{2})).}
By the Taylor expansion of the sine function, one has
\eq{\label{s3varphieqnpre}\varphi_x&=\phi_t+\gamma\cos(\frac{Q}{2})\varphi-\beta\gamma\cos(\frac{Q}{2})\phi+O(|\varphi|^2+|\phi|^2)=\gamma\cos(\frac{Q}{2})\varphi+O(|\phi_t|+\eps|\varphi|+|\phi|),\\
\theta&=\phi_x-\beta\gamma\cos(\frac{Q}{2})\varphi+\gamma\cos(\frac{Q}{2})\phi+O(|\varphi|^2+|\phi|^2)=O(|\varphi|+|\phi|+|\phi_x|).}
Here we use the bound $\norm{(\varphi,\phi)}_{L^\infty}\lesssim \norm{(\varphi,\phi)}_{H^1}\lesssim\eps$.

Now we rewrite the first equation in \eqref{s3varphieqnpre} as  $\varphi_x=\gamma\cos(\frac{Q}{2})\varphi+F$ where $F=F(t,x)=O(|\phi_t|+\eps|\varphi|+|\phi|)$. Note that this ODE for $\varphi$ can be viewed as a linearization of the first equation of the B\"acklund transform \eqref{backeqn} around a kink. We remind our readers of the proof of Lemma \ref{lemexistPhi} where a similar linearized equation, \eqref{leminitsoleqn}, has been studied. By solving the ODE, we obtain an explicit formula of $\varphi$:
\eq{\label{s3varphiid0}\varphi(t,x)&=(\I F)(t,x)+\frac{\varphi(t,\beta t+x(t))}{\cosh(\gamma(x-\beta t-x(t)))}}
where
\eq{\label{s3:defn:mclI}(\I F)(t,x)&=\int_{\beta t+x(t)}^x\frac{\cosh(\gamma(y-\beta t-x(t)))}{\cosh(\gamma(x-\beta t-x(t)))}F(t,y)\ dy.}
Moreover, by \eqref{mthmc1} we have $\int_\R\varphi(t,x)\sech(\gamma(x-\beta t-x(t)))\ dx=0$, so 
\fm{0&=\int_\R(\I F)(t,x)\sech(\gamma(x-\beta t-x(t)))\ dx+\int_\R\frac{\varphi(t,\beta t+x(t))}{\cosh(\gamma(x-\beta t-x(t)))^2}\ dx}
and thus
\fm{\varphi(t,\beta t+x(t))&=-\frac{\gamma}{2}\int_\R(\I F)(t,x)\sech(\gamma(x-\beta t-x(t)))\ dx.}
In summary, we have
\eq{\label{s3varphiid}\varphi(t,x)=\I F(t,x)-\frac{\gamma}{2\cosh(\gamma(x-\beta t-x(t)))}\int_\R (\mcl{I}F)(t,y)\sech(\gamma(y-\beta t-x(t)))\ dy.}

We next state a lemma for $\I$. Though we have set $F=\varphi_x-\gamma\cos(Q/2)\varphi$ above,  the lemma holds for an arbitrary function $F=F(t,x)$. For convenience, we use $J_{t,x}$ to denote the closed interval with endpoints $x$ and $\beta t+x(t)$. 
\lem{\label{lemiest}For any function $F=F(t,x)$, we have
\fm{&|\I F(t,x)|\leq2\int_{J_{t,x}}e^{-\gamma|x-y|}|F(t,y)|\ dy,\\
&\int_\R|\I F(t,x)|\sech(\gamma(x-\beta t-x(t)))\ dx\lesssim_\gamma (e^{-\gamma|\cdot|}*|F(t)|)(\beta t+x(t)).}
The implicit constants depend only on $\gamma$.}
\begin{proof}
Fix an arbitrary $y\in J_{t,x}$. We have $\sgn(y-\beta t-x(t))=\sgn(x-\beta t-x(t))$ unless $y=\beta t+x(t)$. We also have $\sgn(x-y)=\sgn(x-\beta t-x(t))$ unless $x=y$. Since $e^{|z|}/2\leq\cosh(z)\leq e^{|z|}$, we have
\eq{\label{lemiestf}0<\frac{\cosh(\gamma(y-\beta t-x(t)))}{\cosh(\gamma(x-\beta t-x(t)))}\leq 2 e^{\gamma(|y-\beta t-x(t)|-|x-\beta t-x(t)|)}=2e^{-\gamma|x-y|}.}
As a result, we have
\fm{|\mcl{I}F(t,x)|\leq \int_{J_{t,x}}2e^{-\gamma|x-y|}|F(t,y)|\ dy.}

Next, by \eqref{lemiestf}, we have
\fm{
&\int_\R|\I F(t,x)|\sech(\gamma(x-\beta t-x(t)))\ dx\\
&\leq 2\int_\R\int_{J_{t,x}} |F(t,y)|\sech(\gamma(x-\beta t-x(t)))e^{-\gamma|x-y|}\  dydx\\
&\leq 4\int_\R\int_{J_{t,x}}|F(t,y)|\sech(\gamma(y-\beta t-x(t)))e^{-2\gamma|x-y|}\  dydx\\
&\lesssim_\gamma\int_\R |F(t,y)|e^{-\gamma|y-\beta t-x(t)|} dy=(e^{-\gamma|\cdot|}*|F(t)|)(\beta t+x(t)).}
\end{proof}\rm

Now we apply Lemma \ref{lemiest} to \eqref{s3varphiid}. For all $x\in\R$ and $t\geq 0$, we have
\fm{|\varphi(t,x)|&\lesssim_\gamma\int_{J_{t,x}}e^{-\gamma|x-y|}(|\phi_t|+|\phi|+\eps|\varphi|)(t,y)\ dy\\
&\quad+\sech(\gamma(x-\beta t-x(t)))\int_\R e^{-\gamma|y-\beta t-x(t)|}(|\phi_t|+|\phi|+\eps|\varphi|)(t,y)\ dy.}
Since $\sup_{x\in\R}\norm{e^{-\gamma|\cdot-x|}}_{L^1\cap L^2(\R)}\lesssim_\gamma1$, we apply H\"older's inequality to obtain
\fm{\sup_{x\in\R}\int_{\R}e^{-\gamma|x-y|}(|\phi_t|+|\phi|)(t,y)\ dy
&\lesssim_\gamma \norm{(\phi,\phi_t)(t)}_{(L^2+L^\infty)(\R)},\\
\sup_{x\in\R}\int_{\R}e^{-\gamma|x-y|}|\varphi(t,y)|\ dy
&\lesssim_\gamma \norm{\varphi(t)}_{L^\infty(\R)}.}
Therefore, we have
\fm{\norm{\varphi(t)}_{L^\infty(\R)}&\lesssim_\gamma \eps\norm{\varphi(t)}_{L^\infty(\R)}+\norm{(\phi,\phi_t)(t)}_{(L^2+L^\infty)(\R)}.}
The implicit constants depend only on $\gamma$. 
By choosing $\eps\ll_\gamma1$, we can use the left hand side to absorb the term $\eps\norm{\varphi(t)}_{L^\infty(\R)}$. In conclusion, we have
\eq{\label{varphilinftykappat}\norm{\varphi(t)}_{L^\infty(\R)}&\lesssim_\gamma \norm{(\phi,\phi_t)(t)}_{(L^2+L^\infty)(\R)},\qquad \forall t\geq 0. }
Then, by \eqref{s3varphieqnpre}, we have
\fm{|\varphi_x|+|\theta|&\lesssim |\varphi|+|\phi|+|\phi_t|+|\phi_x|.}
By the triangle inequality, we have
\eq{\label{varphilinftykappat2}\norm{(\varphi_x,\theta)(t)}_{(L^2+L^\infty)(\R)}&\leq \norm{\varphi(t)}_{L^\infty(\R)}+\norm{(\phi,\phi_t,\phi_x)(t)}_{(L^2+L^\infty)(\R)}\lesssim_\gamma\norm{(\phi,\phi_t,\phi_x)(t)}_{(L^2+L^\infty)(\R)}.}
We thus obtain \eqref{mthmc2}.

We now derive asymptotics for $(\varphi,\varphi_x,\theta)$. By \eqref{s3varphieqnpre}, for all $(t,x)$ we have
\fm{F=\varphi_x-\gamma\cos(\frac{Q}{2})\varphi=\phi_t-\beta\gamma\cos(\frac{Q}{2})\phi+O(|\varphi|^2+|\phi|^2).}
By \eqref{varphilinftykappat} and Lemma \ref{lemiest}, we have
\eq{\label{s3varphieqnpre1}\varphi(t,x)&=\I F(t,x)-\frac{\gamma}{2\cosh(\gamma(x-\beta t-x(t)))}\int_\R (\mcl{I}F)(t,y)\sech(\gamma(y-\beta t-x(t)))\ dy\\
&=\I(\phi_t-\beta\gamma\cos(\frac{Q}{2})\phi)(t,x)\\&\quad-\frac{\gamma}{2\cosh(\gamma(x-\beta t-x(t)))}\int_\R \mcl{I}(\phi_t-\beta\gamma\cos(\frac{Q}{2})\phi)(t,y)\sech(\gamma(y-\beta t-x(t)))\ dy\\
&\quad+O(\norm{(\phi,\phi_t)(t)}^2_{(L^2+L^\infty)(\R)}).}
By plugging this into \eqref{s3varphieqnpre}, we obtain
\eq{\label{s3varphieqnpre2}\varphi_x&=\phi_t-\beta\gamma\cos(\frac{Q}{2})\phi+\gamma\cos(\frac{Q}{2})\varphi+O(|\varphi|^2+|\phi|^2),\\
\theta&=\phi_x+\gamma\cos(\frac{Q}{2})\phi-\beta\gamma\cos(\frac{Q}{2})\varphi+O(|\varphi|^2+|\phi|^2).}
By \eqref{varphilinftykappat}, we have \fm{|\varphi|^2+|\phi|^2&\lesssim\norm{\varphi}_{L^\infty(\R)}^2+\norm{\phi}^2_{L^\infty(\R)}\lesssim \norm{\phi}^2_{L^\infty(\R)}+\norm{\phi_t}^2_{(L^2+L^\infty)(\R)}.}
We thus obtain \eqref{mthmc4}.

Finally, by \eqref{diffthmc1} and H\"older's inequality, we obtain
\fm{|x'(t)|&\lesssim[|(f_t+\beta f_x)(t)|*\sech(\gamma(\cdot))](\beta t+x(t))
\lesssim[(|\varphi_x(t)|+|\theta(t)|)*e^{-\gamma|\cdot|}](\beta t+x(t))\\
&\lesssim \norm{(\varphi_x,\theta)(t)}_{(L^2+L^\infty)(\R)}\norm{e^{-\gamma|\cdot|}}_{(L^1\cap L^2)(\R)}\lesssim_\gamma\norm{(\phi,\phi_t,\phi_x)(t)}_{(L^2+L^\infty)(\R)} .}
In the second estimate, we use 
\fm{f_t+\beta f_x&=\theta+Q_t(t,x;\beta,x(t))+\beta \varphi_x+\beta Q_x(t,x;\beta,x(t))=\theta+\beta\varphi_x.}
See part a) of Lemma \ref{lemQ}.

\subsection{Proof of part v) of Theorem \ref{mthm}}\label{secpfmthm:final}
It remains to prove \eqref{mthmc21}. Fix $R\in\R$ and $s>1/2$. Suppose that
\eq{\label{cs0defn}1\leq C_{s,0}:=1+\norm{(f,f_t)(0)-(Q,Q_t)(0,\cdot;\beta,x_0)}_{H^{1,s}\times L^{2,s}}<\infty.}
By part c) in Lemma \ref{lemQ}, we have \fm{\norm{\lra{\cdot}^s(\sin(f/2),f_x,f_t)(0)}_{L^2}\lesssim C_{s,0}.} Because of  \eqref{s3varphieqnpre}, it suffices to prove that
\fm{|f(t,x)-Q(t,x;\beta,x(t))|\lesssim_{R,s,C_{s,0}}\min\{t^{-1/4}\lra{|x|-t}^{-1/4},\lra{|x|-t}^{-s}\}}whenever $t\gtrsim_{R,s}1$ and $|x|\geq t+R$. Moreover, since $|x(t)-x_0|+|x'(t)|\lesssim \eps$, we have \fm{|\beta t+x(t)|\leq |\beta t|+|x(t)|\leq (|\beta|+C\eps)t+C\leq \frac{|\beta|+1}{2}(t-|R|)\leq \frac{|\beta|+1}{2}|x|}
whenever $t\gtrsim_R1$ and $|x|\geq t+R$. Set $c_0(x)=2\pi\cdot 1_{x> 0}$. It follows that
\fm{ | Q(t,x;\beta,x(t))-c_0(x)|= \int_{|x-\beta t-x(t)|}^\infty 2\gamma \sech(\gamma y)\ dy\lesssim  \int_{C^{-1}|x|}^\infty e^{-\gamma |y|}\ dy\lesssim e^{-C^{-1}|x|}}
for $|x|\geq t+R$. Since $e^{-C^{-1}|x|}=e^{-C^{-1}(|x|-t)}\cdot e^{-C^{-1}t}\lesssim \lra{|x|-t}^{-s}t^{-1/4}$, it is remains to check 
\eq{\label{goalout}|f(t,x)-c_0(x)|\lesssim_{R,s,C_{s,0}}\min\{t^{-1/4}\lra{|x|-t}^{-1/4},\lra{|x|-t}^{-s}\}.}

Let us first prove that $|f(t,x)-c_0(x)|\lesssim_{R,s} \lra{|x|-t}^{-s}C_{s,0}$ whenever $|x|\geq t+R$. 

\lem{We have $|f(t,x)-c_0(x)|\lesssim_{R,s} \lra{|x|-t}^{-s}C_{s,0}$ whenever $t\gtrsim_{R,s}1$ and $|x|\geq t+R$.}
\begin{proof}
Recall from Proposition \ref{proporbsta} that $\norm{f(t)-Q(t,\cdot;\beta,x(t))}_{L^\infty}\lesssim \eps$ for all $t\geq 0$. Thus, if $t+R\leq |x|\leq t+|R|+1$, we already have $|f(t,x)-c_0(x)|\lesssim 1\lesssim \lra{|x|-t}^{-s}$. 

Fix $(t_0,x_0)\in[0,\infty)\times\R$ such that $t_0\gtrsim_{R,s}1$ and $|x_0|\geq t_0+|R|+1$. Fix a cutoff function $\chi\in C^\infty(\R)$ such that $0\leq \chi\leq 1$, $\chi(s)=0$ whenever $s\leq -1/2$, and $\chi(s)=1$ whenever $s\geq 0$. We now set 
\fm{(g^0(x),g^1(x))&=\chi(|x|-|x_0|+t_0)\cdot (f(0,x),f_t(0,x)).}
Since $\supp(g^0,g^1)\subset \{|x|\geq |x_0|-t_0-1/2\}$ and since $\chi,\chi'=O(1)$ where the implicit constant is independent of $(t_0,x_0)$, we have
\fm{\norm{(\sin(g^0/2),(g^0)',g^1)}_{L^{2}}&\lesssim \lra{|x_0|-t_0}^{-s}\norm{(\sin(f/2),f_x,f_t)(0)}_{L^{2,s}}\lesssim \lra{|x_0|-t_0}^{-s}C_{s,0}.} Let $g$ be the global solution to \eqref{sg} with initial data $(g^0,g^1)$. The previous claim shows that $(g^0,g^1)$ is in the energy space, so such a $g$ exists and belongs to ${\bf E}^0_{\sin}$. By energy conservation, we have
\eq{\label{sinenergylem}\norm{(\sin(g/2),g_x,g_t)(t_0)}_{L^2}\lesssim \lra{|x_0|-t_0}^{-s}C_{s,0}.}
By the finite speed of propagation, we also have $g\equiv f$ whenever $|x|\geq t+|x_0|-t_0-1/4$.

We hope to apply the Sobolev embedding and \eqref{sinenergylem}. To achieve this goal, we define \fm{\wt{g}(x)=\chi(\frac{1}{2}(|x|-|x_0|))(g(t_0,x)-c_0(x)).}
Note that $\supp\wt{g}\subset\{|x|\geq |x_0|-1/4\}$ and $f(t_0)-c_0\equiv g(t_0)-c_0\equiv \wt{g}$ for $|x|\geq |x_0|$. By the choice of $\chi$, we have
\fm{|\wt{g}|+|\wt{g}'|&\lesssim|g(t_0)-c_0|+|g_x(t_0,x)|}
where the implicit constant is independent of $(t_0,x_0)$. We now claim that $|g(t_0)-c_0|\sim|\sin(g(t_0)/2)|$ whenever $|x|\geq |x_0|-1/4$. If this claim holds, then by \eqref{sinenergylem} we conclude that $|\wt{g}(x_0)|\lesssim \norm{\wt{g}}_{H^1}\lesssim \lra{|x_0|-t_0}^{-s}C_{s,0}$ immediately. 

Let us prove the claim. Fix any $x$ with $|x|\geq |x_0|-1/4$. Note that 
\fm{|Q(t_0,x;\beta,x(t_0))-c_0(x)|\leq \int_{|x_0|-1/4}^\infty C\gamma e^{-\gamma |y|}\ dy\leq Ce^{-\gamma(|x_0|-1/4)}.}
Then, whenever $|x|\geq |x_0|-1/4$, we have
\fm{|g(t_0,x)-c_0(x)|&=|f(t_0,x)-c_0(x)|\lesssim \eps+e^{-\gamma(|x_0|-1/4)}\lesssim  \eps+e^{-\gamma(t_0+|R|+3/4)}.}
By choosing $t_0\gtrsim_{R,s}1$ and $\eps\ll_{R,s}1$, we obtain $|g(t_0,x)-c_0(x)|=|f(t_0,x)-c_0(x)|\leq \pi/4$, in which case we have $\sin(g(t_0,x)/2)=|\sin((g(t_0,x)-c_0(x))/2)|\sim |g(t_0,x)-c_0(x)|$.
\end{proof}\rm

The estimate \eqref{goalout} now follows from the following lemma. The proof is essentially the same as that in \cite[Section 4.1]{MR4242133}.

\lem{\label{lem4.7weiyang}We have $|f(t,x)-c_0(x)|\lesssim_{R,s} C_{s,0}t^{-1/4}\lra{|x|-t}^{-1/4}$ whenever $t\gtrsim_{R,s}1$ and $|x|\geq t+R$.}
\begin{proof}

Define the energy-momentum tensor $T_{**}$ by
\fm{T_{00}&=\frac{1}{2}(f_t^2+f_x^2)+1-\cos f,\qquad T_{10}=T_{01}=f_tf_x,\qquad T_{11}=\frac{1}{2}(f_t^2+f_x^2)-1+\cos f.}
For each $R\in\R$, we have
\fm{\partial_tT_{00}=\partial_xT_{10},\quad \partial_tT_{01}=\partial_xT_{11},\quad \partial_t((x+|R|+1)T_{00}+tT_{01})=\partial_x((x+|R|+1)T_{01}+tT_{11}).}
For all $R<a<b$, we integrate the last identity above in the region
\fm{U=\{(t,x):\ t>0,\ x-t>a,x+t<b\}.}
By  Green's theorem, we have
\fm{\int_{\partial  U}((x+|R|+1)T_{00}+tT_{01})\nu^0-((x+|R|+1)T_{01}+tT_{11})\nu^1\ dS=0}
where $(\nu^0,\nu^1)$ is the outer unit normal of $\partial U$. The integral on $\partial U\cap\{t=0\}$ equals 
\fm{\int_a^b(x+|R|+1) T_{00}(0,x)\ dx.}
The assumption \eqref{mthmc21asu} implies that this quantity is $O_R(1)<\infty$. The integral on $\partial U\cap\{x-t=a\}$ equals 
\fm{&\int_0^{(b-a)/2}((x+|R|+1) T_{00}+(x+t+|R|+1)T_{01}+tT_{11})|_{x-t=a}\ dt\\
&=\int_0^{(b-a)/2}\frac{1}{2}(2t+a+|R|+1)|(f_t+f_x)(t,t+a)|^2+(a+|R|+1)(1-\cos f(t,t+a))\ dt.}
The integral on $\partial U\cap\{x+t=b\}$ equals 
\fm{&\int_0^{(b-a)/2}((x+|R|+1) T_{00}-(x-t+|R|+1)T_{01}-tT_{11})|_{x+t=b}\ dt\\
&=\int_0^{(b-a)/2}\frac{1}{2}(b+|R|+1-2t)|(f_t-f_x)(t,b-t)|^2+(b+|R|+1)(1-\cos f(t,b-t))\ dt.}
Each of these three integrals is nonnegative.
Taking the orientation of the integral into account and sending $b\to\infty$, we conclude that for all $a>-R$,
\eq{\label{weiyangkeyest}&\int_0^{\infty}\frac{1}{2}(2t+a+|R|+1)|(f_t+f_x)(t,t+a)|^2+(a+|R|+1)(1-\cos f(t,t+a))\ dt\\
&\leq \int_{R}^\infty(x+|R|+1)T_{00}(0,x)\ dx\lesssim_R C_{s,0}^2.}
Note that the implicit constant is uniform in $a$.
 
We now recall the following lemma which is \cite[Lemma 2.1]{MR4242133}. While the authors there assumed $p>1$ in \cite{MR4242133}, the next lemma holds for $p=1$ (so we have $p+3=4$ in the power on the left side and $p+1=2$ in the power on the right side).
\lem[Wei-Yang, \cite{MR4242133}]{Fix $a_1\geq 0$ and $a_2\geq 1$. For each function $g(t)\in H^1_{\rm loc}(a_1,\infty)$, we have for all $t\geq a_1$, 
\fm{(t+a_2)|g(t)|^4\lesssim \int_{a_1}^\infty|g(s)|^2\ ds\cdot \int_{a_1}^\infty(s+a_2)|g'(s)|^2\ ds.}}\rm

In this lemma, we set $g(t)=f(t,t+a)-2\pi$  and $a_2=a+|R|+1$. The value of $a_1$ will be chosen later. Then \eqref{weiyangkeyest} yields
\fm{\int_{a_1}^\infty (s+a+|R|+1)|g'(s)|^2\ ds=\int_{a_1}^\infty (t+a+|R|+1)|(f_t+f_x)(t,t+a)|^2\ dt\lesssim_R C_{s,0}^2,\\
\int_{a_1}^\infty (1-\cos g(s))\ ds=\int_{a_1}^\infty |\sin (f(t,t+a)/2)|^2\ dt\lesssim (a+|R|+1)^{-1}C_{s,0}^2.}
The implicit constants are uniform in $a$ and $a_1$. Recall that we have proved $|Q(t,x;\beta,x(t))-2\pi|\lesssim e^{-C^{-1}x}$ for $x\geq t+R$. It follows that
\fm{|g(t)|\leq |f(t,t+a)-Q(t,t+a;\beta,x(t))|+|2\pi-Q(t,t+a;\beta,x(t))|\lesssim \eps+e^{C^{-1}t}.}
All the implicit constants  are uniform in $a>R$, so we can choose $a_1\gg 1$  and $\eps\ll1$ (both depending on $R$, not on $a$) to make $|g(t)|$ small enough so that  $1-\cos(g(t))=2\sin^2(g(t)/2)\sim |g(t)|^2$ for all $t\geq a_1$. In summary, we have 
\fm{|g(t)|^4=|f(t,t+a)-2\pi|^4\lesssim (t+a+|R|+1)^{-1}(a+|R|+1)^{-1}C_{s,0}^4,\qquad \forall t\geq a_1,}
or equivalently,
\fm{|f(t,x)-2\pi|\lesssim t^{-1/4}\lra{|x|-t}^{-1/4}C_{s,0},\qquad \forall t\gtrsim_{R,s}1,\ x>t+R.}
The proof for $f$ when $x<-t-R$ is essentially the same.  \end{proof}\rm

\section{Proof of Theorem \ref{thmasysta}}\label{secpfasysta}

\subsection{The asymptotic stability}\label{spfthmlrc}
We first prove part (a) of Theorem \ref{thmasysta}.  Fix $s>1/2$, $\beta_0\in(-1,1)$ and $x_0\in\R$. We assume that $f\in{\bf E}^0_{\sin}$ which is a solution to \eqref{sg} such that
\fm{\norm{(f,f_t)(0)-(K,K_t)(0)}_{H^{1,s}\times L^{2,s}}\leq \eps\ll1}
where $K(t,x)=Q(t,x;\beta_0,x_0)$ is defined by \eqref{kink}. By Corollary \ref{corsolvebackthigh}, we obtain a solution $\phi\in{\bf E}^0_{\sin}$ to \eqref{sg}, a new velocity $\beta\in(-1,1)$ and a new center $y_0\in\R$. Moreover,  $f$ is the B\"acklund transform of $\phi$ by a new parameter $a=\sqrt{\frac{1+\beta}{1-\beta}}$, and we have
\fm{|\beta-\beta_0|+|x_0-y_0|+\norm{(\phi,\phi_t)(0)}_{H^{1,s}\times L^{2,s}}+\norm{(f,f_t)(0)-(Q,Q_t)(0,\cdot;\beta,y_0)}_{H^{1,s}\times L^{2,s}}\lesssim\eps.}

We now apply Theorem \ref{mthm} (with $x_0$ replaced by $y_0$). It remains to check \eq{\label{ineq:phil2linfty}\lim_{t\to\infty}\norm{(\phi,\phi_t,\phi_x)}_{(L^2+L^\infty)(\R)}=0} and \eq{\label{ineq:xtr}\lim_{t\to\infty}\norm{(\phi,\phi_t,\phi_x)(t)}_{L^2(\{x\in\R:\ |x|\geq t+R\})}=0,\qquad \forall R\in\R.}
Note that \eqref{ineq:xtr}, \eqref{mthmc21},  and the estimate (here we use $s>1/2$)
\fm{&\int_{|x|\geq t+R}\min\{t^{-1/2}\lra{|x|-t}^{-1/2},\lra{|x|-t}^{-2s}\}\ dx\\
&\leq \int_{|x|-t\in[R,t^{1/4}]}t^{-1/2}\ dx+\int_{|x|\geq t+t^{1/4}}\lra{|x|-t}^{-2s}\ dx\lesssim t^{-1/4}+t^{(1-2s)/4}\to 0,\qquad t\to\infty}
imply part iv) in Theorem \ref{thmasysta}, part (a).

By H\"older's inequality, we have
\fm{\norm{(\phi,\phi_t,\phi_x)(t)}_{L^2(\{x\in\R:\ |x|-t\in [R_1,R_2]\})}&\leq \norm{(\phi,\phi_t,\phi_x)(t)}_{L^2+L^\infty}\norm{1_{|x|\in[t+R_1,t+R_2]}}_{L^1\cap L^2}\\
&\lesssim_{R_2-R_1}\norm{(\phi,\phi_t,\phi_x)(t)}_{L^2+L^\infty }}
for all constants $R_1<R_2$.
It thus suffices to prove \eqref{ineq:phil2linfty} and \eqref{ineq:xtr} with $R=0$.
We end the proof by applying the next proposition.

\prop{\label{thmlowdata}Fix $s>1/2$. Let $\phi$ be a global solution to the sine-Gordon equation. Then, whenever  $\norm{(\phi,\phi_t)(0)}_{H^{1,s}\times L^{2,s}}\ll_s1$, we have
\eq{\lim_{t\to\infty}\kh{\norm{\phi(t)}_{L^\infty(\R)}+\norm{(\phi_t,\phi_x)(t)}_{(L^2+L^\infty)(\R)}+\norm{(\phi,\phi_t,\phi_x)(t)}_{L^2(\{x\in\R:\ |x|\geq t\})}}=0.}}\rm

We will prove it in Appendix \ref{seclow}. In its proof, we make use of several results from \cite{chen2020long,MR1697487}. We will also need some high order conservation laws for the sine-Gordon equation.

\subsection{The asymptotics of the difference}\label{spfthmhrc}
We move on to the proof of part (b). Fix $\beta_0\in(-1,1)$ and $x_0\in\R$. We assume that $f\in{\bf E}^0_{\sin}$ which is a solution to \eqref{sg} such that
\fm{\norm{(f,f_t)(0)-(K,K_t)(0)}_{H^{m+1,1}\times H^{m,1}}\leq \eps\ll1,\qquad m>3/2}
where $K(t,x)=Q(t,x;\beta_0,x_0)$ is defined by \eqref{kink}. By Corollary \ref{corsolvebackthigh}, we obtain a solution $\phi\in{\bf E}^0_{\sin}$ to \eqref{sg}, a new velocity $\beta\in(-1,1)$ and a new center $y_0$. Moreover, $f$ is the B\"acklund transform of $\phi$ by a new parameter $a=\sqrt{\frac{1+\beta}{1-\beta}}$, and we have
\fm{&|\beta-\beta_0|+|x_0-y_0|+\norm{(\phi,\phi_t)(0)}_{H^{m+1,1}\times H^{m,1}}\\
&\qquad+\norm{(f,f_t)(0)-(Q,Q_t)(0,\cdot;\beta,y_0)}_{H^{m+1,1}\times H^{m,1}}\lesssim\eps
.}
 
We now state a proposition on the global bounds and asymptotics for the sine-Gordon equation with sufficiently small, smooth, and localized data.  
\prop{\label{thmwpt}Consider the Cauchy problem \eqref{sg} with initial data $(\phi_0,\phi_1)$ such that \fm{\norm{(\phi,\phi_t)(0)}_{H^{m+1,1}\times H^{m,1}}\leq\eps,\qquad m>3/2.}
Then, for $\eps\ll_m1$, there exists a unique global $C^1$ solution $\phi$ to this Cauchy problem satisfying \fm{\norm{\phi(t)}_{L^\infty}+\norm{\phi_t(t)}_{L^\infty}+\norm{\phi_x(t)}_{L^\infty}\lesssim\eps\lra{t}^{-1/2},\qquad \forall t\geq 0.}
Moreover, there exists a complex-valued function $W=W(\xi)$ defined on $\R$ and two small constants $\delta,\kappa\in(0,1)$, such that for each $l\in[0,1+\delta]$ and all $(t,x)\in[1,\infty)\times\R$, we have
\eq{\label{thmwpt:c}&(\lra{D}^l\phi+i\lra{D}^{l-1}\phi_t)(t,x)\\&=t^{-1/2} \lra{x/\rho}^lW(x/\rho)\exp(-i\rho+\frac{i}{32\lra{x/\rho}}|W(x/\rho)|^2\ln t)\cdot  1_{|x|<t}+O_\kappa(\eps t^{-1/2-\kappa}).}
Here $\rho=\sqrt{t^2-x^2}$. Moreover, we have $|W(\xi)|\lesssim\eps \lra{\xi}^{-1-2\delta}$.
}\rm
\bigskip

In Appendix \ref{secthmwpt}, we present a sketch of proof of this result. We use the method of testing by wave packets which was first introduced in the context of the cubic nonlinear Schr\"odinger flow (NLS) by Ifrim and Tataru \cite{MR3382579}. We also refer to  \cite{MR4483135,MR3499085,MR3667289,MR3462131,MR3816658,MR4019190,MR3813999} for other applications of this method. In our proof, we follow a recent expository note \cite{MR4693097}. See Appendix \ref{secthmwpt} for more details. We also refer to \cite{MR2457221} where a similar result was proved with a different method.

We now apply Theorem \ref{mthm}. It remains to prove \eqref{thmasysta:asyc} and \eqref{thmasysta:asyc2}. Recall the definitions of $\wt{K}$ and $\wt{K}_0$ in Theorem \ref{thmasysta}. By \eqref{thmwpt:c}, for all $(t,x)\in[1,\infty)\times\R$, we have
\fm{(\phi_t-\beta\gamma\cos(\frac{\wt{K}}{2})\phi)(t,x)&=\mcl{A}_W(t,x;\beta,x(t))+O_\kappa(\eps t^{-1/2-\kappa}),\\
(\phi_x+\gamma\cos(\frac{\wt{K}}{2})\phi)(t,x)&=\mcl{B}_W(t,x;\beta,x(t))+O_\kappa(\eps t^{-1/2-\kappa})}
with
\fm{&\mcl{A}_W(t,x;\beta,x(t))\\
&=t^{-1/2}1_{|x|<t}\cdot \Real\kh{(-i\lra{x/\rho}-\beta\gamma\cos(\frac{\wt{K}}{2})) W(x/\rho)\exp(-i\rho+\frac{i}{32\lra{x/\rho}}|W(x/\rho)|^2\ln t)},\\
&\mcl{B}_W(t,x;\beta,x(t))\\
&=t^{-1/2}1_{|x|<t}\cdot \Real\kh{(ix/\rho+\gamma\cos(\frac{\wt{K}}{2})) W(x/\rho)\exp(-i\rho+\frac{i}{32\lra{x/\rho}}|W(x/\rho)|^2\ln t)}.}
Here $\kappa\in(0,1)$ is a sufficiently small constant and $\rho=\sqrt{t^2-x^2}$. By the second formula in \eqref{mthmc4}, we obtain \eqref{thmasysta:asyc2}.

Next, by Lemma \ref{lemiest}, we have
\fm{\Phi(t,x;\beta,x(t))&=\int_{\beta t+x(t)}^x \frac{\cosh(\gamma(y-\beta t-x(t)))}{\cosh(\gamma(x-\beta t-x(t)))}(\phi_t-\beta\gamma\cos(\frac{\wt{K}}{2})\phi)(t,y)\ dy
\\&=\int_{\beta t+x(t)}^x \frac{\cosh(\gamma(y-\beta t-x(t)))}{\cosh(\gamma(x-\beta t-x(t)))}\mcl{A}_W(t,y;\beta,x(t))\ dy+O_\kappa(\eps t^{-1/2-\kappa}).}
This $\Phi$ comes from part iv), Theorem \ref{mthm}. Next, we have
\fm{&\int_{\R}\Phi(t,x;\beta,x(t))\sech(\gamma(x-\beta t-x(t)))\ dx\\
&=\int_{\R}\int_{\beta t+x(t)}^x\frac{\cosh(\gamma(y-\beta t-x(t)))}{\cosh(\gamma(x-\beta t-x(t)))^2}\mcl{A}_W(t,y;\beta,x(t))\ dy dx+O_\kappa(\eps t^{-1/2-\kappa})\\
&=\int_{\R}\int_{y}^{\sgn(y-\beta t-x(t))\cdot\infty}\frac{\cosh(\gamma(y-\beta t-x(t)))}{\cosh(\gamma(x-\beta t-x(t)))^2}\mcl{A}_W(t,y;\beta,x(t))\ dxdy +O_\kappa(\eps t^{-1/2-\kappa})\\
&=\gamma^{-1}\int_{\R}\sgn(y-\beta t-x(t))\cdot e^{-\gamma|y-\beta t-x(t)|}\mcl{A}_W(t,y;\beta,x(t))\ dy +O_\kappa(\eps t^{-1/2-\kappa}).}
In the second step, we interchange the order of integration. The proof of the second estimate in Lemma \ref{lemiest} can be used to explain why one can do so. To get the last step, we notice that
\fm{&\int_y^{\pm\infty}\sech(\gamma(x-\beta t-x(t)))^2\ dx=\gamma^{-1}\tanh(\gamma(x-\beta t-x(t)))]_{y}^{\pm\infty}\\
&=\gamma^{-1}(\pm1-\tanh(\gamma(y-\beta t-x(t))))=\pm\gamma^{-1}\sech(\gamma(y-\beta t-x(t)))e^{\mp\gamma(y-\beta t-x(t))}.}
We thus obtain \eqref{thmasysta:asyc}.

\subsection{The boundedness of the center}\label{secthmvariant}
We prove part (c) of Theorem \ref{thmasysta}. Fix  $\beta_0\in(-1,1)$ and $x_0\in\R$. We assume that $f\in{\bf E}^0_{\sin}$ which is a solution to \eqref{sg} such that
\fm{\norm{(f,f_t)(0)-(K,K_t)(0)}_{H^{m+1,2}\times H^{m,2}}\leq \eps\ll_m1,\qquad m>3/2}
where $K(t,x)=Q(t,x;\beta_0,x_0)$ is defined by \eqref{kink}. We have proved that there exist a new center $y_0$, a new velocity $\beta$, and a solution $\phi\in {\bf E}_{\sin}^0$ to \eqref{sg} with the properties stated in Section \ref{spfthmhrc}. By Proposition \ref{thmwpt}, we have $\phi\in C^1_{t,x}$ and
\fm{\norm{\phi(t)}_{L^\infty}+\norm{\phi_t(t)}_{L^\infty}+\norm{\phi_x(t)}_{L^\infty}\lesssim\eps\lra{t}^{-1/2},\qquad \forall t\geq 0.}

As mentioned in Section \ref{s1ideapf}, in order to prove this boundedness of $x(t)$, we will construct another $C^1$ path $\wt{x}(t)$ defined for all $t\geq 0$ such that $f(t,\beta t+\wt{x}(t))=\pi$ and that $|x(t)-\wt{x}(t)|+|\wt{x}(t)-x_0|\lesssim\eps$ for all $t\geq 0$. Moreover, we will show $\lim_{t\to\infty}|x(t)-\wt{x}(t)|=0$.

We first check that $f\in C_{t,x}^1([0,\infty)\times\R)$. We will rely on the computations in Section \ref{s3pfmthmiiv}.  Let $x(t)$ be the $C^1$ path defined in part (a) of Theorem \ref{thmasysta} and set $(\varphi,\theta)$ by \eqref{defnvarphitheta}. Since $(t,x)\mapsto (f-\varphi,f_t-\theta)(t,x)= (Q,Q_t)(t,x;\beta,x(t))$ is $C^1$, we only need  to check that $\varphi\in C^{1}_{t,x}$ and $\theta\in C_{t,x}$. In fact, since
\fm{\varphi(t)-\varphi(t_0)=f(t)-f(t_0)+Q(t_0,\cdot;\beta,x(t_0))-Q(t,\cdot;\beta,x(t)),}
and since $f(t)\to f(t_0)$ as $t\to t_0$ in $H^1_x$, 
by part b) of Lemma \ref{lemQ}, we obtain $\varphi(t)\to\varphi(t_0)$ as $t\to t_0$ in $H^1_x$. In other words, we have $\varphi\in C_tH^1_x\subset C_{t,x}$. Then, by \eqref{varphixtheta:eq} and since $\phi\in C_{t,x}^1$, we conclude that $\varphi_x,\theta\in C_{t,x}$.

We then present a lemma on the existence and uniqueness of $\wt{x}(t)$.
\lem{\label{lemwtxicst}For each $t\geq 0$, there exists a unique $\wt{x}(t)\in\R$ with $f(t,\beta t+\wt{x}(t))=\pi$. Besides, we have $|x(t)-\wt{x}(t)|\lesssim \eps \lra{t}^{-1/2}$. Moreover, the map $t\mapsto \wt{x}(t)$ is $C^1$.}
\begin{proof}
The existence of such $\wt{x}(t)$ is guaranteed by the continuity of $f$ and the two limits $\lim_{x\to\infty}f(t,x)=2\pi$ and $\lim_{x\to-\infty}f(t,x)=0$. We now show the uniqueness of such a $\wt{x}(t)$. By part (b) of Theorem \ref{thmasysta}, we have $\norm{(f,f_x)(t)-(Q,Q_x)(t,\cdot;\beta,x(t))}_{L^\infty}\lesssim\eps\lra{t}^{-1/2}$. Thus, if $f(t,y)=\pi$ for some $y\in \R$, we must have
\eq{\label{wtxicst}|Q(t,y;\beta,x(t))-\pi|\lesssim\norm{f(t)-Q(t,\cdot;\beta,x(t))}_{L^\infty}\lesssim \eps \lra{t}^{-1/2}.}
We claim that $|y-\beta t-x(t)|\lesssim\eps\lra{t}^{-1/2}$. In fact, for each $C_0>0$ we have 
\fm{&\pi-Q(t,\beta t+x(t)-C_0\eps\lra{t}^{-1/2};\beta,x(t))=\pi-Q(0,-C_0\eps\lra{t}^{-1/2};\beta,0)\\
&= \int_{-C_0\eps\lra{t}^{-1/2}}^{0}2\gamma\sech(\gamma x)\ dx\geq 2C_0\gamma\sech(C_0\gamma\eps)\cdot \eps\lra{t}^{-1/2},\\
&Q(t,\beta t+x(t)+C_0\eps\lra{t}^{-1/2};\beta,x(t))-\pi=Q(0,C_0\eps\lra{t}^{-1/2};\beta,0)-\pi\\
&= \int^{C_0\eps\lra{t}^{-1/2}}_{0}2\gamma\sech(\gamma x)\ dx\geq 2C_0\gamma\sech(C_0\gamma\eps)\cdot \eps\lra{t}^{-1/2}.}
By choosing $C_0\gg_\gamma1$, we can make the $2C_0\gamma\sech(C_0\gamma\eps)$ larger than the implicit constant in \eqref{wtxicst}. In other words, if $f(t,y)=\pi$, then $|y-\beta t-x(t)|\leq C_0\eps\lra{t}^{-1/2}$.

Now, for each $x$ with $|x-\beta t-x(t)|\leq C_0\eps\lra{t}^{-1/2}$, we have 
\fm{f_x(t,x)=f_x(t,x)-Q_x(t,x;\beta,x(t))+Q_x(t,x;\beta,x(t))\geq 2\gamma\sech(C_0\gamma \eps)-C\eps\lra{t}^{-1/2}.}
Here we use $Q_x(t,x;\beta,x(t))=2\gamma\sech(\gamma(x-\beta t-x(t)))$. Thus, by choosing $\eps\ll1$, we have $f_x(t,x)>0$. This indicates that there exists a unique $y$ with $f(t,y)=\pi$. Now we set $\wt{x}(t)=y-\beta t$ where $y$ is the unique solution to $f(t,y)=\pi$. The claim above indicates that $|x(t)-\wt{x}(t)|\lesssim\eps\lra{t}^{-1/2}$.

Finally, to show the map $t\mapsto \wt{x}(t)$ is $C^1$, we apply the implicit function theorem. We have $f(t_0,\beta t_0+\wt{x}(t_0))=\pi$, and previously we have shown that $f_x(t_0,\beta t_0+\wt{x}(t_0))>0$. Thus, we obtain a  $C^1$ function $y(t)$ defined near $t_0$ such that $f(t,\beta t+y(t))=\pi$. The uniqueness proved in the previous paragraph forces $y(t)=\wt{x}(t)$, so $\wt{x}(t)$ is $C^1$.
\end{proof}\rm

Let us now prove that $|\wt{x}(t)-x_0|\lesssim\eps$. All the estimates in Section \ref{spfthmhrc} still hold. Moreover, we  have
\fm{\norm{(\phi,\phi_t)(0)}_{H^{m+1,2}\times H^{m,2}}+\norm{(f,f_t)(0)-(Q,Q_t)(0,\cdot;\beta,y_0)}_{H^{m+1,2}\times H^{m,2}}\lesssim\eps,\qquad m>3/2.}
For such an initial data, we claim that for each $\lambda\in(0,1/2)$, as long as $\eps\ll_\lambda 1$ we also have
\eq{\label{zuptbd}\norm{Z\phi(t)}_{L^\infty}\lesssim\eps\lra{t}^{-\lambda},\qquad \forall t\geq 0.}
Here  $Z=t\partial_x+x\partial_t$ is the Lorentz boost. The proof of this result relies on that of Proposition \ref{thmwpt}, and we will present it in Section \ref{lastappendixb}. 

Let us temporarily assume that $y_0=0$ (recall that $y_0$ comes from Corollary \ref{corsolvebackthigh}). Differentiate $f(t,\beta t+\wt{x}(t))=\pi$ with respect to $t$, and we get
\fm{f_t(t,\beta t+\wt{x}(t))+(\beta +\wt{x}'(t))f_x(t,\beta t+\wt{x}(t))=0.}
By the B\"acklund transform, we have
\eq{\label{xprimespecial}\wt{x}'(t)&=-\frac{f_t+\beta f_x}{f_x}=-\frac{\phi_x+\beta\phi_t+\gamma^{-1}(\sin(\frac{f+\phi}{2})-\sin(\frac{f-\phi}{2}))}{f_x}(t,\beta t+\wt{x}(t)).} Since $f(t,\beta t+\wt{x}(t))=\pi$, in \eqref{xprimespecial} we have \fm{\sin(\frac{f+\phi}{2})-\sin(\frac{f-\phi}{2})=\sin(\frac{\pi+\phi}{2})-\sin(\frac{\pi-\phi}{2})=0.}We also recall that $f_x(t,\beta t+\wt{x}(t))\geq 2\gamma\sech(C_0\gamma\eps)-C\eps$ from the proof of Lemma \ref{lemwtxicst}. Besides, for each $t\geq 1$, we have
\fm{|(\phi_x+\beta\phi_t)(t,\beta t+\wt{x}(t))|&\leq t^{-1}|Z\phi(t,\beta t+\wt{x}(t))|+t^{-1}|\wt{x}(t)||\phi_t(t,\beta t+\wt{x}(t))|\\
&\lesssim \eps t^{-1-\lambda}+\eps t^{-3/2}|\wt{x}(t)|.}
For $0\leq t\leq 1$, we have
\fm{|(\phi_x+\beta\phi_t)(t,\beta t+\wt{x}(t))|\lesssim \eps\lra{t}^{-1/2}\lesssim \eps.}
In summary, for each $t\geq 0$, we have
\fm{|\wt{x}'(t)|\lesssim \eps|\wt{x}(t)| \lra{t}^{-3/2}+\eps\lra{t}^{-1-\lambda}.}
By  Gronwall's inequality, we conclude that $|\wt{x}(t)|\lesssim |\wt{x}(0)|+\eps$.  It is also easy to show that  $|\wt{x}(0)|\lesssim \eps$. We thus conclude that $|\wt{x}(t)|\lesssim\eps$. For a general center $y_0$, we can simply consider the spatial translation $f(t,\cdot-y_0)$.

\appendix

\section{Proof of Proposition \ref{thmlowdata}}\label{seclow}

In this section, we prove Proposition \ref{thmlowdata}. Let $\phi$ be a global solution to the sine-Gordon equation with $\norm{(\phi,\phi_t)(0)}_{H^{1,s}\times L^{2,s}}=\eps\ll_s1$. We seek to prove that
\eq{\label{seclow:finalest}\lim_{t\to\infty}\kh{\norm{\phi(t)}_{L^\infty(\R)}+\norm{(\phi_t,\phi_x)(t)}_{(L^2+L^\infty)(\R)}+\norm{(\phi,\phi_t,\phi_x)(t)}_{L^2(\{x\in\R:\ |x|\geq t\})}}=0.}

\subsection{Reduction of the problem}

We start our proof with a reduction of the problem. We claim that it suffices to prove the following proposition.

\prop{\label{thmschwdata}Fix $s>1/2$. Let $\phi$ be a global solution to the sine-Gordon equation with $C_c^\infty$ initial data. Then, whenever  $\norm{(\phi,\phi_t)(0)}_{H^{1,s}\times L^{2,s}}\ll_s1$,  we have \eq{\label{seclow:finalest2}\lim_{t\to\infty}(\norm{\phi(t)}_{L^\infty(\R)}+\norm{\phi_t(t)}_{L^\infty(\R)}+\norm{\phi_x(t)}_{L^\infty(\R)})=0.}
}\rm

Note that the smallness of the $H^{1,s}\times L^{2,s}$ norm (with $s>1/2$) of the initial data rules out the possibility of solitons.  We refer to  \cite[Remark 2.2]{chen2020long}. Thus, one can apply Cheng-Venakides-Zhou \cite[Theorem 1]{MR1697487} to prove this proposition.  There is no need to assume that higher derivatives of the initial data are small. Since the choice of $\eps$ in Proposition \ref{thmwpt} depends on the $H^{5/2+,1}\times H^{3/2+,1}$ norm of the initial data, we cannot directly apply that proposition to obtain Proposition \ref{thmschwdata}. 

We postpone the proof of Proposition \ref{thmschwdata} to the next subsection. Here let us prove that Proposition \ref{thmschwdata} implies Proposition \ref{thmlowdata}. We will need the following result from the proof of  \cite[Theorem 11.5]{chen2020long}. There the authors applied the inverse scattering method, so this result relies on the integrability of the sine-Gordon equation. 
\prop[Chen-Liu-Lu, the proof of Theorem 11.5 in \cite{chen2020long}]{\label{chen2020thm}Fix $s>1/2$. Recall the definition of ${\bf E}^{*,*}$ in Section \ref{secnotation}. Let $\phi\in{\bf E}^{0,s}$ and $\phi_n\in{\bf E}^{1,1}$, $n=1,\dots$ be global solutions to the sine-Gordon equation \eqref{sg}.  Suppose that $\norm{(\phi,\phi_t)(0)}_{H^{1,s}\times L^{2,s}}=\eps\ll_s1$ and that
\fm{(\phi_n,\phi_{n,t})(0)\to (\phi,\phi_{t})(0),\qquad \text{as }n\to\infty \text{ in }H^{1,s}\times L^{2,s}.}
Then, we have
\fm{\lim_{n\to\infty}\sup_{t\geq 0}\norm{(\phi-\phi_n,\phi_{t}-\phi_{n,t})(t)}_{H^{1}(\R)\times L^2(\R)}=0.}
}\rm

Now, we take a sequence of $\{(\phi_{0,n},\phi_{1,n})\}$ in $C_c^\infty\times C_c^\infty$ such that
\fm{(\phi_{0,n},\phi_{1,n})\to (\phi_0,\phi_1),\qquad \text{as }n\to\infty \text{ in }H^{1,s}\times L^{2,s}.}
Without loss of generality, we assume that \fm{\sup_n\norm{(\phi_{0,n},\phi_{1,n})}_{H^{1,s}\times L^{2,s}}\leq2 \eps\ll_s1.}
Let $\phi_n$ be the global solution to \eqref{sg} with initial data $(\phi_{0,n},\phi_{1,n})$. Suppose that we have already proved \eqref{seclow:finalest2}, so for each integer $n$, 
\eq{\label{seclow:finalest2:var}\lim_{t\to\infty}\norm{(\phi_n,\phi_{n,t},\phi_{n,x})(t)}_{L^\infty(\R)}=0.}
Moreover, by the finite speed of propagation (see, e.g., \cite[Theorem I.2.2]{MR2455195}), for each $n$ we can find $R_n>0$ such that $\phi_n$ vanishes in $\{|x|>t+R_n\}$. As a result, from \eqref{seclow:finalest2:var} we have
\eq{\label{seclow:finalest2:var2}\lim_{t\to\infty}\norm{(\phi_n,\phi_{n,t},\phi_{n,x})(t)}_{L^2(\{x\in\R:\ |x|\geq t\})} =0.}

To prove \eqref{seclow:finalest}, we fix an arbitrary constant $\eta>0$. By Proposition \ref{chen2020thm}, we can choose a large integer $N$ so that 
\fm{\sup_{t\geq 0}\norm{(\phi-\phi_N,\phi_t-\phi_{N,t},\phi_{x}-\phi_{N,x})(t)}_{ L^2(\R)}< \frac{\eta}{4(1+C_0)}.}
Here we let $C_0>0$ be the constant in the Sobolev embedding $H^1(\R)\subset L^\infty(\R)$. As a result, we have
\fm{\sup_{t\geq 0}\norm{(\phi-\phi_N)(t)}_{ L^\infty(\R)}< \frac{C_0\eta}{4(1+C_0)}<\eta/4.}
By \eqref{seclow:finalest2:var} and \eqref{seclow:finalest2:var2} with $n$ replaced by this $N$, we can choose $T=T_{N,\eta}>0$ so that
\fm{\norm{(\phi_N,\phi_{N,t},\phi_{N,x})(t)}_{L^\infty(\R)}+\norm{(\phi_N,\phi_{N,t},\phi_{N,x})(t)}_{L^2(\{x\in\R:\ |x|\geq t\})}<\eta/4,\qquad \forall t>T.}
In summary, for all $t> T$, we have
\fm{&\norm{ \phi(t)}_{L^\infty(\R)}+\norm{(\phi_t,\phi_x)(t)}_{(L^2+L^\infty)(\R)}+\norm{(\phi,\phi_t,\phi_x)(t)}_{L^2(\{x\in\R:\ |x|\geq t\})}\\
&\leq \norm{(\phi-\phi_N)(t)}_{L^\infty(\R)}+2\norm{(\phi-\phi_N,\phi_t-\phi_{N,t},\phi_x-\phi_{N,x})(t)}_{L^2(\R)}\\
&\quad+\norm{(\phi_N,\phi_{N,t},\phi_{N,x})(t)}_{L^\infty(\R)}+\norm{(\phi_N,\phi_{N,t},\phi_{N,x})(t)}_{L^2(\{x\in\R:\ |x|\geq t\})}<\eta.}
That is, we obtain \eqref{seclow:finalest}.

\subsection{Proof of Proposition \ref{thmschwdata}}
Suppose that $\phi$ is a global solution to the sine-Gordon equation in Proposition \ref{thmschwdata}. We first claim that
\eq{\label{schwptb}\lim_{t\to\infty} \norm{\phi(t)}_{L^\infty(\R)}=0.}
To prove this claim, we recall Cheng-Venakides-Zhou \cite[Theorem 1]{MR1697487}. There the authors proved that \fm{\lim_{t\to\infty} \norm{\phi(t)}_{L^\infty(\{x\in\R:\ |x|\geq t\})}=0.}
Moreover, they proved an asymptotic formula for $\sin(\phi(t,x))$ with $|x|<t$. By setting $z_0=\sqrt{\frac{|x-t|}{|x+t|}}$ and $\tau=\frac{tz_0}{1+z_0^2}$, one has
\fm{|\sin(\phi(t,x))|&\leq\sqrt{\frac{4\ln(1+|r(z_0)|^2)}{\pi \tau}}+C(z_0)\frac{\ln \tau}{\tau},\qquad \text{for } |x|<t,\ \tau\to\infty.}
Here $r(z_0)$ is the reflection coefficient corresponding to the initial data of $\phi$, and $C(z_0)$ is a constant. Under the assumptions of Proposition \ref{thmschwdata}, we have both $r(z_0)$ and $C(z_0)$ decay rapidly as $z_0\to\infty$. Thus, whenever $z_0>1/2$ (or equivalently, $-1<x/t<3/5$), we have $t/(5z_0)\leq \tau\leq t/z_0$ and
\fm{|\sin(\phi(t,x))|&\lesssim \sqrt{ \frac{|r(z_0)|^2z_0}{t}}+C(z_0)\cdot (t/z_0)^{-2/3}\lesssim t^{-1/2}.}
The implicit constants can be chosen to be independent of $z_0$ because of the rapid decays of $r(z_0)$ and $C(z_0)$. For example, we have $|r(z_0)|^2z_0\lesssim z_0^{-2}\cdot z_0=z_0^{-1}<2$ for $z_0>1/2$. By the smallness of the $H^1\times L^2$ norm of the initial data, we have $|\phi|\ll1$ and $|\phi|\sim|\sin(\phi)|$; see Section \ref{secpforbstac3}. It follows that
\fm{\lim_{t\to\infty}\norm{\phi(t)}_{L^\infty(\{x\in\R:\ -t<x< \frac{3}{5}t\})}=0.}
Applying the same result to $\phi(t,-x)$ yields
\fm{\lim_{t\to\infty}\norm{\phi(t)}_{L^\infty(\{x\in\R:\ -\frac{3}{5}t<x<t\})}=0.}
This finishes the proof of our claim.

Next, we claim that
\eq{\label{schwclaim}\norm{\phi_{tt}(t)}_{L^\infty(\R)}+\norm{\phi_{xx}(t)}_{L^\infty(\R)}\lesssim  1,\qquad \forall t\geq 0.}
Let us explain how to finish the proof assuming this claim holds. For any $f\in C^2([-1,1])$, we have
\eq{|f'(0)|&\lesssim \norm{f}_{L^\infty([-1,1])}^{1/2}\norm{f''}_{L^\infty([-1,1])}^{1/2}+\norm{f}_{L^\infty([-1,1])}.}
See, for example, \cite[(1.1) or Exercise 1.6]{MR1269107}. As a result, by setting $f(y)=\phi(t,x+y)$ and setting $f(y)=\phi(t+y,x)$ respectively in the inequality above, we have
\fm{|\phi_x(t,x)|&\lesssim \norm{\phi(t)}_{L^\infty([x-1,x+1])}^{1/2}\norm{\phi_{xx}(t)}_{L^\infty([x-1,x+1])}^{1/2}+\norm{\phi(t)}_{L^\infty([x-1,x+1])},\\
|\phi_t(t,x)|&\lesssim \norm{\phi(\cdot,x)}_{L^\infty([t-1,t+1])}^{1/2}\norm{\phi_{tt}(\cdot,x)}_{L^\infty([t-1,t+1])}^{1/2}+\norm{\phi(\cdot,x)}_{L^\infty([t-1,t+1])}.}
Taking the supreme over all $x\in\R$, we have
\fm{\norm{\phi_x(t)}_{L^\infty(\R)}&\lesssim \norm{\phi(t)}_{L^\infty(\R)}^{1/2}\norm{\phi_{xx}(t)}_{L^\infty(\R)}^{1/2}+\norm{\phi(t)}_{L^\infty(\R)},\\
\norm{\phi_t(t)}_{L^\infty(\R)}&\lesssim \sup_{t'\geq t-1}(\norm{\phi(t')}_{L^\infty(\R)}^{1/2}\norm{\phi_{tt}(t')}_{L^\infty(\R)}^{1/2}+\norm{\phi(t')}_{L^\infty(\R)}).} 
Take $t\to\infty$ and we finish the proof of Proposition \ref{thmschwdata}.

We now return to the proof of the claim \eqref{schwclaim}. Since
$\phi_{tt}=\phi_{xx}-\sin(\phi)=\phi_{xx}+O(1)$,
it suffices to estimate $\phi_{xx}$. By the Sobolev embedding $H^1(\R)\subset L^\infty(\R)$, we only need to prove the following proposition.
\prop{Let $\phi$ be a global solution to the sine-Gordon equation with $C_c^\infty$ initial data. Then, 
\eq{\partial_x^k\partial_t^l\phi\in L^\infty_tL^2_x([0,\infty)\times\R),\qquad k+l\leq 3.}}
\begin{proof}
By the usual conservation law and since $|\sin(\phi)|\sim|\phi|$ whenever $|\phi|\ll1$, we have $\partial_x^k\partial_t^l\phi\in L^\infty L^2$ whenever $k+l\leq 1$. Set $\partial_+=\frac{1}{\sqrt{2}}(\partial_t+\partial_x)$ and $\partial_-=\frac{1}{\sqrt{2}}(\partial_t-\partial_x)$. By following the computations in 
\cite{PhysRevD.13.3440}, we obtain
\eq{\label{propb.1f}\partial_+J_{2n,a}^++\partial_-J_{2n,a}^-=\partial_+J_{2n,b}^++\partial_-J_{2n,b}^-=0,\qquad n=1,2.}
Here $J_{*,*}^\pm$ are called the \emph{energy currents} defined by
\fm{J_{2,a}^+=(\phi_{--})^2-\frac{1}{4}(\phi_-)^2,\qquad J_{2,a}^-=\frac{1}{2}(\phi_{-})^2\cos\phi;}
\fm{J_{4,a}^+&=(\phi_{---})^2+\frac{5}{2}(\phi_-)^2(\phi_{--})^2+\frac{5}{3}(\phi_-)^2\phi_{---}+\frac{1}{8}(\phi_-)^6,\\
J_{4,a}^-&=-\frac{5}{3}(\phi_-)^3\phi_{--+}-\frac{3}{8}(\phi_-)^4\cos\phi+\frac{3}{2}(\phi_-)^2\phi_{--}\sin\phi+\frac{1}{2}(\phi_{--})^2\cos\phi.}
To define $J_{2n,b}^\pm$, we simply replace $a$ with $b$ and interchange  $+$ and $-$ everywhere in the formulas above. We now rewrite \eqref{propb.1f} as 
\fm{\partial_t(J_{2n,a}^++J_{2n,a}^-)+\partial_x(J_{2n,a}^+-J_{2n,a}^-)=\partial_t(J_{2n,b}^++J_{2n,b}^-)+\partial_x(J_{2n,b}^+-J_{2n,b}^-)=0.}
We thus obtain two conserved quantities:
\fm{E_{2n}(t)&=\int_\R (J_{2n,a}^++J_{2n,a}^-+J_{2n,b}^++J_{2n,b}^-)(t,x)\ dx,\qquad n=1,2.}
We also set 
\fm{E_0(t)&=\int_\R \frac{1}{2}(\phi_t^2+\phi_x^2)+1-\cos\phi\ dx}
which is the usual conserved energy.

We now have 
\fm{E_2(t)&=\int_\R (\phi_{--})^2+(\phi_{++})^2-\frac{1}{4}((\phi_-)^2+(\phi_+)^2)+\frac{1}{2}((\phi_-)^2+(\phi_+)^2)\cos\phi\ dx\\
&=\int_\R \frac{1}{2}(\phi_{tt}+\phi_{xx})^2+2\phi_{tx}^2\ dx+O(E_0(t))=\int_\R \frac{1}{2}(2\phi_{xx}-\sin\phi)^2+2\phi_{tx}^2\ dx+O(1).}
Since $\phi(0)$ and $\phi_t(0)$ are $C_c^\infty$, we have $E_2(0)<\infty$ and thus $2\phi_{xx}-\sin\phi,\phi_{tx}\in L^\infty L^2$. But since $|\sin(\phi)|^2\leq 4|\sin(\phi/2)|^2\leq 2(1-\cos\phi)$, we have $\sin\phi\in L^\infty L^2$. It follows that $\phi_{xx}\in L^\infty L^2$ and $\phi_{tt}=\phi_{xx}-\sin\phi\in L^\infty L^2$.

Finally, we compute $E_4(t)$. Note that
\fm{&\int_\R (\phi_-)^2(\phi_{---}-\phi_{--+})+(\phi_+)^2(\phi_{+++}-\phi_{++-})\ dx\\
&=\int_\R -\frac{1}{\sqrt{2}}(\phi_-)^2\phi_{--x}+\frac{1}{\sqrt{2}}(\phi_+)^2\phi_{++x}\ dx=\int_\R \sqrt{2}\phi_-\phi_{-x}\phi_{--}-\sqrt{2}\phi_+\phi_{+x}\phi_{++}\ dx\\
&=O\kh{\norm{\phi_-}_{L^\infty}\norm{\phi_{-x}}_{L^2}\norm{\phi_{--}}_{L^2}+\norm{\phi_+}_{L^\infty}\norm{\phi_{+x}}_{L^2}\norm{\phi_{++}}_{L^2}}=O(1).}The implicit constant in $O(1)$ is independent of $t$. To get the last estimate, we make use of the Sobolev embedding: $\norm{(\phi_t,\phi_x)(t)}_{L^\infty}\lesssim \norm{(\phi_t,\phi_x)(t)}_{H^1}\lesssim1$. Similarly, the integral
\fm{\int_\R-\frac{3}{8}((\phi_-)^4+(\phi_+)^4)\cos\phi+\frac{3}{2}[(\phi_-)^2\phi_{--}+(\phi_+)^2\phi_{++}]\sin\phi+\frac{1}{2}((\phi_{--})^2+(\phi_{++})^2)\cos\phi\ dx}
is also $O(1)$. And since $E_4(0)<\infty$ (because the initial data of $\phi$ is $C_c^\infty$), we conclude that $\phi_{+++},\phi_{---}\in L^\infty L^2$. To finish the proof, we simply notice that
\fm{2\sqrt{2}\phi_{+++}&=\phi_{ttt}+3\phi_{ttx}+3\phi_{txx}+\phi_{xxx}=4\phi_{xxt}-\phi_t\cos\phi+4\phi_{xxx}-3\phi_x\cos\phi,\\
2\sqrt{2}\phi_{---}&=\phi_{ttt}-3\phi_{ttx}+3\phi_{txx}-\phi_{xxx}=4\phi_{xxt}-\phi_t\cos\phi-4\phi_{xxx}+3\phi_x\cos\phi.}
Since $\phi_t\cos\phi,\phi_x\cos\phi\in L^\infty L^2$, we conclude that $\phi_{xxx},\phi_{xxt}\in L^\infty L^2 $. We finally notice that $\phi_{ttt}=\phi_{txx}-\phi_t\cos\phi\in L^\infty L^2$ and $\phi_{xtt}=\phi_{xxx}-\phi_x\cos\phi\in L^\infty L^2$.
\end{proof}\rm

\section{Proof of Proposition \ref{thmwpt}}
\label{secthmwpt}
In this section, we sketch the proof of Proposition \ref{thmwpt} by following the method in \cite{MR4693097}. Since the sine-Gordon equation is globally well-posed in the energy space $H^1_{\sin}\times L^2$, we already obtain a global solution $\phi\in {\bf E}^0_{\sin}$ to \eqref{sg}.  It remains to check the bound
\eq{\label{thmwptbd}\norm{\phi(t)}_{L^\infty}+\norm{\phi_t(t)}_{L^\infty}+\norm{\phi_x(t)}_{L^\infty}\lesssim\eps\lra{t}^{-1/2},\qquad \forall t\geq 0}
and to show the asymptotic for $\phi$ and its derivatives.

By setting $u=\phi+i\lra{D}^{-1}\phi_t$, we convert \eqref{sg} to a first order PDE
\eq{\label{sg1}
\left\{
\begin{array}{l}
     \displaystyle(\partial_t+i\lra{D})u=i\lra{D}^{-1}(\Real u-\sin(\Real u)),\\[1em]
     u|_{t=0}=u_0\in H^{m,1},\quad \norm{u_0}_{H^{m,1}}\lesssim \eps\ll1.
\end{array}\right.
}
Here $m>5/2$ is an arbitrary fixed real constant. Without loss of generality, we also assume $m\leq 3$.

To prove \eqref{thmwptbd}, we set up a continuity argument. Fix some time $T\geq 0$, and suppose we have
\eq{\label{appca}\norm{\lra{D}^{1+\delta}u(t)}_{L^\infty}\leq A\eps \lra{t}^{-1/2},\qquad \forall t\in[0,T].} Here $A>1$ is a  large constant and $0<\delta<\frac{m-5/2}{4}$ is a  small constant. Their values will be chosen later.   Because of the Sobolev embedding, the bound \eqref{appca} holds for $T=0$ for sufficiently large $A$ depending on the initial data. Now we seek to show that $t\mapsto \lra{D}^{1+\delta}u(t)$ is continuous from $\R$ to $L^\infty$ and that \eqref{appca} holds with $A$ replaced by $A/2$ as long as $\eps\ll1$. In Section \ref{wtpenergyest}, we derive several necessary energy estimates for $u$. Then, we prove \eqref{thmwptbd} in Section \ref{sec:b.2:ptwisebdd}.

After proving \eqref{thmwptbd}, we derive the asymptotics for the solution in Section \ref{secb.3:asy}. In the proof, we follow the discussion in \cite{MR3382579}. Finally, in Section \ref{lastappendixb}, we prove \eqref{zuptbd}. 

\subsection{The energy estimates}\label{wtpenergyest} Before applying the method of testing by wave packets, we need to first derive several energy bounds for $u$. Set $L=x-tD\lra{D}^{-1}$ for each time $t\geq 0$. Define a function space $X=X_t$ equipped with the norm \fm{\norm{u(t)}_X=\norm{u(t)}_{H^m}+\norm{Lu(t)}_{H^m}.}
Our goal in this subsection is to show $\norm{u(t)}_X\lesssim\eps \lra{t}^{C_A\eps^2}$ for all $t\geq 0$. Here the implicit constant in ``$\lesssim$'' is independent of the choice of $A$. 

By the chain rule, we have
\fm{\norm{(\partial_t+i\lra{D})u}_{H^m}&\lesssim \norm{\Real u-\sin(\Real u)}_{H^{m-1}}\lesssim \norm{u}_{L^\infty}^2\norm{u}_{H^{m-1}}\lesssim_A \eps^2\lra{t}^{-1}\norm{u}_{H^m}.}
Here we use \eqref{appca} and the estimate $\norm{fg}_{H^{m-1}}\lesssim \norm{f}_{L^\infty}\norm{g}_{H^{m-1}}+\norm{g}_{L^\infty}\norm{f}_{H^{m-1}}$. By  Gronwall's inequality and the initial condition in \eqref{sg1}, we obtain $\norm{u(t)}_{H^m}\lesssim\eps\lra{t}^{C_A\eps^2}$. 
Next, since $[\psi(D),x]=\psi'(D)$ and $m-1\leq 2$, we have
\fm{&\norm{(\partial_t+i\lra{D})(xu)}_{H^m}\lesssim \norm{x(\partial_t+i\lra{D})u}_{H^m}+\norm{u}_{H^m}\lesssim \norm{x[\Real u-\sin(\Real u)]}_{H^2}+\norm{u}_{H^m}\\
&\lesssim \norm{u}_{L^\infty}^2\norm{xu}_{H^2}+\norm{u}_{L^\infty}\norm{u_x}_{L^\infty}\norm{xu}_{H^1}+\norm{u}_{H^m}\lesssim A^2\eps^2\lra{t}^{-1}\norm{xu}_{H^m}+\eps\lra{t}^{C_A\eps^2}. }
By  Gronwall's inequality and the initial condition, we obtain $\norm{xu(t)}_{H^3}\lesssim\eps\lra{t}^{1+C_A\eps^2}$.
Finally, we define the Lorentz boost $Z=t\partial_x+x\partial_t$. We have $[\partial_t+i\lra{D},Z]=D\lra{D}^{-1}(\partial_t+i\lra{D})$ and $[Z,\lra{D}^{-1}]=D\lra{D}^{-3}\partial_t$. Since $m-2\leq 1$ and $m-3\leq 0$, we have
\fm{&\norm{(\partial_t+i\lra{D})(Zu)}_{H^{m-1}}\lesssim \norm{Z(\partial_t+i\lra{D})u}_{H^{m-1}}+\norm{(\partial_t+i\lra{D})u}_{H^{m-1}}\\
&\lesssim \norm{Z(\Real u-\sin(\Real u))}_{H^1}+\norm{\partial_t(\Real u-\sin(\Real u))}_{L^2}+\norm{(\partial_t+i\lra{D})u}_{H^{m-1}}\\
&\lesssim (\norm{u}_{L^\infty}+\norm{u_x}_{L^\infty})^2(\norm{Zu}_{H^1}+\norm{u_t}_{L^2}+\norm{u}_{H^{m-2}})\lesssim A^2\eps^2\lra{t}^{-1}(\norm{Zu}_{H^{m-1}}+\norm{u}_{H^m}).}
Here we use $u_t=(\partial_t+i\lra{D})u-i\lra{D}u$.  Add this estimate with the estimate for $\norm{(\partial_t+i\lra{D})u}_{H^m}$ and then apply  Gronwall's inequality. We conclude that $\norm{Zu(t)}_{H^{m-1}}\lesssim \eps \lra{t}^{C_A\eps^2}$. By noticing that $Z=(\partial_t+i\lra{D})x-i\lra{D}L$, we conclude
\eq{\norm{u(t)}_{X}\lesssim \eps \lra{t}^{C_A\eps^2},\qquad \forall t\in[0,T].}

Note that this inequality implies that 
\fm{\norm{\lra{D}^{1+\delta}u(t)}_{L^\infty}\lesssim \norm{u(t)}_{H^{3}}\lesssim \eps\lra{t}^{C_A\eps^2}.}
Thus, we only need to improve \eqref{appca} for  $t\geq 1$. Moreover, since $\phi=\Real u\in C_tH^{m}(\R\times\R)$ and $\phi_t=\Imag \lra{D}u\in C_tH^{m-1}(\R\times\R)$, we have $\phi\in C_{t,x}^1(\R\times\R)$ by the Gagliardo-Nirenberg inequality.

\subsection{The pointwise estimates}\label{sec:b.2:ptwisebdd}
\subsubsection{Setup}\label{wptsetup}
Fix a small constant $\mu\in(0,1/4)$, and choose $n_0\in\Z_{<0}$  so that $(1+\mu)^{n_0}<\mu$. Now, we set
\fm{\mcl{I}=\{(1+\mu)^n:\ n\geq n_0\}.}
For each $\lambda\in\mcl{I}$, we define 
\fm{I_\lambda=\{\eta\in\R:\ \frac{1}{1+\mu}<|\eta/\lambda|\leq 1+\mu\},\qquad \lambda>(1+\mu)^{n_0},\\
I_\lambda=\{\eta\in\R:\ |\eta/\lambda|< 1+\mu\},\qquad \lambda=(1+\mu)^{n_0}.}
We also define a family of   Littlewood-Paley multipliers $P_{\leq \lambda}$ and $P_\lambda$ which are associated to $\{I_\lambda\}_{\lambda\in\mcl{I}}$. Fix  a real-valued even bump function $\varphi$ such that $0\leq \varphi\leq 1$, $\supp\varphi\subset(-1-\mu,1+\mu)$ and $\varphi|_{[-1,1]}\equiv 1$. Set 
\fm{\F(P_{\leq \lambda }f)(\xi)=\varphi(\xi/\lambda)\wh{f}(\xi),\qquad \F(P_{ \lambda}f)(\xi)=(\varphi(\xi/\lambda)-\varphi((1+\mu)\xi/\lambda))\wh{f}(\xi).}

We also set
\fm{J_\lambda=\{\xi\lra{\xi}^{-1}:\ \xi\in I_\lambda\}.}
It is clear that the map $\xi\mapsto\xi/\lra{\xi}$ is  bijective  from $I_\lambda$ to $J_\lambda$ with inverse $v\mapsto \xi_v=v/\sqrt{1-v^2}$.
Next we set 
\eq{\mcl{D}=\bigcup_{\lambda\in\mcl{I}}(\lambda,\infty)\times J_\lambda.}
For each $(t,v)\in \mcl{D}$, we define a wave packet
\eq{\Psi_v(t,x)&=\lra{\xi_v}^{3/2}\chi(t^{-1/2}\lra{\xi_v}^{3/2}(x-vt))e^{-i\sqrt{t^2-x^2}}.}
Here the cutoff function $\chi\in C_c^\infty(\R)$  satisfies $\int\chi=1$ and  $\chi\geq 0$.  Moreover, we  choose $\chi$ so that $|x|<t$ in $\supp\Psi$. In fact, since $(t,v)\in\mcl{D}$, we have $t>\lambda$ and $v\in J_\lambda$ for some $\lambda\in\mcl{I}$. Then we have $\xi_v\in I_\lambda$ and thus $\lra{\xi_v}\sim \lra{\lambda}\sim\lambda$. In other words, by setting $y=t^{-1/2}\lra{\xi_v}^{3/2}(x-vt)$, we have
\fm{|x/t-v|=t^{-1/2}\lra{\xi_v}^{-3/2}|y|\lesssim t^{-1/2}\lambda^{-3/2}|y|\lesssim\lambda^{-2}|y|\lesssim\lra{\lambda}^{-2}|y|.}
Meanwhile, since $(\xi/\lra{\xi})'=\lra{\xi}^{-3}>0$ and since $\xi_v\in I_\lambda$, we have
\fm{1-|v|\geq 1-\frac{(1+\mu)\lambda}{\lra{(1+\mu)\lambda}}=\frac{1}{\sqrt{1+(1+\mu)^2\lambda^2}\cdot (\sqrt{1+(1+\mu)^2\lambda^2}+(1+\mu)\lambda)}\sim \lra{\lambda}^{-2}.}
Thus, if $\supp\chi\subset(-c,c)$ for some sufficiently small $c>0$, for $(t,x)\in\supp\Psi_v$ we have
\fm{|x/t-v|\leq Cc\lra{\lambda}^{-2}\leq \frac{1-|v|}{2}.}
This forces $|x|<t$. Here we emphasize that the choice of $c$ is independent of $\lambda$ and $v$.

Define the asymptotic profile $\gamma$ for a given function $\psi=\psi(t,x)$.  For $(t,v)\in\mcl{D}$, we set
\eq{\label{gammadefn}\gamma(t,v)&=\lra{\psi,\Psi_v}=\int_\R \psi(t,x)\overline{\Psi}_v(t,x)\ dx.}
Since $\chi\in C_c^\infty$,  the integral here is well-defined and finite.

For each  $\lambda\in\mcl{I}$, we write
\eq{\label{phildefn}u_\lambda=P_{\lambda}u,\qquad\text{if }\lambda> (1+\mu)^{n_0};\\
u_\lambda=P_{\leq\lambda}u,\qquad\text{if }\lambda= (1+\mu)^{n_0}.}
Here $P_{\leq \lambda}$ and $P_\lambda$ are the corresponding Littlewood-Paley multipliers defined above.
As a result, $u=\sum_{\lambda\in\mcl{I}}u_\lambda$. Moreover, we  set
\eq{\label{wtphidefn}\wt{u}_\lambda=P_{\leq \lambda(1+\mu)}u-P_{\leq \lambda(1+\mu)^{-2}}u,\qquad\text{if }\lambda> (1+\mu)^{n_0+1};\\
\wt{u}_\lambda=P_{\leq\lambda(1+\mu)}u,\qquad\text{if }\lambda\leq (1+\mu)^{n_0+1}.}
An important fact is that the symbol of the Fourier multiplier $u\mapsto \wt{u}_\lambda$ equals $1$ in an open neighborhood of $I_\lambda$.

\subsubsection{Reduction of the problem}\label{secBreduction}
We now prove \eqref{appca} with $A$ replaced by $A/2$. Let $M>1$ be a sufficiently large constant to be chosen. We claim that it suffices to prove
\eq{\label{gammaclaim}\sup_{v\in J_\lambda}|\gamma(t,v)|\lesssim \eps \lambda^{-1-2\delta},\qquad  \forall \lambda^M+1\lesssim t\leq T.}
Here we recall that $\gamma$ is defined by \eqref{gammadefn} and that the upper bound $T$ comes from the assumption \eqref{appca}.
Since the proof of this claim will also be useful in Section \ref{secb.3:asy}, we summarize the main result in the next lemma.

\lem{\label{sec:appB:keylemma}Fix $\delta\in (0,\frac{m-5/2}{4})$. Assume that \eqref{gammaclaim} holds for a constant $M>1$ depending on $\delta$. Then, there exists a sufficiently large constant $M>1$ depending on $\delta$ (such that \eqref{gammaclaim} still holds for this new constant $M$) and a sufficiently small constant $\kappa>0$ depending on $\delta$ and $M$,  such that for all $l\in[0,1+\delta]$, $\lambda\in\mcl{I}$, $\lambda^M+1\lesssim t\leq T$,  and $v\in J_\lambda$, we have 
\eq{\label{sec:appB:keylemma:c}\lra{D}^lu(t,vt)&=t^{-1/2}e^{-it\sqrt{1-v^2}}\lra{\xi_v}^{l}\gamma(t,v)+O_\kappa(\eps t^{-1/2-\kappa}),\qquad \text{with }\xi_v=v/\sqrt{1-v^2}}
as long as $\eps\ll_{A,M,\delta,\kappa,m}1$. Moreover, for all $l\in[0,1+\delta]$ and all $(t,v)\in ([1,T]\times\R)\setminus\bigcup_{\lambda\in\mcl{I}}\{t\gtrsim\lambda^M+1,\ v\in J_\lambda\}$, we have
\eq{\label{sec:appB:keylemma:c2}\lra{D}^lu(t,vt)&=O_\kappa(\eps t^{-1/2-\kappa}),\qquad \text{as long as }\eps\ll_{A,M,\delta,\kappa,m}1.}
Note that none of the constants depend on $T$.
}
\rm
\bigskip

Assuming that \eqref{gammaclaim} and Lemma \ref{sec:appB:keylemma} have been proved, we now prove \eqref{appca} with $A$ replaced by $A/2$. Fix an arbitrary $(t,x)\in[1,T]\times\R$. If $(t,x/t)\notin\bigcup_{\lambda\in\mcl{I}}\{t\gtrsim\lambda^M+1,\ v\in J_\lambda\}$, we obtain from \eqref{sec:appB:keylemma:c2} that $\lra{D}^{1+\delta}u(t,x)=O(\eps t^{-1/2-\kappa})$ where the implicit constant is independent of~$A$. If we have $x/t\in J_\lambda$ and $t\gtrsim \lambda^M+1$ for some $\lambda\in\mcl{I}$,  by \eqref{gammaclaim} and \eqref{sec:appB:keylemma:c}, we have
\fm{|\lra{D}^{1+\delta}u(t,x)|&\lesssim t^{-1/2}\lambda^{1+\delta}\cdot \eps \lambda^{-1-2\delta}+\eps t^{-1/2-\kappa}\lesssim \eps t^{-1/2}.}
The implicit constant in $\lesssim$ is independent of $A$ again. Thus, by choosing $A\gg1$, we obtain \eqref{appca} with $A$ replaced by $A/2$.

Let us prove Lemma \ref{sec:appB:keylemma}. We first recall a Klainerman-Sobolev type inequality from  \cite[Lemma 3.2]{MR2457221}:
\eq{\label{ptwbdlem3.2}\norm{\psi(t)}_{L^\infty}\lesssim \lra{t}^{-1/2}\norm{\psi(t)}_{H^{3/2}}^{1/2}(\norm{\psi(t)}_{H^{3/2}}^{1/2}+\norm{L\psi(t)}_{H^{3/2}}^{1/2}).}
Recall that $L=x-tD\lra{D}^{-1}$. It follows that for all $t\in[1,T]$ and $\lambda\in \mcl{I}$,
\eq{\norm{u_\lambda}_{L^\infty}\lesssim \lra{t}^{-1/2}(\norm{u_\lambda}_{H^{3/2}}+\norm{u_\lambda}_{H^{3/2}}^{1/2}\norm{Lu_\lambda}_{H^{3/2}}^{1/2})\lesssim \eps t^{-1/2}\lambda^{3/2-m}\norm{u}_X\lesssim \eps t^{-1/2+C_A\eps^2}\lambda^{3/2-m}.}
Here we use the fact that $Lu_\lambda$ is localized at frequency $\lambda$. 
If $t\lesssim \lambda^M+1$, we have \fm{\norm{u_\lambda}_{L^\infty}\lesssim \eps t^{-1/2}\lambda^{3/2-m+MC_A\eps^2}\lesssim \eps t^{-1/2}\lambda^{-1-2\delta}}for $\delta\ll_m1$ and $\eps\ll_{A,M,\delta}1$.
In fact, since $m-3/2>1$, we first choose $\delta\in(0,1)$ so that $1+4\delta<m-3/2$. The choice of $\delta$ depends only on $m$ and the value of $\delta$ will not be changed from now on. We then choose $\eps\ll_{A,M,\delta}1$ so that $MC_A\eps^2<2\delta$. This gives us the last estimate above.
By the Bernstein inequalities, for each $l\in[0,1+\delta]$ and $t\in[1,T]$, we have
\eq{\label{est:b.13:D}\sum_{\lambda: \lambda^M+1\gtrsim t}\norm{\lra{D}^{l}u_\lambda(t)}_{L^\infty}\lesssim \sum_{\lambda:\lambda^M+1\gtrsim t}\lambda^{l}\cdot \eps t^{-1/2}\lambda^{-1-2\delta}\lesssim\sum_{\lambda:\lambda^M+1\gtrsim t}\eps t^{-1/2}\lambda^{-\delta} \lesssim \eps t^{-1/2-\delta/M}.}
The implicit constants depend on $\mu$ (see Section \ref{wptsetup}) and $\delta$ but not on $M$. For example, one can show that $\lambda\gtrsim t^{1/M}$ whenever $\lambda^M+1\gtrsim t$ and $t\geq 1$. Here the implicit constants are independent of $M$. As a result, we have
\fm{\sum_{\lambda:\lambda^M+1\gtrsim t}\lambda^{-\delta}\lesssim \frac{(Ct^{1/M})^{-\delta}}{1-(1+\mu)^{-\delta}}\lesssim_{\mu,\delta} t^{-\delta/M}.}

Next, we compute the sum of $(\lra{D}^lu_\lambda)(t,x)$ over $\lambda\in\mcl{I}$ satisfying $\lambda^M+1\lesssim t$. 
We recall \cite[Proposition 2.4]{MR4693097}.
\prop{\label{propb.1:psidecay}Fix $\lambda$ and fix $\xi_0\in\R$ with $\lra{\xi_0}\sim\lambda$. Set $x_0=t\xi_0/\lra{\xi_0}$.  Suppose that $\psi=\psi(t,x)$ is a function with $\supp\wh{\psi}(t,\cdot)\subset (-\infty,\xi_0]\cap\{\lra{\xi}\sim\lambda\}$. Also fix a constant $0<\kappa\leq 1/8$. Then, for $t\gtrsim_\kappa 1+\lambda^{2/\kappa}$ and $x\geq x_0-t^{1-\kappa}$, we have 
\eq{|\psi(t,x)|\lesssim_\kappa\lambda^{C_\kappa}t^{-(1+\kappa)/2}(\norm{L\psi}_{L^2}+\norm{\psi}_{L^2}).}
The same estimate holds if $\supp\wh{\psi}(t,\cdot)\subset[\xi_0,\infty)\cap\{\lra{\xi}\sim\lambda\}$ and for $x\leq x_0+t^{1-\kappa}$.}\rm
\bigskip

By setting $\psi=\lra{D}^{l}u_\lambda^\pm$ where $u_\lambda^\pm$ denotes the localization of $u_\lambda$ to positive or negative frequency, we apply Proposition \ref{propb.1:psidecay} with $\kappa=1/8$ to conclude that \eq{\norm{\lra{D}^{l}u_\lambda^\pm}_{L^\infty(\R\setminus (t\cdot \wt{J}_\lambda^\pm))}&\lesssim \eps \lambda^{C}t^{-9/16+C_A\eps^2}\lesssim \eps \lambda^{-2}t^{-9/16+C_A\eps^2+(C+2)/M}\\
&\lesssim \eps \lambda^{-2}t^{-17/32},\qquad \forall t\gtrsim \lambda^M+1.}
Here \fm{\wt{J}^\pm_\lambda=(J_\lambda\cup J_{\lambda/(1+\mu)}\cup J_{\lambda(1+\mu)})\cap\{\pm v>0\}.}
To get the last estimate, we choose $M\gg1$ and $\eps\ll_A1$ so that $C_A\eps^2+(C+2)/M\leq 1/32$. 
Fix $(t,x)\in [1,T]\times\R$. It is now clear that
\fm{\sum_{\lambda:\ \lambda^M+1\lesssim t,\ x/t\notin \wt{J}_\lambda}|\lra{D}^lu_\lambda(t,x)|\lesssim \sum_{\lambda}\eps\lambda^{-2}t^{-17/32}\lesssim \eps t^{-17/32}.}
If $(t,x/t)\notin\bigcup_{\lambda\in\mcl{I}}\{t\gtrsim \lambda^M+1,v\in J_\lambda\}$, one can remove the restriction $x/t\notin\wt{J}_\lambda$ in the sum on the left hand side. By combining this estimate with \eqref{est:b.13:D}, we obtain \eqref{sec:appB:keylemma:c2}. If there exists a $\lambda\in\mcl{I}$ such that $t\gtrsim\lambda^M+1$ and $x/t\in J_\lambda$, one can check that $x/t\notin\wt{J}_{\lambda'}$ unless $\lambda/\lambda'\in[(1+\mu)^{-1},1+\mu]$. As a result,  it remains to compute $\lra{D}^l\wt{u}_\lambda$ in $t\cdot J_\lambda$ where $\wt{u}_\lambda$ is defined by \eqref{wtphidefn}.

If  we set $\gamma^{l}(t,v)=\lra{u,\lra{D}^{l}\Psi_v}$, then by \cite[Lemma 5.10]{MR4693097}, we have
\fm{|\gamma^{l}(t,v)-t^{1/2}\lra{D}^{l}\wt{u}_\lambda(t,vt) e^{it\sqrt{1-v^2}}|\lesssim \eps t^{-1/4+C_A\eps^2}\lambda^C,\qquad v\in J_\lambda,\ t\gtrsim \lambda^M+1.}
By setting $\xi_v=v/\sqrt{1-v^2}$, for each $v\in J_\lambda$ and $ \lambda^M+1\lesssim t\leq T$, we also have
\fm{|\lra{\xi_v}^{l}\gamma-\gamma^{l}|&\lesssim \eps\lambda^Ct^{-\kappa_0},\qquad \text{where } \kappa_0\in(0,1)\text{ is a small constant.}}
This estimate holds because $\wh{\Psi}_v$ is essentially supported at $\xi_v$. Finally, since $t\gtrsim\lambda^M+1$, we have $\lambda^C\lesssim t^{C/M}$. By choosing $M\gg_{\kappa_0}1$ and $\eps\ll_{A,\kappa_0}1$, we have
\fm{t^{-1/4+C_A\eps^2}\lambda^C+t^{-\kappa_0}\lambda^C\lesssim t^{-\kappa},\qquad \kappa=\min\{\kappa_0/2,1/8\}.}
That is, for all $v\in J_\lambda$ and $\lambda^M+1\lesssim t\leq T$, we have
\eq{\label{est:approx:b.17}\lra{D}^{l}\wt{u}_\lambda(t,vt)&=t^{-1/2}e^{-it\sqrt{1-v^2}}\cdot\lra{\xi_v}^l\gamma(t,v)+O(\eps t^{-1/2-\kappa}).}
We thus obtain \eqref{sec:appB:keylemma:c}.

\subsubsection{Estimates for $\gamma$}\label{secb.2.3:estgamma}
Let us now prove \eqref{gammaclaim}. Fix $v\in J_\lambda$ and $\lambda^M+1\lesssim t\leq T$. We have
\fm{\dot{\gamma}&=\lra{(\partial_t+i\lra{D})u,\Psi_v}+\lra{u,(\partial_t+i\lra{D})\Psi_v}=\frac{i}{12}\lra{\xi_v}^{-1}\lra{(\Real u)^3,\Psi_v}+O(\eps t^{-5/4+C_A\eps^2}\lambda^C).}
Here we refer to  \cite[Lemma 4.4]{MR4693097} to control $(\partial_t+i\lra{D})\Psi_v$. Following the discussions in  \cite[Section 5]{MR4693097}, we conclude
\fm{\dot{\gamma}&=\frac{i}{12}\lra{\xi_v}^{-1}\lra{(\Real u_\lambda)^3,\Psi_v}+O_\kappa(\eps t^{-1-\kappa}\lambda^C),\qquad 0<\kappa\ll1.}
This estimate holds because $P_{\neq \lambda}\Psi_v$ decays sufficiently fast in time, and because $u_{\lambda'}$ has a good decay in the support of $\Psi_v$ if $v\notin J_{\lambda'}$. For simplicity, we skip the detailed proof here.

By  \cite[Lemma 5.10]{MR4693097}, we can replace $u_\lambda(t,x)$ with $t^{-1/2}\gamma(t,x/t) e^{-i\sqrt{t^2-x^2}}$. Besides, we can show that $|x/t-v|\lesssim\lambda^Ct^{-1/2}$ in $\supp\Psi_v$ and that
\fm{\norm{\partial_v\gamma(t)}_{L^2_v(J_\lambda)}\lesssim \lambda^C\norm{u}_{X}+\eps t^{-\kappa}\lambda^C\lesssim \eps\lambda^C t^{C_A\eps^2}.}
It follows that 
\fm{|\gamma(t,x/t)-\gamma(t,v)|&\leq|\int_{v}^{x/t}\partial_v\gamma(t,z)\ dz|\leq |x/t-v|^{1/2}\norm{\partial_v\gamma(t)}_{L^2_v(J_\lambda)}\lesssim\eps\lambda^C t^{-1/4+C_A\eps^2}. }
So we can also replace $\gamma(t,x/t)$ with $\gamma(t,v)$. We thus conclude that
\eq{\label{gammaeqn}\dot{\gamma}&= \frac{1}{96}i\lra{\xi_v}^{-1}t^{-1}(\gamma^3h_{-3}+\bar{\gamma}^3h_3+3|\gamma|^2\gamma +3|\gamma|^2\bar{\gamma}h_1)+O_\kappa(\eps t^{-1-\kappa}\lambda^C),\qquad 0<\kappa\ll1}
where \fm{h_k(t,v)&=t^{-1/2}\lra{e^{ik\sqrt{t^2-x^2}},\Psi_v}=\int_\R t^{-1/2}\lra{\xi_v}^{3/2}\chi(y)e^{i(1+k)\sqrt{t^2-x^2}}\ dx,\ y=t^{-1/2}\lra{\xi_v}^{3/2}(x-vt).}
The term $|\gamma|^2\gamma$ should be multiplied by $h_{-1}$, but we have  $h_{-1}=1$ (recall that $\int \chi=1$). 

Here we follow \cite{MR2457221}  and call the term $|\gamma|^2\gamma$ the \emph{resonant} cubic term.  We  call other cubic terms the \emph{nonresonant} cubic terms. In \cite{MR4693097}, only the resonant cubic term is involved. To handle the nonresonant cubic terms, we define 
\eq{R_{t,v}(y)&=t^{-1/2}\lra{\xi_v}^{-9/2}y^3\int_0^1\frac{3(t^{-1/2}\lra{\xi_v}^{-3/2}yh+v)(1-h)^2}{2(1-(t^{-1/2}\lra{\xi_v}^{-3/2}yh+v)^2)^{5/2}} dh.}
Then,  we have
\fm{-\sqrt{t^2-x^2}&=-t\lra{\xi_v}^{-1}+\xi_v(x-vt)+\frac{y^2}{2}+R_{t,v}(y).}
It follows that for $k\neq -1$, 
\fm{h_k&=t^{-1/2}\lra{e^{ik\sqrt{t^2-x^2}},\Psi_v}=t^{-1/2}\int \lra{\xi_v}^{3/2}\chi(y) e^{i(1+k)\sqrt{t^2-x^2}}\ dx\\
&=t^{-1/2}\int\lra{\xi_v}^{3/2}\chi(y)e^{-i(1+k)t\lra{\xi_v}^{-1}+i(1+k)\xi_v(x-vt)+i(1+k)\frac{y^2}{2}+i(1+k)R_{t,v}(y)}\ dx\\
&=e^{-i(1+k)t\lra{\xi_v}^{-1}}\int\chi(y)e^{i(1+k)\xi_vt^{1/2}\lra{\xi_v}^{-3/2}y+i(1+k)\frac{y^2}{2}+i(1+k)R_{t,v}(y)}\ dy.}
It follows that $h_k=O(1)$ and $\partial_t(h_ke^{i(1+k)t\lra{\xi_v}^{-1}})=O(t^{-1/2})$. The second estimate here gives us a better bound when we integrate $h_k$ with respect to the time. For any times $t_1\leq t_2$ such that $t_j\gtrsim \lambda^M+1$, by integrating by parts we have \fm{\int_{t_1}^{t_2}h_k(t,v)\ dt&=\int_{t_1}^{t_2}e^{-i(1+k)t\lra{\xi_v}^{-1}} \cdot e^{i(1+k)t\lra{\xi_v}^{-1}}h_k(t,v)\ dt\\
&=-\frac{\lra{\xi_v}(h_k(t_2)-h_k(t_1))}{i(1+k)}-\int_{t_1}^{t_2}-\frac{\lra{\xi_v}e^{-i(1+k)t\lra{\xi_v}^{-1}}}{i(1+k)}\partial_t(e^{i(1+k)t\lra{\xi_v}^{-1}}h_k)\ dt\\
&=O(\lambda+\int_{t_1}^{t_2}\lambda t^{-1/2}\ dt)=O(\lambda\sqrt{t_2}).}

We now finish the proof of \eqref{gammaclaim} by  a continuity argument. At $t_1\sim\lambda^M+1$, the estimate \eqref{gammaclaim} holds. Now assume that for some $T\geq t_1$ we have
\eq{\sup_{t\in[t_1,T]}|\gamma(t,v)|\leq B\eps\lambda^{-1-2\delta}.}
Here $B\geq 1$ is a large constant to be chosen. The equation \eqref{gammaeqn} implies
\eq{\label{gammaeqn:norm}\partial_t|\gamma|^2&=2\Real(\bar{\gamma}\dot{\gamma})=\Real\kh{\frac{1}{48}i\lra{\xi_v}^{-1}t^{-1}\bar{\gamma}(\gamma^3h_{-3}+\bar{\gamma}^3h_3+3|\gamma|^2\bar{\gamma}h_1)}+O_\kappa(\eps t^{-1-\kappa}\lambda^C|\gamma|).}
For  $t_2\geq t_1\sim\lambda^M+1$, by choosing $M\gg1$ we have
\fm{\int_{t_1}^{t_2}\eps t^{-1-\kappa}\lambda^C|\gamma|\ dt\lesssim B\eps^2\lambda^Ct_1^{-\kappa}\lesssim B\eps^2\lambda^{-2-4\delta}.}
 For each $\alpha=0,1,3$ with $\lambda^M+1\lesssim t_1\leq t_2$,  by setting $H_k(t,v)=\int_{t_1}^th_k(s,v)\ ds$, we have
\eq{\label{est:gammabar4}\int_{t_1}^{t_2}t^{-1}\bar{\gamma}^{4-\alpha}\gamma^\alpha h_{3-2\alpha}\ dt&=(t^{-1}\bar{\gamma}^{4-\alpha}\gamma^\alpha H_{3-2\alpha})|_{t=t_2}-\int_{t_1}^{t_2} H_{3-2\alpha}\partial_t(t^{-1}\bar{\gamma}^{4-\alpha}\gamma^\alpha )\ dt\\
&=O(\lambda|\gamma(t_2)|^4 t_2^{-1/2})+\int_{t_1}^{t_2}O(\lambda(t^{-3/2}|\gamma|^4+t^{-1/2}|\dot{\gamma}||\gamma|^3))\ dt\\&=O(B^4\eps^4\lambda^{-3-8\delta}t_1^{-1/2}+B^3\eps^4\lambda^Ct_1^{-1/2}).}
Since $t_1\sim\lambda^M+1$, by choosing $M\gg1$ and $\eps\ll1$, we can replace the right hand side by $O(\eps^2\lambda^{-2-4\delta} )$ with no $B$ involved.
Thus, whenever $v\in J_\lambda$ and $t_2\gtrsim\lambda^M+1$, we have
\fm{|\gamma(t_2)|^2&\lesssim \eps^2\lambda^{-2-4\delta}+\int_{t_1}^{t_2} \eps t^{-1-\kappa}\lambda^C|\gamma|\ dt+\abs{\int_{t_1}^{t_2}\lra{\xi_v}^{-1}t^{-1}\bar{\gamma}(\gamma^3h_{-3}+\bar{\gamma}^3h_3+3|\gamma|^2\bar{\gamma}h_1)\ dt}\\
&\lesssim (B+1)\eps^2\lambda^{-2-4\delta}.}
That is,
\fm{\sup_{t\in[t_1,T]}|\gamma(t,v)|\leq \sqrt{C(B+1)}\eps\lambda^{-1-2\delta}.}
By choosing $B\gg1$ so that $\sqrt{C(B+1)}\leq B/2$, we finish the proof.

\subsection{Asymptotics}\label{secb.3:asy}

We now prove the asymptotics for $\lra{D}^lu$ for $l\in[0,1+\delta]$. We shall follow the argument in \cite{MR3382579}.

We start with the following lemma on the asymptotic behavior of $\gamma(t,v)=\lra{u,\Psi_v}$ as $t\to\infty$. Since we have closed the continuity argument, we now have
\eq{\label{secb.3:gammabound}|\gamma(t,v)|\lesssim\eps \lambda^{-1-2\delta},\qquad \forall v\in J_\lambda,\ t\gtrsim \lambda^M+1.}
Besides, all the estimates in Section \ref{secb.2.3:estgamma} hold for all $t\gtrsim \lambda^M+1$.

\lem{\label{lemb.2:gammatvasy}There exists a complex-valued function $W=W(\xi)$ defined on $\R$, a sufficiently large constant $M>1$, and a sufficiently small constant $\kappa\in(0,1)$,  such that for each $\lambda\in\mcl{I}$, we have \eq{\label{lemb.2:gammatvasy:c}\gamma(t,v)&=W(\xi_v)\exp(\frac{1}{32}i\lra{\xi_v}^{-1}|W(\xi_v)|^2\ln t)+O_\kappa(\eps t^{-\kappa}),\qquad \forall t\gtrsim\lambda^M+1,\ v\in J_\lambda.}
We recall that $\xi_v=v/\sqrt{1-v^2}$. Moreover, we have $|W(\xi)|\lesssim\eps\lra{\xi}^{-1-2\delta}$.}
\begin{proof}
We first show that $\lim_{t\to\infty}|\gamma(t,v)|$ exists. Fix $\lambda\in\mcl{I}$ and $v\in J_\lambda$. Set
\fm{\mcl{R}&=\Real\kh{\frac{1}{48}i\lra{\xi_v}^{-1}t^{-1}\bar{\gamma}(\gamma^3h_{-3}+\bar{\gamma}^3h_3+3|\gamma|^2\bar{\gamma}h_1)}.}
By \eqref{gammaeqn:norm} and \eqref{est:gammabar4}, for all $t_2\geq t_1\gtrsim\lambda^M+1$, we have
\fm{||\gamma(t_2,v)|^2-|\gamma(t_1,v)|^2|&\lesssim_\kappa|\int_{t_1}^{t_2}\mcl{R}\ dt|+\int_{t_1}^{t_2}\eps t^{-1-\kappa}\lambda^C\ dt\\
&\lesssim \eps^4\lambda^{C}t_1^{-1/2}+\eps^2\lambda^Ct_1^{-\kappa}\lesssim \eps^2\lambda^Ct_1^{-\kappa}.}
for some small constant $\kappa\in(0,1)$. Then, for each $v\in J_\lambda$ (or equivalently $\xi_v\in I_\lambda$), the limit
\fm{K(\xi_v):=\lim_{t\to\infty}|\gamma(t,v)|} exists. By \eqref{secb.3:gammabound}, we have $|K(\xi_v)|\lesssim \eps\lambda^{-1-2\delta}\sim\eps\lra{\xi_v}^{-1-2\delta}$ for all $\xi_v\in\R$. We also have
\fm{||\gamma(t,v)|^2-K(\xi_v)^2|&\lesssim_\kappa \eps^2\lambda^Ct^{-\kappa},\qquad \forall t\gtrsim\lambda^M+1,\ v\in J_\lambda.}

Next, we return to \eqref{gammaeqn}. Set \fm{\wt{\mcl{R}}=\frac{1}{96}i\lra{\xi_v}^{-1}t^{-1}(\gamma^3h_{-3}+\bar{\gamma}^3h_3+3|\gamma|^2\bar{\gamma}h_1).}
It follows from \eqref{gammaeqn} that
\fm{&\partial_t(\gamma\exp(-\frac{1}{32}i\lra{\xi_v}^{-1}K(\xi_v)^2\ln t))\\
&=\exp(-\frac{1}{32}i\lra{\xi_v}^{-1}K(\xi_v)^2\ln t)\kh{\frac{i}{32}\lra{\xi_v}^{-1}t^{-1}(|\gamma|^2-K(\xi_v)^2)\gamma+\wt{\mcl{R}}+O_\kappa(\eps\lambda^Ct^{-1-\kappa})}\\
&=\exp(-\frac{1}{32}i\lra{\xi_v}^{-1}K(\xi_v)^2\ln t)\wt{\mcl{R}}+O_\kappa(\eps\lambda^Ct^{-1-\kappa}).}
Following the proof of \eqref{est:gammabar4}, we also have
\fm{|\int_{t_1}^{t_2}\exp(-\frac{1}{32}i\lra{\xi_v}^{-1}K(\xi_v)^2\ln t)\wt{\mcl{R}}\ dt|\lesssim\eps^3\lambda^Ct_1^{-1/2},\qquad \forall t_2\geq t_1\gtrsim\lambda^M+1.}
As a result, for all $t_2\geq t_1\gtrsim \lambda^M+1$, we have
\fm{&|\gamma(t_2,v)\exp(-\frac{1}{32}i\lra{\xi_v}^{-1}K(\xi_v)^2\ln t_2)-\gamma(t_1,v)\exp(-\frac{1}{32}i\lra{\xi_v}^{-1}K(\xi_v)^2\ln t_1)|\lesssim_\kappa\eps\lambda^C t_1^{-\kappa}.}
As a result, for each $v\in J_\lambda$, the limit
\fm{W(\xi_v)&=\lim_{t\to\infty}\gamma(t,v)\exp(-\frac{1}{32}i\lra{\xi_v}^{-1}K(\xi_v)^2\ln t)}
exists. Moreover, we have $K(\xi_v)=|W(\xi_v)|$ and
\eq{\label{lem:gammatvasy:est}|\gamma(t,v)-W(\xi_v)\exp(\frac{1}{32}i\lra{\xi_v}^{-1}|W(\xi_v)|^2\ln t)|\lesssim_\kappa\eps\lambda^C t^{-\kappa},\qquad \forall t\gtrsim\lambda^M+1,\ v\in J_\lambda.}

Finally, we notice that the constant $C$ in $\lambda^C$ is independent of the choice of $M$. Thus, by replacing $M$ with $M+2C/\kappa$ if necessary, for all $t\gtrsim\lambda^M+1$, we have $\lambda^C t^{-\kappa/2}\lesssim \lambda^{C-M\kappa/2}\lesssim 1$. In other words, we can remove the factor $\lambda^C$ in \eqref{lem:gammatvasy:est} by replacing $t^{-\kappa}$ with $t^{-\kappa/2}$.
\end{proof}
\rmk{\rm The function $W$ obtained is the same as that in \cite[Theorem 6]{MR4693097}. We thus have $W\in H^{1-C\eps}_{\rm loc}(\R)$ and bounds such as \fm{\norm{\lra{\xi}^{m-3/2-C\eps^2}W}_{L^2}+\norm{\lra{\xi}^{m-3/2-C\eps^2}W}_{H^{1-C\eps}}\lesssim \eps.}
For simplicity, we skip their proofs and refer to \cite{MR4693097}.}\rm
\bigskip

We now apply Lemma \ref{sec:appB:keylemma}. Choose a pair of constants $(M,\kappa)$ so that both Lemma \ref{sec:appB:keylemma} and Lemma \ref{lemb.2:gammatvasy} hold. By \eqref{sec:appB:keylemma:c} and \eqref{lemb.2:gammatvasy:c}, for all $l\in[0,1+\delta]$, $\lambda\in\mcl{I}$, $t\gtrsim \lambda^M+1$ and $v\in J_\lambda$, we have
\fm{\lra{D}^lu(t,vt)&=t^{-1/2}e^{-it\sqrt{1-v^2}}\lra{\xi_v}^{l}\gamma(t,v)+O_\kappa(\eps t^{-1/2-\kappa})\\
&=t^{-1/2}e^{-it\sqrt{1-v^2}}\lra{\xi_v}^{l}W(\xi_v)\exp(\frac{i}{32}\lra{\xi_v}^{-1}|W(\xi_v)|^2\ln t)+O_\kappa(\eps t^{-1/2-\kappa}(\lambda^{1+\delta}+1)).}
By replacing $(M,\kappa)$ with $(M+
(2+2\delta)/\kappa,\kappa/2)$, we can remove the factor $\lambda^{1+\delta}$ in the last estimate. This is because $t^{-\kappa/2}\lambda^{1+\delta}\lesssim \lambda^{1+\delta-\kappa M/2}\lesssim 1$ for all $t\gtrsim \lambda^M+1$.

Now suppose $(t,v)\in[1,\infty)\times\R$ but $(t,v)\notin\bigcup_{\lambda\in\mcl{I}}\{t\gtrsim\lambda^M+1, v\in J_\lambda\}$. In this case, we either have $|v|\geq 1$, or $v\in J_\lambda$ but $t\lesssim \lambda^M+1$. In the first case, we apply \eqref{sec:appB:keylemma:c2} to get $\lra{D}^lu(t,vt)=O_\kappa(\eps t^{-1/2-\kappa})$. In the second case, we notice that $\lra{D}^lu(t,vt)=O_\kappa(\eps t^{-1/2-\kappa})$ and that
\fm{|t^{-1/2}e^{-it\sqrt{1-v^2}}\lra{\xi_v}^{l}W(\xi_v)\exp(\frac{i}{32}\lra{\xi_v}^{-1}|W(\xi_v)|^2\ln t)|&\lesssim t^{-1/2}\lra{\xi_v}^{1+\delta}|W(\xi_v)|\lesssim \eps t^{-1/2}\lambda^{-\delta}\\
&\lesssim \eps t^{-1/2-\delta/M}.}
In the second last step, we use the pointwise estimate for $W$ in Lemma \ref{lemb.2:gammatvasy}. In the last step, we recall that $\lambda\gtrsim t^{1/M}$ with an implicit constant independent of $M$ whenever $t\lesssim \lambda^M+1$. In summary, by shrinking $\kappa$ if necessary, we conclude that for all $l\in[0,1+\delta]$, $t\geq 1$, and $x\in\R$, 
\eq{&(\lra{D}^lu)(t,x)\\&=t^{-1/2}\lra{x/\rho}^{l}W(x/\rho)\exp(-i\rho+\frac{i}{32\lra{x/\rho}}|W(x/\rho)|^2\ln t)\cdot 1_{|x|<t}+O_\kappa(\eps t^{-1/2-\kappa}).}
Here $\rho=\sqrt{t^2-x^2}$. This finishes the proof of Proposition \ref{thmwpt}.

\subsection{A final remark}\label{lastappendixb}
Let us now prove \eqref{zuptbd}. We still consider the Cauchy problem \eqref{sg1} but now we assume that $\norm{u_0}_{H^{m,2}}\lesssim\eps\ll1$ with $m\in(5/2,3]$. We seek to prove that for each constant $\lambda\in(0,1/2)$, as long as $\eps\ll_\lambda1$,
\fm{\norm{Zu(t)}_{L^\infty}\lesssim \eps\lra{t}^{-\lambda},\qquad \forall t\geq 0.}
Recall that $Z=t\partial_x+x\partial_t$ is the Lorentz boost.

We first summarize some estimates from Section \ref{wtpenergyest}. 
\lem{For each $t\geq 0$, we have
\fm{\norm{u(t)}_{H^m}+\norm{u_t}_{H^{m-1}}+\norm{Zu(t)}_{H^{m-1}}+\lra{t}^{-1}\norm{xu(t)}_{H^m}&\lesssim \eps \lra{t}^{C\eps^2},\\
\lra{t}\norm{(\partial_t+i\lra{D})u(t)}_{H^m}+\lra{t}\norm{(\partial_t+i\lra{D})(Zu)(t)}_{H^{m-1}}+\eps^2\norm{(\partial_t+i\lra{D})(xu)(t)}_{H^{m-1}}&\lesssim \eps^3 \lra{t}^{C\eps^2},\\
\norm{\lra{D}^{1+\delta}u(t)}_{L^\infty}&\lesssim\eps\lra{t}^{-1/2},\\
\lra{t}^{-1}\norm{xZu(t)}_{H^{m-2}}+\norm{(\partial_t+i\lra{D})(xZu)(t)}_{H^{m-2}}&\lesssim\eps\lra{t}^{C\eps^2}.}
Here $\delta>0$ is a small constant.
}
\begin{proof}
The only estimate which requires a proof is the last one. Using the commutators and bounds in Section \ref{wtpenergyest}, we have
\fm{&\norm{(\partial_t+i\lra{D})(xZu)(t)}_{H^{m-2}}\lesssim \norm{x(\partial_t+i\lra{D})(Zu)(t)}_{H^{m-2}}+\norm{Zu(t)}_{H^{m-2}}\\
&\lesssim \norm{xZ(\partial_t+i\lra{D})u(t)}_{H^{m-2}}+\norm{x(\partial_t+i\lra{D})u(t)}_{H^{m-2}}+\eps\lra{t}^{C\eps^2}\\
&\lesssim \norm{\lra{D}(xZ\lra{D}^{-1}(\Real u-\sin(\Real u))(t)}_{L^2}+\eps\lra{t}^{C\eps^2}\\
&\lesssim \norm{xZ(\Real u-\sin(\Real u))(t)}_{L^2}+ \norm{Z(\Real u-\sin(\Real u))(t)}_{L^2}\\
&\quad+ \norm{x\partial_t(\Real u-\sin(\Real u))(t)}_{L^2}+ \norm{\partial_t(\Real u-\sin(\Real u))(t)}_{L^2}+\eps\lra{t}^{C\eps^2}.
}
To get the last bound, we keep applying the formulas of commutators such as $[\lra{D},x]$, $[\lra{D},Z]$, etc. By the estimates proved in Section \ref{wtpenergyest}, we have
\fm{\abs{\frac{d}{dt}\norm{xZu(t)}_{H^{m-2}}}\lesssim \norm{(\partial_t+i\lra{D})(xZu)(t)}_{H^{m-2}}\lesssim\eps^2\lra{t}^{-1}\norm{xZu(t)}_{L^2}+\eps \lra{t}^{C\eps^2}.}
Since $\norm{xZu(0)}_{H^{m-2}}\lesssim \norm{u_0}_{H^{m,2}}\lesssim\eps $, we conclude that $\norm{xZu(t)}_{H^{m-2}}\lesssim \eps \lra{t}^{1+C\eps^2}$ by  Gronwall's inequality.
\end{proof}
\rm

Since $\norm{Zu(t)}_{H^2}\lesssim \eps \lra{t}^{C\eps^2}$, there is nothing to prove when $t\lesssim1$. From now on we assume $t\gtrsim1$. For simplicity, we now set
\fm{Q_1(t)&=\norm{Zu(t)}_{L^\infty},\qquad Q_2(t)=\norm{Z^2u(t)}_{H^{m-2}}.}
By $Z^2$, we mean applying the Lorentz boost twice. Since $m>5/2$, by \eqref{ptwbdlem3.2} we have
\eq{\label{q1est}\lra{t}^{1/2}Q_1(t)&\lesssim\norm{Zu(t)}_{H^{3/2}}^{1/2}\cdot (\norm{Zu(t)}_{H^{3/2}}^{1/2}+\norm{LZu(t)}_{H^{3/2}})^{1/2}\\
&\lesssim \eps \lra{t}^{C\eps^2}+\eps^{1/2}\lra{t}^{C\eps^2}(\norm{Z^2u(t)}_{H^{1/2}}+\norm{(\partial_t+i\lra{D})(xZu)(t)}_{H^{1/2}})^{1/2}\\
&\lesssim \eps \lra{t}^{C\eps^2}+\eps^{1/2}\lra{t}^{C\eps^2}Q_2(t)^{1/2}.}
In the second identity, we use $\lra{D}^{-1}Z=\lra{D}^{-1}(\partial_t+i\lra{D})x-iL$.

Moreover, to estimate $Q_2(t)$, we notice that $[\partial_t+i\lra{D},Z]=D\lra{D}^{-1}(\partial_t+i\lra{D})$, we have
\fm{&\abs{\frac{d}{dt}Q_2(t)}\lesssim\norm{(\partial_t+i\lra{D})(Z^2u)}_{H^{m-2}}\\&\lesssim\norm{Z^2(\partial_t+i\lra{D})u}_{H^{m-2}}+\norm{ZD\lra{D}^{-1}(\partial_t+i\lra{D})u}_{H^{m-2}}+\norm{(\partial_t+i\lra{D})(Zu)}_{H^{m-2}}\\
&\lesssim \norm{Z^2\lra{D}^{-1}(\Real u-\sin(\Real u))}_{H^{m-2}}+\norm{ZD\lra{D}^{-2}(\Real u-\sin(\Real u))}_{H^{m-2}}+\eps^3 \lra{t}^{-1+C\eps^2}\\
&\lesssim \norm{Z^2(\Real u-\sin(\Real u))}_{L^2}+\eps^3\lra{t}^{-1+C\eps^2}.
}
The proof is straightforward as we keep applying the commutators. Now using  Leibniz's rule, we have
\fm{\norm{Z^2(\Real u-\sin(\Real u))}_{L^2}&\lesssim
\norm{u}^2_{L^\infty}\norm{Z^2u}_{L^2}+\norm{u}_{L^\infty}\norm{Zu}_{L^2}\norm{Zu}_{L^\infty}\\&\lesssim \eps^2\lra{t}^{-1}\norm{Z^2u}_{L^2}+\eps^2\lra{t}^{-1/2+C\eps^2}Q_1(t).}
Since $Q_2(0)\lesssim\norm{u_0}_{H^{m,2}}\lesssim\eps$, by Gronwall's inequality, we have
\eq{\label{q2est}Q_2(t)\lesssim \eps\lra{t}^{C\eps^2}+\lra{t}^{C\eps^2}\int_0^t\eps^2\lra{s}^{-1/2+C\eps^2}Q_1(s)\ ds.}
It follows from \eqref{q1est} and \eqref{q2est} that
\fm{\lra{t}Q_1(t)^2&\lesssim \eps^2\lra{t}^{C\eps^2}+\eps^2\lra{t}^{C\eps^2}\int_0^t\eps \lra{s}^{-1/2+C\eps^2}Q_1(s)\ ds.}
Set $Q_3(t)=\lra{t}^{1-\kappa}Q_1(t)^2$ where $\kappa\in(0,1/10)$ is a small constant to be chosen. By choosing $\eps\ll_\kappa1$, we have $C\eps^2<\kappa$, so it follows that
\fm{Q_3(t)&\lesssim \eps^2+\eps^2\int_0^t \eps\lra{s}^{-1+\kappa/2+C\eps^2}\sqrt{Q_3(s)}\ ds
\\&\lesssim \eps^2+\eps^2\int_0^t \eps\lra{s}^{-1+\kappa/2+C\eps^2}(\lra{s}^{-\kappa}Q_3(s)+\lra{s}^{\kappa})\ ds\\
&\lesssim \eps^2(1+\eps\lra{t}^{3\kappa/2+C\eps^2})+\eps^2\int_0^t\eps\lra{s}^{-1-\kappa/2+C\eps^2}Q_3(s) \ ds.}
In the first estimate, we use the inequality $2\sqrt{AB}\leq A+B$. Since $Q_3(0)\lesssim \eps^2$, we conclude by  Gronwall's inequality that 
\fm{Q_3(t)&\lesssim  \eps^2\lra{t}^{3\kappa/2+C\eps^2}\Longrightarrow Q_1(t)\lesssim \eps\lra{t}^{-1/2+2\kappa+C\eps^2}.}Here we use $-1-\kappa/2+C\eps^2<-1-\kappa/4$ for $\eps\ll_\kappa$1.
This finishes our proof.

\bibliography{paper_sg}{}
\bibliographystyle{plain}

\end{document}